\documentclass[12pt,a4paper,twoside]{report}

\usepackage{epsf,graphicx}
\usepackage{amsmath,latexsym,amssymb,amsthm}
\usepackage{hyperref}
\usepackage{setspace,cite}
\usepackage{lscape,fancyhdr,fancybox}
\usepackage{a4}
\usepackage{dsfont,texdraw}
\usepackage[hmarginratio=1:1, vmarginratio =6:5,
textheight=23cm,bindingoffset=1.6cm, textwidth=15.3cm]{geometry}
\numberwithin{equation}{section}

\theoremstyle{plain}
\newtheorem{thm}{Theorem}[section]
\newtheorem{cor}[thm]{Corollary}
\newtheorem{lem}[thm]{Lemma}
\newtheorem{prop}[thm]{Proposition}
\theoremstyle{definition}
\newtheorem{defn}{Definition}[section]
\theoremstyle{remark}
\newtheorem{rem}{Remark}[section]
\newtheorem{ex}{Example}[section]

\newcommand{\seq}[1]{\left<#1\right>}
\newcommand{\si}{\sigma}
\newcommand{\Si}{\Sigma}
\newcommand{\mbf}{\mathbf}
\newcommand{\mbb}{\mathbb}
\newcommand{\dsty}{\displaystyle}
\newcommand{\dsum}{\sum\limits}
\newcommand{\dlim}{\lim\limits}
\newcommand{\dsup}{\sup\limits}
\newcommand{\dinf}{\inf\limits}
\newcommand{\dmin}{\min\limits}
\newcommand{\dmax}{\max\limits}
\newcommand{\dprod}{\prod\limits}
\newcommand{\ra}{\rightarrow}
\newcommand{\Ra}{\Rightarrow}
\newcommand{\ol}{\overline}
\newcommand{\ul}{\underline}
\newcommand{\al}{\alpha}
\newcommand{\be}{\beta}
\newcommand{\de}{\delta}
\newcommand{\ep}{\epsilon}
\newcommand{\ga}{\gamma}
\newcommand{\la}{\lambda}
\newcommand{\om}{\omega}
\newcommand{\tri}{\triangle}
\newcommand{\itri}{\triangledown}

\begin{document}
\pagestyle{empty}
\vspace*{1cm}
\begin{center}
\begin{huge}
\textbf{Contributions to Random Energy Models\\}
\end{huge}
\end{center}

\vspace{2cm}

\begin{center}
\begin{large}
\textbf{Nabin Kumar Jana\\}
\end{large}
\end{center}

\vspace{3cm}
\begin{figure}[h]
\begin{center}
\includegraphics[width=4cm]{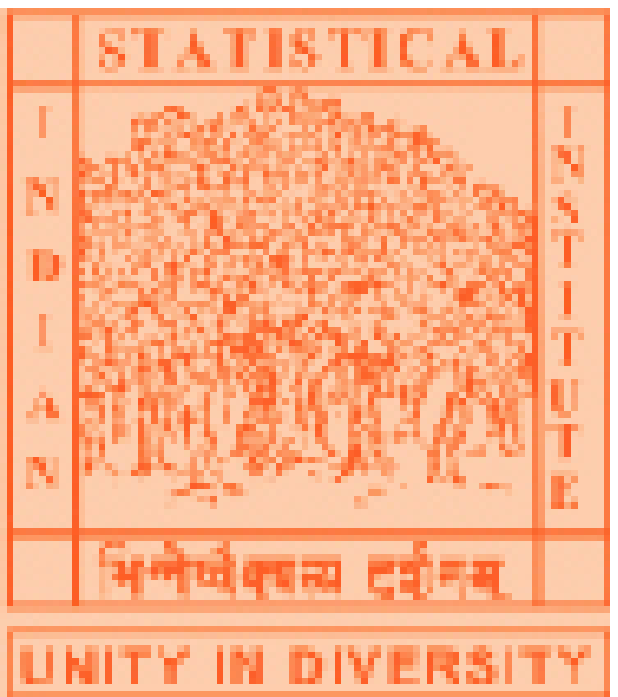}
\end{center}
\end{figure}
\vspace{2cm}

\begin{center}
\begin{large}
\textbf{Indian Statistical Institute\\
Kolkata\\
2007}
\end{large}
\end{center}

\cleardoublepage

\pagestyle{empty}
\vspace*{1cm}
\begin{center}
\begin{huge}
\textbf{Contributions to Random Energy Models\\}
\end{huge}
\end{center}

\vspace{3cm}

\begin{center}
\begin{large}
\textbf{Nabin Kumar Jana\\}
\end{large}
\end{center}

\vspace{4cm}
\begin{center}
Thesis submitted to the Indian Statistical Institute\\
in partial fulfillment of the requirements\\
for the award of
the degree of\\
Doctor of Philosophy.\\
October, 2007
\end{center}

\vspace{3cm}
\begin{center}
\begin{large}
\textbf{Indian Statistical Institute\\
203, B.T. Road, Kolkata, India.}
\end{large}
\end{center}

\cleardoublepage
\pagenumbering{roman}
\pagestyle{plain}
\pagestyle{empty} 
\vspace*{8cm}
\begin{center}
\em{To My Parents\\}
\end{center}
\cleardoublepage
\vspace*{1cm}
\begin{center}
\begin{Large}
\textbf{Acknowledgements}\\
\end{Large}
\end{center}
\vspace{2cm}
\def\baselinestretch{1.6}

A simple thanks will not be enough to convey my gratitude towards my
supervisor, Professor B. V. Rao for his careful guidance, constant
encouragements and the time he has given me during this project.

I am grateful to Professor Rahul Roy for introducing several models
in statistical physics and providing scope to visit ISI, Delhi
campus.

I would like to take this opportunity to thank all my teachers in
Stat-Math Unit of Kolkata, Delhi and Bangalore center for their
illuminating courses and fruitful discussions.

I am thankful to Dr. T. Mukherjee, Principal of Bijoy Krishna Girls
College and my departmental colleagues of this college for their
constant co-operations and encouragements to complete this work.

I thank all my friends, seniors, juniors as well as the members of
Stat-Math Unit for providing such a wonderful atmosphere to work on.

My thanks also goes to my friends Sourabh, Sohini, Soumenda, Dola
and dada Ac' Devatmananda Avt. for their constant encouragements.
Finally, I want to thank my wife who shouldered a lot of
responsibilities so that I could concentrate only to my research
work.

\cleardoublepage
\tableofcontents
\cleardoublepage

\addtolength{\headheight}{15pt}
\pagenumbering{arabic}
\setcounter{page}{1}
\pagestyle{myheadings}
\setcounter{chapter}{-1}
\pagestyle{fancy}
\renewcommand{\chaptermark}[1]{\markboth{\chaptername
\ \thechapter:\,\ #1}{}}
\renewcommand{\sectionmark}[1]{\markright{\ #1}}
\fancyhf{}
\fancyhead[LE]{\sl\leftmark}
\fancyhead[LO,RE]{\rm\thepage}
\fancyhead[RO]{\sl\rightmark}
\fancyfoot[C,L,E]{}
\chapter{Introduction}

In this introductory chapter, we begin with a brief description of
spin glasses in section 1. We are not physicists. The purpose of
this section is to trace the history of the models. Section 2 gives
a brief summary of the thesis and section 3 recalls certain known
facts which will be used later in the thesis.

\section{Origin of the problem}
The models considered in this thesis have their origin in spin glass
theory. Roughly, {\em spin glass} is a glassy state in a spin system
or a disordered material exhibiting high magnetic frustration. The
origin of this behavior can be either a disordered structure (such
as that of a conventional, chemical glass) or a disordered magnetic
doping in an otherwise regular structure. But what is a glass?
Loosely speaking, it is a state of spins with local ordering (in
solid state physics, this is called local `freezing' - locally, the
system looks more like an ordered solid rather than a disordered
liquid) but no global ordering. Spin glass can not remain in a
single lowest energy state (the ground state). Rather it has many
ground states which are never explored on experimental time scales.
The freezing of the spins, in spin glasses, is not a deterministic
one like ferromagnetic materials. Rather they freeze in random with
some memory effect.

Experiments show that the susceptibility obtained by cooling the
spin glass system in the presence of a magnetic field yielded a
higher value than that obtained by first cooling in zero field and
then applying the magnetic field. If the spin glass is cooled below
$T_c$ (a certain critical temperature) in the absence of an external
field, and then a magnetic field is applied, there is a rapid
increase towards a value, called the zero-field-cooled
magnetization. This value is less than the field-cooled
magnetization. The following phenomenon has also been observed in
the measurement of remanent magnetization (the permanent
magnetization that remains after the external field is removed). We
can cool in the presence of external field, remove the external
field and then measure the remanent magnetization. Alternatively,
first cool with out the external field, then apply the external
field and measure the remanent magnetization after removing the
external field. The first value is larger than the second one.

The other peculiarity of the spin glasses is its time dependence,
which will be explained now, that makes it different from other
magnetic systems. Above the spin glass transition temperature,
$T_c$, the spin glass exhibits typical magnetic behavior. In other
words, at temperature above $T_c$, if an external magnetic field is
applied and the magnetization is plotted versus temperature, it
follows the typical Curie law (in which magnetization is inversely
proportional to temperature). This happens until $T_c$ is reached,
at which point the magnetization becomes virtually constant. This is
the onset of the spin glass phase. When the external field is
removed, the spin glass has a rapid decrease of magnetization to a
value called the remnant magnetization, and then a slow decay as the
magnetization approaches zero (or some small fraction of the
original value). This decay is non-exponential and no single
function can fit the curve of magnetization versus time adequately
below $T_c$. This slow decay is particular to spin glasses.
If a similar procedure was followed for a ferromagnetic substance,
when the external field is removed, there would be a rapid change to
a remnant value, but this value is a constant in time. For a
paramagnetic material, when the external field is removed, the
magnetization rapidly goes to zero. In each case, the change is very
rapid and if carefully examined it is exponential decay.

Behind this strange behaviour of spin glasses, according to
physicists, there are essentially two major causes. These are {\em
quenched disorder} and {\em frustration}. The term ``quenched disorder''
refers to constrained disorder in the interactions between the spins
and/or their locations but does not evolve with time. In statistical
physics, a system is said to present quenched disorder when some
parameters defining its behaviour are random variables which do not
evolve with time, i.e., they are quenched or frozen. This is in
contrast to annealed disorder, where the random variables are
allowed to evolve themselves. Usually the spin orientations depend
on several facts such as the interactions, external fields and
thermal fluctuations. Their dynamics or thermodynamics will suggest
whether to order or not. The spin glass phase is an example of
spontaneous cooperative freezing (or order) of the spin orientations
in the presence of the constrained disorder of the interactions or
spin locations. It is thus ``order in the presence of disorder''. On
the other hand, ``frustration'' refers to conflicts between
interactions and the spin-ordering forces, and not all can be obeyed
simultaneously. Frustration arises when pairs of spins get different
ordering instructions through the various paths which link them,
either directly or via intermediate spins. The relevance of
frustration is that it leads to degeneracy or multiplicity of
compromises forcing the system to have several ground states.

Keeping these two in mind, in 1975, S. F. Edwards and P. W. Anderson~\cite{EA}
produced a paper, which in the words of Sherrington~\cite{Sh06}, at one fell swoop recognized the importance of
the combination of frustration and quenched disorder as fundamental
ingredients, introduced a more convenient model, a new and novel
method of analysis, new types of order parameters, a new mean field
theory, new approximation techniques and the prediction of a new
type of phase transition apparently explaining the observed
susceptibility cusp. This paper was a watershed. Edwards and
Anderson's new approach was beautifully minimal, fascinating and
attractive but also their analysis was highly novel and
sophisticated, involving radically new concepts and methods but also
unusual and unproven ans\"{a}tze, as well as several different
approaches. In their model, two spins interact if they are neighbour
to each other. The same year Sherrington and Kirkpatrick~\cite{SK} proposed
their model with mean field interaction. In this model all spins
interact with each other. In both the cases the interaction among
the spins were random and driven by Gaussian random variables. Due
to rich and complicated correlation structure among the energy over
the configuration space of the spins, initially the models were not
easy to study analytically. To get some insight into these models,
in 1980, B. Derrida~\cite{De80} proposed a system without any correlation
structure over the configuration space. He proposed a solvable model
called Random Energy Model (REM) for spin glass theory. In REM, all
the random variables are independent and identically distributed but
the distribution depends on the number of particles. Like
Edwards-Anderson model and SK-model, he considered these random
variables to be Gaussian. But this is a toy model since the energy
of the system does not depend on the configuration. Amazingly he
could show that though this is a very simple model, it exhibits
phase transition.

REM has no correlations at all. But the correlation structure in the
Edwards-Anderson model and SK-model were very complicated. So the
next idea is to study a system which exhibits correlations, but
their structure is simple enough to explicitly solve the model. B.
Derrida~\cite{D2} proposed another model for spin glass theory in 1985, by
bringing correlations through a tree structure. The tree structure
comes from the configuration space. Simply put, he identifies the
configuration space as the branches of a tree. This is called
Generalized Random Energy Model (GREM), a generalization of the REM.
Here also the driving distributions were Gaussian. In this project
we will focus ourselves on REM and GREM and some related models.

\medskip
\section{Setup and Summary}

For an $N$ particle system with classical spins $+1$ or $-1$, a
sequence of +1 and -1 of length $N$ gives a configuration of the
system. A typical configuration is denoted by $\sigma(N)$ or by
$\si$ when $N$ is understood. That is, $\sigma$ is a sequence of +1
and -1 of length $N$. The space of all possible configurations
$\sigma$ of a system is called configuration space and denoted by
$\Sigma_N$ or simply by $2^N$ since $\Sigma_N$ is nothing but $\{+1,
-1\}^N$. Now depending on the configuration, the system possesses
some energy called Hamiltonian. For a configuration $\sigma$, it is
denoted by $H_N(\sigma)$. The model is defined through the
Hamiltonian. So different models have different Hamiltonian
structures. In spin glass theory, the Hamiltonian is considered to
be random.

When the system is cooled, it settles down at a configuration where
the Hamiltonian is minimized. Hence it is very essential to get
information about the configurations where the infimum of the the
Hamiltonian is attained and its value. In statistical physics one
analyzes this problem via the partition function of the system. The
partition function, denoted by $Z_N(\beta)$, is defined as follows:

$$Z_N(\beta)=\sum\limits_{\sigma\in \Sigma_N} e^{-\beta
H_N(\sigma)}.$$ Here $\beta\geq 0$ is a parameter, represents the
inverse temperature. Sometimes when the Hamiltonian depends on an
external field $h$, we will denote the partition function as
$Z_N(\beta,h)$. Now note that among all the summands in the above
sum if one takes large $\beta$, only that summand will contribute
where the Hamiltonian attains the minimum among all possible
configurations. On the other hand, if the focus is on maximum, then
instead of $-\beta$ one has to consider $\beta$ in the exponent.

But the information in partition function about the minimum energy
is in exponential scale. So it is customary to study the logarithm
of the partition function. Also the energy of the system depends on
the number of particles in the system and becomes large when $N$ is
large. To get some asymptotic result on $\log Z_N(\be)$, one has to
normalize it properly. In this case, $\frac{1}{N}$ is the correct
normalization (in some sense). According to statistical physics,
$-\frac{1}{\beta N} \log Z_N(\beta)$ is called the {\em free energy}
of the system. Since one is interested in the asymptotic of the free
energy, that is, in $-\frac{1}{\beta}\lim\limits_{N\ra \infty}
\frac{1}{N} \log Z_N(\beta)$, for mathematical purpose we can forget
about the $-\frac{1}{\beta}$ term in the definition of free energy.
And from now on, we will call $\lim\limits_{N\ra \infty} \frac{1}{N}
\log Z_N(\beta)$ as the {\bf free energy} of the system.

In statistical physics, there is another important concept called Gibbs'
distribution. This is a distribution on the configuration space.
According to this, the probability of a
configuration $\sigma$ is proportional to $e^{-\beta
H_N(\sigma)}$. In particular, if $G_N(\sigma)$ denotes the Gibbs'
probability for a configuration $\sigma\in \Sigma_N$, then
$$G_N(\sigma)=\frac{e^{-\beta
H_N(\sigma)}}{Z_N(\beta)}.$$ It is worth noting that, since
$H_N(\sigma)$'s are random, the Gibbs' distribution is also random.
Note that, Gibbs' distribution is so defined as to give maximum
weight to that configuration which has minimum energy. We shall not
deal with Gibbs' distributions in this thesis.

Generalized random energy model (GREM) is one model in this theory
proposed by B. Derrida~\cite{D2} in 1985. To describe a version of
this fix an integer $n\geq 1$. For $N$ particle system, consider a
partition of $N$ into integers $k(i,N)\geq 0$ for $1\leq i \leq n$
so that $\dsum_i k(i,N)=N$. The configuration space $2^N$, naturally
splits into the product, $\prod 2^{k(i,N)}$ and $\si \in 2^N$ can be
written as $\si_1\si_2\cdots\si_n$ with $\si_i\in 2^{k(i,N)}$. An
obvious $n$-level tree structure can be brought in the configuration
space. Consider an $n$ level tree with $2^{k(1,N)}$ many edges at
the first level. These edges are denoted by $\si_1$, with
$\si_1\in2^{k(1,N)}$. In general, below a typical edge
$\si_1\si_2\cdots\si_{i-1}$ of the $(i-1)$-th level there are
$2^{k(i,N)}$ edges at the $i$-th level denoted by
$\si_1\si_2\cdots\si_{i-1}\si_i$ for $\si_i\in 2^{k(i,N)}$. Thus a
typical branch of the tree reads like $\si_1\si_2\cdots\si_n$ making
a one one correspondence with $2^N$, the configuration space. For
each $i$, $1\leq i \leq n$ and edge $\si_1\cdots\si_i$, associate a
random variables $\xi(\si_1\cdots\si_i)$. All these random variables
are i.i.d. $\mathcal{N}(0,N)$. One non random weight, $a_i>0$ for
each level is fixed. In GREM, Hamiltonian for a configuration
$\si=\si_1\cdots\si_n$ is defined as

\begin{equation}
H_N(\si) = \dsum_{i=1}^{n} a_i\xi(\si_1\si_2\cdots\si_i).
\end{equation}

When $n=1$, GREM reduces to REM, another model proposed as a
solvable model by B. Derrida~\cite{De80} in 1980. If $a_1=1$ then
Hamiltonians of REM are nothing but $2^N$ many i.i.d.
$\mathcal{N}(0,N)$ random variables.

Though it was just a toy model, with correct but heuristic arguments
Derrida~\cite{D1} showed phase transition occurs in REM and in the
low temperature the system got completely frozen.  In 1986, B.
Derrida and E. Gardner~\cite{DG86} gave the solution for the
averaged free energy for GREM and in 1987, Capocaccia et
al~\cite{CCP} gave a rigorous mathematical justification. Indeed,
the convergence holds almost surely as well as in $L_p$ for $1\leq p
<\infty$. In 1989, Galves et al~\cite{GMP89} studied the detailed
fluctuation of free energy for both the models  and further analysis
was carried out in 2002 for REM and other models by Bovier et
al~\cite{BKL02}. In a different direction, Dorlas and
Wedagedera~\cite{DW}, in 2001 used the large deviation principle
(LDP)~\cite{DZ,Va84} to study the free energy for REM. In the next
year, Dorlas and Dukes~\cite{DD} extended this technique to GREM.
Though GREM is a little complicated than REM, it is not a realistic
model for spin glasses. More realistic models were proposed earlier
in 1975 by Edwards and Anderson~\cite{EA} (EA-model) through nearest
neighbour interaction and another by Sherrington and
Kirkpatrick~\cite{SK} (SK-model) by mean field correction in the
same year. These are the most complicated models in this theory.
Though several heuristic arguments and conjectures~\cite{MPV87} were
made and several rigorous results were
proved~\cite{ALR87,FZ87,GG98,Ta00,Ta01,Ta02b}, it was only in 2002,
Guerra and Toninelli~\cite{GT} showed the almost sure existence of
the free energy via interpolation technique and convexity argument.
A discussion of the SK-model using stochastic calculus was initiated
by Comets and Neveu~\cite{CN95} continued in~\cite{Co98,AI}. For
EA-model very little has been known till know. In 2003, the idea of
Guerra and Toninelli has been generalized to the GREM cases by
Contucci et al~\cite{CEGG}. We thank these authors for clarifying
their setup.

Note that all this analysis was done with Gaussian driving
distributions. In 2004, Carmona and Hu~\cite{CH} considered
non-Gaussian distributions and showed that the free energy of the
SK-model does not depend on the driving distribution. Rather, under
some moment condition on the driving distributions the free energy
of SK-model is universal (see also~\cite{Ch05}). It should be noted
that earlier already in 1983, Eisele~\cite{E1} considered a class of
distributions with exponentially decaying tails for the REM. He is
the first to identify the relevance of LDP to study free energy for
REM. He studied completely different types of phase transitions --
some kind of iterated large deviation phenomena. For the analysis to
go through, he assumed the existence of exponential moments of all
orders for the driving distributions. The last two articles are the
starting point for this thesis. Now the natural question to ask is,
whether there is any universality of free energy in REM as well as
in GREM? Moreover, is the existence of exponential moments of all
orders necessary?

To answer the above questions our first successful attempt~\cite{N1}
via LDP argument was with double exponential driving distributions.
In \cite{N1}, we provided negative answer to the above questions.
First of all, considering $H_N$ to be i.i.d. double exponential
driving distribution with parameter 1, we show that the nontrivial
free energy is different from that of the Gaussian REM. Though the
Hamiltonian does not depend on $N$, it is interesting, the system
exhibits phase transition. Secondly, note that in this case $E e^{t
H_N}$ does not exists for $t\geq1$. Here in the first chapter, we
extract the essence of the argument in \cite{N1} and state as
\begin{thm}
Let $\{\la_N\}$ satisfies LDP with a strictly quasi-convex rate
function $\mathcal{I}(x)$. For a.e. $\omega$, the sequence of
empirical measure $\{\mu_N(\omega)\}$ of $2^N$ i.i.d. random
variables having law $\la_N$ satisfies LDP with rate function
$\mathcal{J}$ given by,
$$\mathcal{J}(x)=
\begin{cases}
\mathcal{I}(x)& \text{if $\mathcal{I}(x)\leq \log2$}\\
\infty& \text{if $\mathcal{I}(x)>\log2$}.
\end{cases}$$
\end{thm}
We apply this theorem to the known Gaussian case~\cite{D1,DW,OP84},
as well as to double exponential case and further to Weibull type
exponentially decaying tail distributions. We also show that the
energy in REM is not distribution specific rather {\em rate
specific}. In the compact distribution section we give some partial
results when there is no non-trivial rate function for the driving
distributions. In the concluding section, we apply the above theorem
to discrete distributions -- Poisson and Binomial. There we show
that even the existence of phase transition depends on the parameter
of the underlying distributions. For example, if the Hamiltonian
$H_N(\si)$ is Binomial with parameter $N$ and $p$, phase transition
takes place only when $p>\frac{1}{2}$.

For GREM, once again our first attempt~\cite{JR1} was with the
double exponential driving distributions along with the LDP
arguments~\cite{DD}. The original formulation of GREM in the
literature is slightly different from the formulation we mentioned
above. In the second chapter, we start with a discussion of this
reformulation. Then we bring a general tree structure in GREM and
prove a basic fact which is used in the analysis of this chapter as
well as for several other models considered in the next chapter. The
details are in chapter 2. Briefly, we consider trees all of whose
branches extend up to $n$-th level. Let $B_{iN}$ be the total number
of edges at the $i$-th level and $B_N$ be the number of leaves of
the tree. Let $s_{iN}^2$ be the sum of the squares of the numbers of
leaves at the $n$-th level below each edges at the $i$th level. If
$\xi$ denotes a random variable having the common distribution of
the $\xi(\si_1\si_2\cdots\si_i)$, then we have the following.
\begin{thm}\label{t0.2.2}
Let $\tri = \tri_1 \times \cdots \times \tri_n \subset
\mathbb{R}^n$. Denote $q_{iN} = P(\xi \in \tri_i)$ for $1\leq i \leq
n$.

a) If  $\sum\limits_{N \geq n}B_{iN} q_{1N}\cdots q_{iN} < \infty$,
for some $i, 1\leq i \leq n$ then a.s. eventually, $$\mu_N(\tri)=
0.$$

b) If for all $i=1,\cdots,n$, $\sum\limits_{N\geq n}
\frac{s_{iN}^2}{B_N^2 q_{1N}\cdots q_{iN}} < \infty$, then for any
$\ep >0$ a.s. eventually,
\[(1-\ep)\mbf{E}\mu_N(\tri) \leq \mu_N(\tri) \leq (1+\ep)
\mbf{E}\mu_N(\tri).\]
\end{thm}

In section 2.4, we use this result for GREM with a general family of
driving distributions. For fixed $\gamma>0$, we consider the driving
distributions of $\xi(\si_1\cdots\si_i)$ having density
\begin{equation}
\phi_{N,\ga}(x) =
\frac{1}{2\Gamma(\frac{1}{\ga})}\left(\frac{\ga}{N}\right)^{\frac{\ga-1}{\ga}}
e^{-\frac{|x|^{\ga}}{\ga N^{\ga-1}}} \quad -\infty<x<\infty.
\end{equation}
Note that when $\gamma=2$, this is the Gaussian case. We discuss
this case systematically in section 2.5 and bring out the essence of
this model. Here it is.

For each $j$, $1\leq j\leq n$, we have a sequence of probabilities
$\{\la_N^j, N\geq 1\}$ on $\mbb{R}$ satisfying LDP with a good,
strictly quasi convex rate function $\mathcal{I}_j$ and
$\xi(\si_1\cdots\si_i)\sim\la_N^j$. Define for each $\omega$,
$\mu_N(\omega)$ to be the empirical measure on $\mbb{R}^n$, namely,
$$\mu_N(\omega)=\frac{1}{2^N}\sum_\si
\delta\left<\xi(\si_1,\omega),\xi(\si_1\si_2,\omega),\cdots,\xi(\si_1\cdots\si_n,\omega)\right>$$
where $\delta\left<x\right>$ denotes the point mass at
$x\in\mbb{R}^n$.
\begin{thm}\label{t0.2.3}
Suppose $\frac{k(j,N)}{N}\ra p_j>0$ for $1\leq j\leq n$. Then for
a.e. $\omega$, the sequence $\{\mu_N(\omega), N\geq 1\}$ satisfies
LDP with rate function $\mathcal{J}$ given as follows:

Supp$(\mathcal{J})=\{(x_1,\cdots,x_n): \dsum_{k=1}^j
\mathcal{I}_k(x_k) \leq \dsum_{k=1}^j p_k\log2\; \mbox{ for }\,
1\leq j\leq n \}$

and $$\mathcal{J}(x) = \begin{cases}\dsum_{k=1}^n
\mathcal{I}_k(x_k)& \mbox{if }
x\in \mbox{Supp}(\mathcal{J})\\
\infty&\mbox{otherwise}.
\end{cases}$$
\end{thm}

This result, with the help of Varadhan's Integral
lemma~\cite{Va66,DZ}, reduces the problem of free energy to merely
calculation of certain infimum. In section 2.6, we solve this
variational problem for general $n$ and produce the explicit energy
expression in the case of $\gamma>1$ and $\ga=1$ by different
arguments. Further, for $\gamma\geq 1$, we characterize the energy
function for GREM and show that the energy function is continuous in
$\gamma$. For $0<\gamma<1$, we only give the energy expression for
$n=2$. The beauty of the above theorem (Theorem \ref{t0.2.3}) is
that, it allows us to consider different distributions at different
levels of the under lying tree. This we considered in \cite{JR2} and
here in section 2.7. Even the simple case, $n=2$ the model exhibits
a lot of interesting phenomena. For example, consider a 2 level GREM
with exponential driving distribution at the first level and
Gaussian in the second, and give equal weights at the two levels,
that is, $a_1=a_2$. Then even if $p_2=0.00001$ (very small) the
system reduces to a Gaussian REM. On the other hand, if we consider
a 2 level GREM with Gaussian driving distribution at the first level
and exponential in the second, the system will never reduce to a
Gaussian REM. Moreover, in either case, the system will never reduce
to that of an exponential
REM.\\

In the third chapter, we randomize the underlying trees. To keep the
same number of furcations for all edges at given level, for fixed
$N$, we take one Poisson random variable at each level to determine
the number of furcations. We called this model as {\em regular
Poisson GREM}. On the other hand, it is possible to keep the number
of furcations depend on the edge. In other words, for each edge we
can associate a Poisson random variable to determine the number of
furcations for this edge. This we called {\em Poisson GREM}. We
discussed multinomial variation also using results
from~\cite{HMT89}. These are all different methods to randomize the
tree. Note that the configuration space is no longer $2^N$. These
models are interesting and Theorem \ref{t0.2.2} above is powerful
enough to handle these models. However, in all these cases the free
energy remains same as in the usual GREM. Whether there are other
interesting tree structures that exhibit peculiar phenomena is not
clear to us. As far as our knowledge goes, the GREM with randomized
(or even nonrandomized but general) trees is not discussed in the
literature.

In 2006 Bolthausen and Kistler~\cite{BK} proposed a model (BK-GREM)
bypassing the {\em ultrametricity} in the configuration space. Even
in this model, they have shown that the energy of the system is
again a suitable GREM energy. In section 3.4, we provide a proof via
LDP arguments. Then in section 3.5, we construct $n!$ many GREMs,
one corresponding to each permutation of the set $\{1,2,\cdots,n\}$
by manipulating the weights from BK-GREM. We characterize a class of
permutations so that (1) the corresponding GREM energy will be the
same for all the permutations in that class, (2) this energy is the
minimum over all possible $n!$ many GREMs and (3) this is the energy
of the BK-GREM. Bolthausen and Kistler~\cite{BK}, have shown that
the energy of BK-GREM is the infimum over GREM energies
corresponding to all possible {\em chains}. Our analysis shows that
instead of considering all chains, one needs to consider $n!$ many
GREMs. This is still a large number. We conclude this chapter, by
defining one model, called block tree GREM, where the free energy is
maximum of all possible $n!$ many $n$ level GREMs, rather than
minimum as in the BK-GREM.

In the last chapter, we introduce a new model, word GREM. This
brings out the crucial role played by LDP in all the earlier models.
Here we start with a distribution having finite mean and $\la_N$
denotes the law of the sample mean (of size $N$). Then Cramer's
theorem~\cite{DZ} suggests that the sequence $\{\la_N\}$ satisfies
LDP with a convex rate function given in terms of the
Fenchel-Legendre transformation of the starting distribution. We
consider a set of $n$ symbols
$I=\{\varsigma_1,\varsigma_2,\cdots,\varsigma_n\}$ and take $S$ to
be any collection of finite number of words formed by these $n$
symbols. As earlier, we consider $k(i,N)\geq 0, \; 1\leq i\leq n$ as
a partition of $N$ so that the configuration space $2^N$ splits as
$\prod_{i=1}^n 2^{k(i,N)}$. For
$s=\varsigma_{i_1}\varsigma_{i_2}\cdots\varsigma_{i_l}\in S$ and a
configuration $\si=\seq{\si^1,\cdots,\si^n}$ where $\si^i\in
2^{k(i,N)}$, we denote
$\si(s)=\seq{\si^{i_1},\si^{i_2},\cdots,\si^{i_l}}$. For $s\in S$,
let $\la^s$ be a probability on $\mbb{R}$ having finite mean. Let
$\la_N^s$ denote the distribution of the mean of the first $N$
random variables of an i.i.d. sequence with common law $\la^s$. For
the $N$ particle system, we have the following. For each $s\in S$
and each $\si\in2^N$, we have a random variable $\xi(s,\si(s))$.
These are independent random variables. For fixed $s$, they are
identically distributed and the common distribution is $\la_N^s$.
Then for a configuration $\si=\seq{\si_1,\cdots,\si_N}\in 2^N$, we
define the Hamiltonian in word GREM as
\begin{equation}
H_N(\si)=Nf(\xi(\si)) + h\dsum_{i=1}^N\si_i,
\end{equation}
where $f:\mbb{R}^{S}\ra\mbb{R}$ is a continuous function,
$\xi(\si)=(\xi(s,\si(s)))_{s\in S}$ and $h\geq 0$ is the intensity
of the external field.

We present a large deviation proof for the existence of the free
energy for this model and apply the analysis to known~\cite{D1} REM
with external field.

This model includes REM, GREM and BK-GREM and may perhaps include
models truly more general than these. Further, it allows external
field. Moreover, different driving distributions can be used at
different words in the collection.

In this project, we did not consider the analysis of Gibbs'
distribution. With Gaussian driving distribution, there are several
results for REM~\cite{T1,Bo06} and for GREM~\cite{BK04a,BK04b}. See
also~\cite{Ru87,Ko90,KP91,AR05}. For exponential driving
distribution we verified that for REM, in the high temperature
regime, Gibbs's distribution converges to the uniform
distribution~\cite{N1} where as in the low temperature regime it
converges to the Poisson-Dirichlet distribution. This is similar to
that of Gaussian REM. So is it true for any other distributions
considered in this thesis? Since we do not have anything substantial
to say regarding this issue, we have not considered.

\medskip
\section{Large Deviation Terminology}

Recall that $$Z_N(\beta)= \sum_\sigma e^{-\beta H_N(\sigma)}=
2^N \mathbf{E}_\sigma e^{-\beta H_N(\sigma)}$$ where
$\mathbf{E}_\sigma$ is expectation w.r.t uniform probability on
$2^N$ space. And hence,

$$\frac{1}{N}\log Z_N(\beta)=\log 2+ \frac{1}{N}\log\mathbf{E}_\sigma
e^{-N\beta\frac{H_N(\sigma)}{N}}.$$

The last term in the above equation is well known expression in
the Laplace's principle. It is indeed
$$\frac{1}{N}\log\int e^{-Nf}d\mu_N,$$
where $\mu_N$ is the uniform probability on the space $2^N$. The
only trouble is $\mu_N$ are on different spaces. If we transport
$\mu_N$ to $\mbb{R}$ by the map $\si \mapsto\frac{H_N(\si)}{N}$,
then we will arrive at exactly the Laplace type situation, where
Varadhan's integral lemma comes to rescue. Since $H_N$ depends on
$\om$, the transported probability will be random. So the
application of LDP needs careful attention.

Since there are several terminologies (using $\ep$ or using $N$ etc)
for large deviations, we fix our terminology now and recall some
known facts. Let $\mbb{X}$ be a Polish space.
\begin{defn}
A function $\mathcal{I}: \mbb{X}\to \mbb{R}$ is called a lower
semicontinuous function if for any $a\in \mbb{R}$, the set
$\{x:\;\mathcal{I}(x)\leq a\}$ is a closed set. It will be called
good  if the set $\{x:\;\mathcal{I}(x)\leq a\}$ is a compact set.
\end{defn}

The following two properties of lower semicontinuous function are
worth mentioning.
\begin{prop}
Let $f$ be a lower semicontinuous function. Then for any $x$,
$$\dsup_{G:\, \mbox{\em{neighbourhood of} } x} \inf_{y\in G} f(y)=f(x).$$
\end{prop}

\begin{prop}\label{p0.3.2}
Let $f$ be a good function and $\{F_n\}_n$ be a sequence of closed
sets so that $F_{n+1}\subseteq F_n$ for every $n$ and $\cap_n
F_n=\{x_0\}$. Then
$$f(x_0)=\lim_n \inf_{y\in F_n} f(y).$$
\end{prop}

\begin{defn}
Let $\{\mu_N\}$ be a sequence of probabilities on $\mbb{X}$.
$\{\mu_N\}$ is said to satisfy large deviation principle with rate
function $\mathcal{I}$ if
\begin{enumerate}
\item {$\mathcal{I}: \mbb{X}\to [0,\infty]$ is a lower semicontinuous function,}
\item{for any Borel set $B$,
$$-\inf_{x\in B^0} \mathcal{I}(x)\leq \liminf_N\frac{1}{N}\log\mu_N(B)\leq
\limsup_N\frac{1}{N}\log\mu_N(B) \leq -\inf_{x\in \ol{B}} \mathcal{I}(x).$$}
\end{enumerate}
Further, if for $0\leq a<\infty$, the set $\{x:\;\mathcal{I}(x)\leq
a\}$ is a compact set, then $\mathcal{I}$ is called a good rate
function.
\end{defn}

A sufficient condition for the existence of LDP is the following:
\begin{prop}\label{p0.3.3}
Let $\mbb{X}$ be a Polish space. Let  $\mathcal{A}$ be an open base
for $\mbb{X}$. Let $\{\mu_N\}$ be a sequence of probabilities on
$\mbb{X}$. For each $A\in\mathcal{A}$, let $L_* (A)=-\liminf_N
\frac{1}{N} \log \mu_N(A)$ and $L^* (A)=-\limsup_N \frac{1}{N} \log
\mu_N(A)$. Suppose for every $x\in \mbb{X}$, $$\sup\limits_{x\in A
\in\mathcal{A}}L_*(A)=\sup\limits_{x\in A
\in\mathcal{A}}L^*(A)=\mathcal{I}(x)\quad\text{ (say)}.$$ Assume
moreover that either $\{\mu_N\}$ is eventually supported on a
compact set or the sequence is exponentially tight, that is, given
any $\al<\infty$, there is a compact set $K$ such that
$\limsup_N\frac{1}{N} \log \mu_N (K^c)<-\al$.

Then the sequence $\{\mu_N\}$ satisfies LDP with rate function $\mathcal{I}$.
\end{prop}

The next proposition is a variation of well-known Varadhan's
integral lemma, which will suggest that only we have to calculate
some infimum to get the free energy limit.
\begin{prop}\label{p0.3.5}
Suppose the sequence of probabilities $\{\mu_N\}$ on a Polish space
$\mbb{X}$ satisfies LDP with rate function $\mathcal{I}$ and,
moreover, $\mu_N$ are eventually supported on a compact set $C$.
Then for any continuous function $f: \mbb{X} \rightarrow
\mathbb{R}$,
$$\lim_{N\rightarrow \infty} \frac{1}{N}\log \int e^{-Nf} d\mu_N = -
\inf_{x\in C} \{f(x)+\mathcal{I}(x)\}.$$
\end{prop}

We need the following known as Cramer's Theorem.
\begin{thm}\label{p0.3.6}
Let $X_1,X_2,\cdots$ be i.i.d. real valued random variables with $E
X_1<\infty$ and $\la_N$ be the law of their sample mean. Then the
sequence of probabilities $\{\la_N\}$ satisfies LDP with a convex rate
function $\mathcal{I}$ given by
\begin{equation}
\mathcal{I}(x)=\sup_{\mathcal{D}} \{\theta x - \log E e^{\theta
X_1}\},
\end{equation}
where $\mathcal{D}=\{\theta:\;\log E e^{\theta
X_1}<\infty\}\subseteq\mbb{R}$.
\end{thm}

The following is also true.
\begin{lem}
Let $\mathcal{I}$ be as in the above theorem. If
$\mathcal{D}_\mathcal{I}=\{x:\;\mathcal{I}(x)<\infty\}$ and
$\ol{x}=E X_1<\infty$, then

1. $\mathcal{I}(\ol{x})=0$,

2. $\mathcal{I}$ is strictly decreasing on $\{x\leq
\ol{x}\}\cap\mathcal{D}_\mathcal{I}$,

3. $\mathcal{I}$ is strictly increasing on $\{\ol{x}\leq
x\}\cap\mathcal{D}_\mathcal{I}$.

\end{lem}


\cleardoublepage

\chapter{The Random Energy Model}



In this chapter we discuss a toy model of spin glass theory, called
Random Energy Model (REM). In the literature\cite{D1,T1}, this is
driven with Gaussian distributions. In \cite{E1}, Eisele discussed
the model with more general distributions, particularly with
regularly varying distributions and the relevance of large deviation
methods in this context. We study the model with other types of
distributions. For instance, the driving distributions could be
exponential or more generally Weibull. Or they could be compactly
supported etc. Our discussion mainly relies on the idea of Dorlas
and Wedagedera. In \cite{DW}, they first used the large deviation
techniques to get the asymptotics of the free energy.

After defining the model in the first section, we give a general
large deviation result in the section 1.2 and apply the results in
REM with diverse distributions in section 1.3. In section 1.4, we
give partial results with compact distributions where we could not
use the large deviation results. We conclude this chapter by
considering the model driven by some discrete distributions.

\goodbreak
\section{Setup}

In this model,  proposed originally by B. Derrida\cite{De80}, for
each $N$, the Hamiltonian $H_N(\sigma)$ are independent over
$\si\in\Si_N$. Derrida considered them to be centered Gaussian with
variance $N$. In spite of the simplicity of this model, in
\cite{D1}, he showed the existence of phase transition. Using the
entropy energy equation, he evaluated the limiting annealed free
energy $\lim_{N}\frac{1}{N}E\log Z_N(\be)$ and showed that for low
temperature, that is, for $\be$ large, free energy becomes linear in
$\be$. It is known that, in fact, $\frac{1}{N}\log Z_N(\be)$
converges a.s.\cite{DW}.

\section{Main Results}
Let us consider a sequence of probabilities $\left(\la_N, N\geq1
\right)$ on $\mbb{R}^n$. Assume that $\{\la_N\}$ satisfies large
deviation principle (LDP) with a strictly quasi-convex good rate
function $\mathcal{I}(x)$. An extended real valued function $f$,
which may take the value $+\infty$ but not $-\infty$, defined on a
convex set will be called a strictly quasi-convex function if for
any two distinct $x_1$ and $x_2$ in $\{f(x)<\infty\}$ and for each
$\theta\in (0,1)$ we have $f(\theta x_1 +
(1-\theta)x_2)<\max\{f(x_1), f(x_2)\}$. For every $N$, let $\xi_i$,
$1\leq i\leq 2^N$ be i.i.d. random variables ($\mathbb{R}^n$ valued)
with distribution $\la_N$. These random variables, of course, depend
on $N$ but to ease the notation we are suppressing their dependence
on $N$. For every sample point $\omega$, we define $\mu_N(\omega)$
to be the empirical measure on $\mbb{R}^n$, namely
$\mu_N(\omega)=\frac{1}{2^N}\sum \delta\left<\xi_i(\omega)\right>$.
Here $\delta\left<x\right>$ denote the point mass at $x$. Now we are
ready to state our first theorem.

\begin{thm}\label{t1.2.1}
For a.e. $\omega$ the sequence $\{\mu_N(\omega)\}$ is supported on a
compact set and satisfies LDP with rate function $\mathcal{J}$ given
by,
\begin{equation*}
\mathcal{J}(x) =\begin{cases}\mathcal{I}(x) &\mbox{if } \mathcal{I}(x)\leq \log2\\
\infty &\mbox{if } \mathcal{I}(x)>\log2.\end{cases}
\end{equation*}

\end{thm}

\begin{proof} {\bf Step 1:} {\em Let $\tri$ be an open subset of
$\mbb{R}^n$. If $\sum 2^N \la_N(\tri)<\infty$, then almost surely
eventually
$\mu_N(\tri)=0.$}\\
Indeed, using $P$ for the probability on the space where the random
variables are defined,
$$P(\mu_N(\tri)>0)= P(\xi_i \in \tri \mbox{ for some } i) \leq2^N \la_N(\tri).$$
Now Borel - Cantelli completes the proof.

{\bf Step 2:} {\em Let $\tri$ be an open subset of $\mbb{R}^n$. If
$\sum \frac{1}{2^N \la_N(\tri)} <\infty$, then for any $\ep>0$,
almost surely eventually $$(1-\ep) \la_N(\tri) \leq \mu_N(\tri) \leq
(1+\ep)\la_N(\tri).$$} Indeed, $$ \mbox{Var }\mu_N(\tri) =
\mbf{E}\left(\frac{1}{2^N}\sum 1_{\tri}(\xi)\right)^2 -
\la_N^2(\tri) \leq \frac{1}{2^N} \la_N(\tri).$$ Now Chebyshev yields
$$P\left\{|\mu_N(\tri)-\la_N(\tri)|> \ep \la_N(\tri)\right\} \leq \frac{1}{\ep^2 2^N \la_N(\tri)}$$
and Borel-Cantelli completes the proof.\\

Since $\mathcal{I}$ is strictly quasi-convex, the set
$\{\mathcal{I}(x)=\log2\}$ does not contain any line segment, we can
choose a countable open base $\mathfrak{B}$ such that for every
$\tri\in \mathfrak{B}$ either
$\ol{\tri}\cap\{\mathcal{I}(x)\leq\log2\}=\varnothing$ or
$\tri\cap\{\mathcal{I}(x)<\log2\}\neq\varnothing$. For instance, we
could choose $\mathfrak{B}$ to be the collection all open boxes such
that (i) $\mathcal{I}$ value at a corner point is different from
$\log 2$; (ii) each co-ordinate of a corner point is either rational
or $\pm \infty$.

{\bf Step 3:} {\em Let $\mathcal{I}(x)>\log2$. Then almost surely,
$\sup\limits_{x\in\tri\in\mathfrak{B}}\{-\liminf\frac{1}{N}\log\mu_N(\tri)\}$
as well as
$\sup\limits_{x\in\tri\in\mathfrak{B}}\{-\limsup\frac{1}{N}\log\mu_N(\tri)\}$
are $\infty$.}

Since $\mathcal{I}(x)>\log2$, pick $\tri_0\in\mathfrak{B}$ such that
$x\in \tri_0$ and $\ol{\tri}_0\cap
\{\mathcal{I}(x)\leq\log2\}=\varnothing$. Then
$\limsup\frac{1}{N}\log \la_N(\tri_0)\leq
-\inf\limits_{y\in\ol{\tri}_0}\mathcal{I}(y)= -L<-\log2$. Fix
$\al>0$ such that $-L<-\log2-\al$. Then for sufficiently large $N$,
$\frac{1}{N}\log \la_N(\tri_0) \leq -\log2-\al$, that is,
$\la_N(\tri_0)\leq 2^{-N}e^{-N\al}$. In other words,
$2^N\la_N(\tri_0) \leq e^{-N\al}$ for all large $N$. Thus by Step 1,
a.s. eventually
$\mu_N(\tri_0)=0$ and hence the claim.

Incidentally, this also shows the following. If $K$ is the compact
set $\{x:\, \mathcal{I}(x)\leq\log2\}$ then consider a bounded open
box $\tri\in\mathfrak{B}$ such that $K\subset \tri$. Clearly,
$\ol{\tri}^c$ is union of $2n$ many boxed from $\mathfrak{B}$. The
above argument shows that $\mu_N$ is a.s. eventually zero for each
of these $2n$ boxes. This shows that the sequence $\{\mu_N\}$ is
a.s. eventually
supported on a compact set, namely, $\ol{\tri}$.\\

{\bf Step 4:} {\em Let $\mathcal{I}(x)\leq\log2$. Then almost
surely,
$\sup\limits_{x\in\tri\in\mathfrak{B}}\{-\liminf\frac{1}{N}\log\mu_N(\tri)\}$
as well as
$\sup\limits_{x\in\tri\in\mathfrak{B}}\{-\limsup\frac{1}{N}\log\mu_N(\tri)\}$
are $\mathcal{I}(x)$.}

Fix $\tri\in \mathfrak{B}$ such that $x\in\tri$. Then
$\liminf\frac{1}{N}\log \la_N(\tri)\geq
-\inf\limits_{y\in\tri}\mathcal{I}(y)= -L>-\log2$, where the last
inequality is a consequence of the strict quasi-convexity of
$\mathcal{I}$. Fix $\al>0$ so that $-L>-\log2 +\al$. So for large
$N$, $\frac{1}{N}\log \la_N(\tri)>-\log2 +\al$, that is,
$\la_N(\tri)>2^{-N}e^{N\al}$. In other words, $2^N\la_N(\tri)>
e^{N\al}$. Now use Step 2, for any $\ep\in (0,1)$, eventually
$$(1-\ep) \la_N(\tri) \leq \mu_N(\tri) \leq
(1+\ep)\la_N(\tri).$$
Hence by definition of LDP, we have eventually
\begin{align}\label{e1.2.0}
- \mathcal{I}(\tri)\leq \liminf\limits_{N\ra \infty}\frac{1}{N} \log
\mu_N(\tri)\leq \limsup\limits_{N\ra \infty}\frac{1}{N} \log
\mu_N(\tri)\leq -\mathcal{I}(\ol{\tri}),\end{align} where as usual
$\mathcal{I}(A)=\inf\limits_{x\in A} \mathcal{I}(x)$.

From the first part of the
above inequality we have,
\begin{equation}\label{e1.2.0a}
\sup_{\tri\in\mathfrak{B}:x\in \tri}\{-\liminf_N\frac{1}{N}\log
\mu_N(\tri)\}\leq \sup_{\tri\in\mathfrak{B}:x\in \tri}
\mathcal{I}(\tri)\;\leq\; \mathcal{I}(x).
\end{equation}

Moreover, for every $\tri \in \mathfrak{B}$ such that $x\in\tri$ using the right side inequality of
(\ref{e1.2.0}), we have
$$\limsup_N\frac{1}{N}\log \mu_N(\tri)\leq -\mathcal{I}(\ol{\tri}).$$

Let $\mathfrak{B}_{x}=\{\tri_k \in \mathfrak{B}: k\geq
1\}$ be a subclass of $\mathfrak{B}$ so that
$\ol{\tri}_{k+1}\subset\tri_k$ for every $k$ and
$\cap_k\tri_k=\{x\}.$ Then
\begin{align}\label{e1.2.0b}
\sup_{\tri\in\mathfrak{B}:x\in \tri}\{-\limsup_N\frac{1}{N}\log
\mu_N(\tri)\}&\geq \sup_{\tri\in\mathfrak{B}:x\in \tri}
\mathcal{I}(\ol{\tri})\notag\\
&\geq \sup_{\tri\in \mathfrak{B}_{x}}
\mathcal{I}(\ol{\tri})\notag\\
&=\lim_k \mathcal{I}(\ol{\tri}_k)\notag\\ &=
\mathcal{I}(x).
\end{align}
The last equality follows as $\mathcal{I}$ is a good lower
semicontinuous function (see Proposition \ref{p0.3.2}).

From (\ref{e1.2.0a}) and (\ref{e1.2.0b}), it follows that
$$\sup_{\tri\in\mathfrak{B}:x\in \tri}\{-\liminf_N\frac{1}{N}\log
\mu_N(\tri)\}=\sup_{\tri\in\mathfrak{B}:x\in
\tri}\{-\limsup_N\frac{1}{N}\log \mu_N(\tri)\}=\mathcal{I}(x).$$

Now proof of Theorem \ref{t1.2.1} is completed by appealing to
Proposition \ref{p0.3.3} and observing that $\{\mu_N\}$ is
eventually supported on a compact set.
\end{proof}

\begin{rem}\label{r1.2.1}
The fact that $\mathcal{I}$ is a {\em good} rate function is
essential in the above theorem to conclude that almost surely
eventually the sequence $\{\mu_N(\om)\}$ is supported on a compact
subset of $\mbb{R}^n$.
\end{rem}

\begin{rem}\label{r1.2.2}
Observe that the strict quasi-convexity of the rate function in the
above theorem is a technical assumption. On real line that
assumption can be replaced by the assumption:  $\mathcal{I}$ is
strictly monotone on $\{x: \mathcal{I}(x)\in (0, \infty)\}$ or by
the assumption that $\{x: \mathcal{I}(x)=\log2\}$ is a nowhere dense
set. Such a condition is needed only to ensure that there exists a
countable base as mentioned in the above proof.
\end{rem}

The implication of the above theorem in REM is amazing. To see this,
let us assume that $\{\la_N\}$ is a sequence of probabilities on
$\mbb{R}$ and satisfies large deviation principle with a good rate
function $\mathcal{I}$. Let us assume also that $\mathcal{I}$ be
strictly quasi-convex or satisfies any one of the conditions in
remark \ref{r1.2.2}. For fixed $N$, let us consider $2^N$ i.i.d.
random variables $\xi(\si)$, $1\leq \si \leq 2^N$ distributed like
$\la_N$. We can identify these $2^N$ many $\si$ with the elements of
$\Si_N=\{+1, -1\}$. Let us define the Hamiltonian for $\si \in
\Si_N$ to be
$$H_N(\si)=N\xi(\si).$$
Now note that the partition function can be written as
$$Z_N(\be)=2^NE_\si e^{-\be H_N(\si)},$$ where $E_\si$ is the
expectation with respect to uniform probability on the $\Si_N$
space. Hence $$\lim_N \frac{1}{N}\log Z_N(\be)= \log2 + \lim_N
\frac{1}{N} \log \int e^{-N\be x}d\mu_N(x).$$

By Theorem \ref{t1.2.1}, the induced probabilities $\{\mu_N\}$ are
a.s. eventually supported on a compact set. That is, for almost
every $\omega$, there is a compact set $K_\omega$ such that
$\{\mu_N(\omega),\, N\geq 1\}$ are all supported on $K_\omega$.
Moreover, by previous theorem they satisfy LDP with a good rate
function a.s. Hence to find the existence of $\lim_N \frac{1}{N}\log
Z_N(\be)$, we can use Varadhan's integral lemma with any continuous
function, in particular, $f(x)=\be x$ on $\mbb{R}$. This will lead
to the following:
\begin{thm}\label{t1.2.2}
If the sequence of probabilities $\{\la_N\}$ satisfies LDP with
strictly quasi-convex good rate function $\mathcal{I}$, then almost
surely,
\begin{equation}\label{e1.2.1a}
\dlim_{N\ra \infty} \frac{1}{N} \log Z_N(\be) = \log2\, -
\dinf_{\{\mathcal{I}(x)\leq \log 2\}} \left(\be x +
\mathcal{I}(x)\right).
\end{equation}
\end{thm}

Thus in the REM, the existence of limiting free energy is just a
corollary of large deviation principle. To get the the expression of
the free energy one has to solve the variational formula. Hence the
calculation of asymptotics of free energy reduces to calculation of
the above infimum.

\begin{rem}
In the literature\cite{D1,DW,T1}, for REM, the Hamiltonian
$H_N(\si)$ is defined as $\sqrt{N}\xi(\si)$, with $\xi(\si)\thicksim
\mathcal{N}(0,N)$. In our case with the Gaussian driving
distribution it is same as $N\xi(\si)$ where $\xi(\si)\thicksim
\mathcal{N}(0,\frac{1}{N})$. But the large deviation technique allows
us to consider $H_N$ to be any continuous function of $\xi(\si)$. In
other words, if $f(x)$ is a continuous function on $\mbb{R}$ then
one can define the random Hamiltonian $H_N(\si)=Nf(\xi(\si))$ where
$\xi(\si)\thicksim \la_N$. In that case, the above theorem will
reduce to
\end{rem}
\begin{thm}\label{t1.2.3}
If a sequence of probabilities $\{\la_N\}$ on $\mbb{R}$ satisfies
LDP with strictly quasi-convex good rate function $\mathcal{I}$ and
$H_N(\si)=Nf(\xi(\si))$ where $\xi(\si)\thicksim \la_N$ and $f$ is a
continuous function on $\mbb{R}$, then almost surely,
\begin{equation}\label{e1.2.1b}
\dlim_{N\ra \infty}
\frac{1}{N} \log Z_N(\be) = \log2\, - \dinf_{\{\mathcal{I}(x)\leq
\log 2\}} \left(\be f(x) + \mathcal{I}(x)\right).
\end{equation}
\end{thm}

Of course, the appearance of $f$ above makes it more general, but
this could be obtained from (\ref{e1.2.1a}) by contraction principle
of large deviation techniques. Different choices of functions $f$
allows us to consider the Hamiltonian driven by other distribution
which can be obtained as a function of known distributions. For
instance, if we consider $f(x)=x^2$ then we can get the information
of the model when its Hamiltonian is  an appropriate $\chi^2$ if
$\la_N$ as $\mathcal{N}(0,\frac{1}{N})$. Also we can consider
several other functions, for which we do not know the corresponding
closed form expression of the distribution of the Hamiltonian. For
example, $f(x)=x \cos(1000\pi |x|^{|x|})$ etc.

\begin{rem}
A close look at $\dinf_{\{\mathcal{I}(x)\leq \log 2\}} \left(\be
f(x) + \mathcal{I}(x)\right)$ suggests the following: If we consider
the Hamiltonian of an $N$-particle system to be an odd function $f$
of random variables $\xi(\si)\thicksim \la_N$ and if the sequence
$\{\la_N\}$ satisfies LDP with a quasi-convex good rate function
taking value 0 at the origin, then the contribution for the limiting
free energy $\dlim_{N\ra \infty} \frac{1}{N} \log Z_N(\be)$ comes
from only that part where the function $f(x)$ is negative. More
precisely, $f$ being an odd function $f(0)=0$. Since
$\mathcal{I}(0)=0$, infimum in (\ref{e1.2.1b}) is non-positive. So
only points $x$ where $f(x)\leq 0$ need to be considered while
calculating the infimum. Thus, we have, almost surely, $$\dlim_{N\ra
\infty} \frac{1}{N} \log Z_N(\be) = \log2\, -
\dinf_{\{\mathcal{I}(x)\leq \log 2\}^-} \left(\mathcal{I}(x) + \be
f(x)\right),$$ where $\{\mathcal{I}(x)\leq \log
2\}^-=\{\mathcal{I}(x)\leq \log 2\}\cap \{f(x)\leq 0\}$. For
example, when $f(x)=x$  we have the following
\begin{cor}\label{c1.2.3}
If $f(x)=x$ and $\mathcal{I}(0)=0$, then almost surely,
$$\dlim_{N\ra \infty} \frac{1}{N} \log Z_N(\be) = \log2\, -
\dinf_{\{\mathcal{I}(x)\leq \log 2\}^-} \left(\mathcal{I}(x) + \be x
\right).$$
\end{cor}

One can see that the contribution to the free energy is only from
the negative values of the random variable. This can be made precise
as follows: Let $\{\la_N\}$ and $\{\nu_N\}$ be two sequences of
probabilities satisfying LDP with a good strictly quasi-convex rate
functions $\mathcal{I}_1$ and $\mathcal{I}_2$ respectively so that
$\mathcal{I}_1(0)=\mathcal{I}_2(0)=0$ and
$\mathcal{I}_1(x)=\mathcal{I}_2(x)$ for $x\leq 0$. Then
consideration of $\xi_N\sim\la_N$ or $\xi_N\sim\nu_N$ will lead to
the same limiting free energy.  In other words, symmetry of the
random variables does not play any role in the evaluation of
limiting free energy. To illustrate, if we consider the density of
$\la_N$ given by,
\begin{equation}\label{e1.2.1}
\phi_N(x)=\begin{cases}\frac{1}{2}\sqrt{\frac{N}{2\pi}}e^{-\frac{1}{2}N
x^2}& \text{for $x\geq 0$}\\ \frac{N}{2} e^{Nx}& \text{for
$x<0$}\end{cases},\end{equation} then from the discussion of our
next section, it will be clear that this sequence $\{\la_N\}$
satisfies LDP with rate function
$$\mathcal{I}(x)=\begin{cases}\frac{x^2}{2}& \text{for  $x\geq 0$},\\
-x& \text{for $x<0$}\end{cases}$$ and hence
$\{\mathcal{I}(x)\leq\log2\}=[-\log2,\sqrt{2\log2}]$. Here the
distribution of $\la_N$ is of Gaussian form in the positive part of
the real line whereas on the negative part it is of exponential
nature.  If $\xi_N\thicksim \la_N$, then as $\dinf_{\{0< x\leq
\sqrt{2\log2}\}} \left(\frac{1}{2}x^2 + \be x \right)\geq0$ the
above corollary will reduce to
\begin{cor}
If $\la_N$ has density $\phi_N$ given by (\ref{e1.2.1}), then almost
surely, $$\dlim_{N\ra \infty} \frac{1}{N} \log Z_N(\be) = \log2\, -
\dinf_{\{-\log 2\leq x\leq 0\}} \left(-x + \be x \right).$$
\end{cor}
Hence the Gaussian part of the random variables does not
contribute to the limiting free energy. Similarly, if we
consider the density of $\la_N$ to be
\begin{equation}\label{e1.2.2}
\phi_N(x)=\begin{cases}\frac{N}{2} e^{-Nx}& \text{for $x\geq
0$}\\
\frac{1}{2}\sqrt{\frac{N}{2\pi}}e^{-\frac{1}{2}N x^2} & \text{for
$x<0$}\end{cases},\end{equation} then the rate function will
be $$\mathcal{I}(x)=\begin{cases}x& \text{for $x\geq 0$}\\
\frac{1}{2}x^2& \text{for $x<0$}.\end{cases}$$ In this case
Corollary \ref{c1.2.3} will reduce to
\begin{cor}
If $\la_N$ has density $\phi_N$ given by (\ref{e1.2.2}), then almost
surely, $$\dlim_{N\ra \infty} \frac{1}{N} \log Z_N(\be) = \log2\, -
\dinf_{\{-\sqrt{2\log2}\leq x\leq 0\}} \left(\frac{1}{2}x^2 + \be x
\right).$$\end{cor}
Here the exponential nature of the random
variable on the positive side does not play any role.
\end{rem}

\begin{rem}\label{r1.2.6}
Suppose that the sequence of probabilities $\{\la_N\}$ is supported
on $[0,\infty)$ and satisfies LDP with a quasi-convex rate function
$\mathcal{I}$ with $\mathcal{I}(0)=0$. For example, fix a number
$\ga>0$ and put $\mathcal{I}(x)= x^\gamma$ for $x\geq 0$ and
$\infty$ for $x<0$ is such a rate function. It follows from
Corollary \ref{c1.2.3} that, when $f(x)=x$, then almost surely
$$\dlim_{N\ra \infty} \frac{1}{N} \log Z_N(\be) = \log2\, -
\dinf_{\{\mathcal{I}(x)\leq \log 2\}^-} \left(\mathcal{I}(x) + \be x
\right).$$ As the sequence $\{\la_N\}$ is supported on non-negative
real line, $\mathcal{I}(x)=\infty$ for $x<0$. Hence
$\{\mathcal{I}(x)\leq \log 2\}^-=\{\mathcal{I}(x)\leq \log 2\}\cap
\{f(x)\leq 0\}=\{0\}$ and
$$\dlim_{N\ra \infty} \frac{1}{N} \log Z_N(\be) =
\log2$$ almost surely. In this case, the system will not show phase
transition.
\end{rem}

The examples given above are rather artificial. Of course, there are natural examples of random
variables $\xi(\si)$ whose distributions satisfy large deviation
principle with a good convex rate function. In the following
sections, we discuss some examples.

\section{Distribution with exponentially decaying Tail}

In this section, we consider the driving sequence of
distributions ($\la_N$) such that for $x>0$, $\la_N[-x,x]^c\sim
e^{-Nx^\ga}$ for some $\ga>0$.
\subsection{Gaussian Distribution}\label{s1.3.1}

Our first natural example is the Gaussian distribution well studied
in the literature~\cite{D1,OP84,DW,T1}. Let $\la_N$ be the centered
Gaussian probability with variance $\frac{1}{N}$. That is, having
density $\sqrt{\frac{N}{2\pi}}e^{-N\frac{x^2}{2}}$, for $-\infty < x
< \infty$. It is obvious that $\la_N\Ra 0$ as $N\ra \infty$. The
following is well known. It can also be obtained from Cramer's
theorem~\ref{p0.3.6}. Since the proof is simple, we give it.

\begin{prop} The sequence $\{\la_N\}$ satisfies LDP with
rate function $\mathcal{I}= \frac{x^2}{2}$ on $\mbb{R}$.
\end{prop}

\begin{proof}
Let $\tri\subset \mbb{R}$ be an open interval. Let $m = \dsty\inf_{x
\in \tri} |x|$, $M = \dsty\sup_{x \in \tri}|x|$, and $q_N=
\la_N(\tri)$. With this notation, note that, we have
\begin{equation}
q_N \leq 2\sqrt{\frac{N}{2\pi}}\dsty\int_{m}^{M} e^{-N\frac{x^2}{2}}
dx <\dsty\int_{\sqrt{N}m}^{\sqrt{N}M} e^{-\frac{x^2}{2}} dx<
\dsty\int_{\sqrt{N}m}^{\infty} e^{-\frac{x^2}{2}} dx \leq
\frac{1}{\sqrt{N}m}e^{-\frac{Nm^2}{2}},
\end{equation}
with the understanding that when $m_i=0$, the last expression is
$\frac{1}{2}$ and
\begin{equation}
q_N \geq \frac{1}{\sqrt{2\pi}} \dsty\int_{\sqrt{N}m}^{\sqrt{N}M}
e^{-\frac{x^2}{2}} dx
> \frac{1}{2} \dsty\int_{\sqrt{N}m}^{\sqrt{N}(m + \de)} e^{-\frac{x^2}{2}} dx >
\frac{\sqrt{N}\de}{2} e^{-\frac{N}{2}(m + \de)^2},
\end{equation}
for any $0<\de < M - m$.

From above two inequalities, we can conclude that for any open
interval $\tri$, the limit, $\dlim_{N\ra \infty} \frac{1}{N} \log
\la_N(\tri)= - \frac{m^2}{2}.$ Once again, Proposition \ref{p0.3.3}
completes the proof.
\end{proof}

\begin{rem}
Note that, here $\mathcal{I}(x)=\frac{x^2}{2}$ is a continuous
function with compact level sets. Not only that, it is a convex
function and hence quasi-convex.
\end{rem}

As a consequence of Theorem \ref{t1.2.3}, if for $\si\in\Si_N$ the random
variables $\xi(\si)\thicksim\la_N$ and the Hamiltonian $H_N(\si)=N
f(\xi(\si))$ with any continuous function $f$ on $\mbb{R}$, we get
the following:

\begin{cor}
If $\la_N\thicksim \mathcal{N}(0,\frac{1}{N})$, then almost surely,
$$\dlim_{N\ra \infty} \frac{1}{N} \log Z_N(\be) = \log2\, -
\dinf_{x^2\leq 2\log 2} \left(\be f(x) + \frac{x^2}{2}\right).$$
\end{cor}

Taking $f(x)=x$, we will get the classical case where the
Hamiltonian for $N$-particle system $H_N$ is Gaussian with mean 0
and variance $N$. Note that, $\dinf_{x^2\leq 2\log 2} \left(\be x +
\frac{x^2}{2}\right)=\dinf_{0\leq x \leq \sqrt{2\log 2}}
\left(\frac{x^2}{2}-\be x\right)$. Let us denote the function
$g(x)=\frac{x^2}{2}-\be x$ so that $g'(x)=x-\be$ and $g''(x)=1>0$.
Therefore, at $x=\be$ the function $g$ attains its infimum. So as
long as $\be \leq \sqrt{2\log2}$, the $\dinf_{0\leq x \leq
\sqrt{2\log 2}} \left(\frac{x^2}{2}-\be x\right)$ is attained at
$x=\be$. Moreover, as $g$ is a decreasing function on $[0,\be]$, for
$\be>\sqrt{2\log2}$ the $\dinf_{0\leq x \leq \sqrt{2\log 2}}
\left(\frac{x^2}{2}-\be x\right)$ is attained at $x=\sqrt{2\log 2}$.
Hence we get the following

\begin{thm}
If $H_N(\si)$ are independent $\mathcal{N}(0,N)$, then almost
surely,
$$\begin{array}{llll}
\dlim_N \frac{1}{N}\log Z_N(\be)&=&\log2+\frac{\be^2}{2} &\mbox{ if
} \be<
\sqrt{2\log2}\\
&=&\be\sqrt{2\log2 }&\mbox{ if } \be\geq \sqrt{2\log2}.
\end{array}$$
\end{thm}
As we mentioned at the beginning, this is classical.

\subsection{Exponential Distribution}\label{s1.3.2}
Another simple but interesting example is the exponential
distribution. Let $\la_N$ be two sided exponential probability with
scale parameter $\frac{1}{N}$. That is, having density
$\frac{1}{2}Ne^{-N|x|}$, for $-\infty< x < \infty$. Once again, it
is obvious that $\la_N\Ra 0$ as $N\ra \infty$. Now we show the
following

\begin{prop} The sequence $\{\la_N\}$ satisfies LDP with
rate function $\mathcal{I}= |x|$ on $\mbb{R}$.
\end{prop}
\begin{proof}
Let $\tri\subset \mbb{R}$ be an interval. Let $m =
\inf\{|x|:\, x \in \tri\}$, $M = \sup\{|x|:\, x \in \tri\}$, and $q_N=
\la_N(\tri)$. With this notation, we have
$$q_N =\frac{N}{2}\dsty\int_\tri e^{-N|x|} dx < \dsty\int_{\sqrt{N}m}^{\infty}
e^{-x} dx \geq e^{-Nm},$$
and
$$q_N \geq \dsty\int_{Nm}^{NM} e^{-x} dx
> \dsty\int_{Nm}^{N(m + \de)} e^{-x} dx >
N\de e^{-N(m + \de)},$$
for any $0<\de < M - m$.

From above two inequalities, we can conclude that for any interval
$\tri$, the limit, $\dlim_{N\ra \infty} \frac{1}{N} \log
\la_N(\tri)= - m.$ Once again, Proposition \ref{p0.3.3} completes
the proof.
\end{proof}

\begin{rem}
As in the Gaussian case, here also, note that, $\mathcal{I}(x)=|x|$
is a convex continuous function with compact level sets.
\end{rem}

\begin{cor}
If $\xi_N(\si)$ are independent (over $\sigma$) two sided
exponential variables with scale parameter $\frac{1}{N}$, $f$ is
a continuous function on $\mathbb{R}$ and $H_N(\si)=Nf(\xi(\si))$
then almost surely,
$$\dlim_{N\ra \infty} \frac{1}{N} \log Z_N(\be) = \log2\, -
\dinf_{|x|\leq 2\log 2} \left(\be f(x) + |x|\right).$$
\end{cor}

If we take $f(x)=x$, the Hamiltonian for $N$-particle system $H_N$
is, of course, two sided exponential random variables with scale
parameter 1, that is have density $\frac{1}{2}e^{-|x|}$ on
$\mbb{R}$. In that case, for the limiting free energy, we only need
to calculate $\dinf_{|x| \leq 2\log 2} \left(|x|+\be x\right)$, that
is, $\dinf_{0\leq x \leq 2\log 2} x(1-\be)$. A simple calculation
yields the following

\begin{thm}
If $H_N(\si)$ are independent two sided exponential random variables
with scale parameter 1, then almost surely,
$$\begin{array}{llll}
\dlim_N \frac{1}{N}\log Z_N(\be)&=&\log 2 &\mbox{ if } \be< 1\\
&=&\be\log2 &\mbox{ if } \be\geq 1.
\end{array}$$
\end{thm}
\begin{rem}
Interesting observation in this analysis is that, the random
Hamiltonian, for $N$-particle system being exponential random
variables with scale parameter 1, does not depend on the number of
particles. Even then, the system shows a phase transition.
\end{rem}

\subsection{Weibull Distribution}
A more general class that can be considered are the Weibull
distributions.  Let $\la_N$ be the probabilities on $\mathbb{R}$
having density
\begin{equation}\phi_{N,\ga}(x) = \frac{N}{2}
|x|^{\ga-1}e^{-N\frac{|x|^\ga}{\ga}}, \quad -\infty<
x<\infty.\end{equation} This is known as Weibull distribution with
shape parameter $\ga>0$ and scale parameter
$\left(\frac{\ga}{N}\right)^\frac{1}{\ga}$. Clearly, $\phi_{N,1}$ is
the usual two sided exponential density, considered in the previous
subsection. We show that,

\begin{prop}
If $\la_N$ has density $\phi_{N,\ga}$, then $\{\la_N\}$ satisfy LDP
with rate function $\mathcal{I}= \frac{|x|^\ga}{\ga}$ on $\mbb{R}$.
\end{prop}
\begin{proof}
 Let $\tri\subset \mbb{R}$ be an interval. Let $m =
\inf\{|x|:\, x \in \tri\}$, $M = \sup\{|x|:\, x \in \tri\}$, and $q_N=
\la_N(\tri)$. With this notation, we have
\begin{equation}
q_N =\dsty\int_\tri\frac{N}{2} |x|^{\ga-1} e^{-N\frac{|x|^\ga}{\ga}} dx
\leq\dsty\int_{N\frac{m^\ga}{\ga}}^{N\frac{M^\ga}{\ga}} e^{-x} dx<
\dsty\int_{N\frac{m^\ga}{\ga}}^{\infty}e^{-x}= e^{-N\frac{m^\ga}{\ga}} dx,
\end{equation}
 and
\begin{equation}
q_N \geq \frac{1}{2}\dsty\int_{N\frac{m^\ga}{\ga}}^{N\frac{M^\ga}{\ga}} e^{-x} dx
> \frac{1}{2}\dsty\int_{N\frac{m^\ga}{\ga}}^{N\frac{(m+\de)^\ga}{\ga}} e^{-x} dx
>\frac{\de}{2}N(m+\theta\de)^{\ga-1} e^{-N\frac{(m + \de)^\ga}{\ga}},
\end{equation}
for any $0<\de < M - m$ and some $\theta$, $0<\theta<1$. Mean value
theorem is used here.

The above two inequalities imply that for any interval $\tri\subset
\mbb{R}$,  $\dlim_{N\ra \infty} \frac{1}{N} \log \la_N(\tri) = -
\frac{m^\ga}{\ga}$. Thus Proposition \ref{p0.3.3} completes the
proof.
\end{proof}

\begin{rem}
In this case, rate function $\mathcal{I}(x)=\frac{|x|^\ga}{\ga}$ is
a continuous function with compact level sets. But $\mathcal{I}(x)$
is convex, only when $\ga\geq1$. Note that, for $0<\ga<1$,
$\mathcal{I}(x)$ is not convex but clearly quasi-convex and hence
Theorem~\ref{t1.2.2} is applicable.
\end{rem}

\begin{cor}
If $\xi_N(\si)$ are independent (over $\sigma$) having density
$\phi_{N,\ga}$, $f$ is a continuous function on $\mbb{R}$ and
$H_N(\si)=Nf(\xi_N(\si))$ then almost surely,
$$\dlim_{N\ra \infty} \frac{1}{N} \log Z_N(\be) = \log2\, -
\dinf_{|x|^\ga\leq \ga\log 2} \left(\be f(x) +
\frac{|x|^\ga}{\ga}\right).$$
\end{cor}

As earlier, if we take, $f(x)=x$, the Hamiltonian for $N$-particle
system $H_N$ is a two sided Weibull distribution with shape
parameter $\ga$ and scale parameter
$\ga^\frac{1}{\ga}N^\frac{\ga-1}{\ga}$. In this case the problem of
limiting free energy reduces to the calculation of $\dinf_{|x|^\ga
\leq \log 2} \left(\frac{|x|^\ga}{\ga} +\be x\right)$, that is,
$\dinf_{0\leq x \leq (\ga\log 2)^\frac{1}{\ga}}
\left(\frac{x^\ga}{\ga} -\be x\right)$.

For $\ga>1$, to calculate the above infimum, we imitate the Gaussian
case. Let us denote $g(x)=\frac{x^\ga}{\ga} -\be x$ so that
$g'(x)=x^{\ga-1}-\be$ and $g''(x)=(\ga-1)x^{\ga-2}\geq0$ on
$[0,(\ga\log 2)^\frac{1}{\ga}]$. So $g$ being twice differentiable
convex function, the infimum of $g$ will attain where $g'=0$. But
$g'$ will be 0 on $[0,(\ga\log 2)^\frac{1}{\ga}]$ only when
$\beta\leq (\ga\log2)^\frac{1}{\ga}$. In that case, the infimum will
occur at $x=\be^{\frac{1}{\ga-1}}$. When
$\be>(\ga\log2)^\frac{1}{\ga}$, then $g'<0$ on $[0,(\ga\log
2)^\frac{1}{\ga}]$ and hence infimum occur at $x=(\ga\log
2)^\frac{1}{\ga}$.

For $\ga\leq 1$, the function $g(x)=\frac{x^\ga}{\ga} -\be
x=x^\ga\left(\frac{1}{\ga}-\be x^{1-\ga}\right)$ is a product of two
functions. Here $x^\ga$ is a positive increasing function on
$[0,(\ga\log 2)^\frac{1}{\ga}]$. On the other hand,
$\frac{1}{\ga}-\be x^{1-\ga}$ is a decreasing function taking the value
$\frac{1}{\ga}$ at $0$. If this function always remains positive then
clearly the minimum of $g$ is $0$ attained at $x=0$. On the other
hand, if this function takes negative value in $[0,(\ga\log
2)^\frac{1}{\ga}]$ then the infimum of $g$ is attained at $(\ga\log
2)^\frac{1}{\ga}$. This situation occurs only when
$\frac{1}{\ga}-\be x^{1-\ga}=0$ for some $x$ in $[0,(\ga\log
2)^\frac{1}{\ga}]$. This happens only when  $\be \geq
\ga^{-\frac{1}{\ga}}(\log2)^{-\frac{1-\ga}{\ga}}$. Hence the infimum
of $g$ on $[0,(\ga\log 2)^\frac{1}{\ga}]$ is attained at $x=0$ for
$\be<\ga^{-\frac{1}{\ga}}(\log2)^{-\frac{1-\ga}{\ga}}$ and at
$x=(\ga\log 2)^\frac{1}{\ga}$ for $\be \geq
\ga^{-\frac{1}{\ga}}(\log2)^{-\frac{1-\ga}{\ga}}$.

We can combine the above arguments as
\begin{thm}
If $\{H_N(\si),\, \si\in \Si_N\}$ are independent having two sided
Weibull distribution with shape parameter $\ga>0$ and scale
parameter $\ga^\frac{1}{\ga}N^\frac{\ga-1}{\ga}$, then almost
surely,
$$\dlim_N \frac{1}{N}\log
Z_N(\be)=\begin{cases}\log2+\frac{\ga-1}{\ga}\be^{\frac{\ga}{\ga-1}}
&\text{ if $\be<(\ga\log2)^\frac{1}{\ga}$},\\
(\ga\log2)^\frac{1}{\ga}\be &\text{ if $\be\geq
(\ga\log2)^\frac{1}{\ga}$}\end{cases}$$ if $\ga>1$

and
$$\dlim_N \frac{1}{N}\log Z_N(\be)=
\begin{cases}\log2 &\text{ if $\be< \ga^{-\frac{1}{\ga}}(\log2)^{-\frac{1-\ga}{\ga}}$},\\
(\ga\log2)^\frac{1}{\ga}\be &\text{ if $\be\geq
\ga^{-\frac{1}{\ga}}(\log2)^{-\frac{1-\ga}{\ga}}$}
\end{cases}$$ if $\ga\leq 1$.
\end{thm}

\begin{rem}
It is easy to verify that, if $\la_N$ has density,
$$\phi_{N,\ga}(x) =
\mbox{Const.} e^{-N\frac{|x|^\ga}{\ga}} \quad -\infty<x<\infty,$$
more precisely,
$$\phi_{N,\ga}(x) =
\frac{1}{2\Gamma(\frac{1}{\ga})}\ga^{\frac{\ga-1}{\ga}}
N^\frac{1}{\ga} e^{-N\frac{|x|^{\ga}}{\ga}} \quad -\infty<x<\infty$$
then $\{\la_N\}$ satisfies LDP with rate function
$\mathcal{I}(x)=\frac{1}{\ga}|x|^\ga$. Note that here $\ga=2$ is the
Gaussian distribution. Hence in REM, if we consider $H_N$ to
$\mathcal{N}(0,N)$ or two sided Weibull distribution with shape
parameter $\ga=2$ and scale parameter $\sqrt{2N}$, they will produce
the same limiting free energy. So the limiting free energy of REM is
not entirely distribution specific, but it is {\em 'rate-specific'}.
\end{rem}

\bigskip

\section{Compact Distributions}

In the previous section we observed that, for the existence and
evaluation of free energy, we concentrated our attention on the set
$\{x:\,\mathcal{I}(x)\leq \log2\}$. That is, the entire support of
the random variables are not contributing to the system. To do that,
we used a variant of Varadhan's lemma. In general Varadhan's lemma
is applicable to the class of bounded continuous functions. In our
case the functions used in the previous section are rather
unbounded. As suggested by Proposition \ref{p0.3.5}, if the
underlying sequence of probabilities are eventually supported on a
compact set, we can overcome this little technicality. Our
assumption that the rate function $\mathcal{I}$ is a good rate
function will ensure that the sequences of induced probabilities are
almost surely eventually supported on a compact set. Since
$\mathcal{I}$ is a good rate function $\{\mathcal{I}(x)\leq \al\}$
is a compact set for every $\al\in\mbb{R}$. In particular,
$\{\mathcal{I}(x)\leq \log2\}$ is a compact set. For example, if
$\{\la_N\}$ is a sequence of probabilities satisfying LDP with rate
function $\mathcal{I}$ so that $\mathcal{I}(x)\leq \log2$ for all
$x\in \mbb{R}$, then we may not be able to apply Varadhan's lemma to
get the free energy of the system. In particular, if
$\mathcal{I}(x)=0$ for all $x\in \mbb{R}$, we can not infer anything
about the existence of the free energy of the system by large
deviation techniques.

To start with, let us note that, if $\mathcal{I}$ takes two value
$0$ and $\infty$ so that $\{\mathcal{I}(x)=0\}$ is a compact set,
say $C$. In this case, in view of Remark \ref{r1.2.2}, we can apply
Theorem \ref{t1.2.3} with $f(x)=\be x$. This will ensure the almost
sure existence of the limiting free energy and is equal to
$\log2-\dinf_{C} \be x=\log 2 - \be x_0$, where $x_0=\inf\{x: \,
x\in C\}$. Note that, if $C\subset [0,\infty)$ then clearly the
limiting free energy becomes negative for large $\be$ if $0\notin C$
whereas if $0\in C$ then it will be just a constant, $\log2$. So we
will not get any phase transition here.

Though we do not have a clear picture when $\mathcal{I}$ is
identically 0, we have some partial results. First of all, note that
{\em $\mathcal{I}(x)=0$ for all $x \in \mbb{R}$ iff
$\lim\limits_{N\ra \infty} \frac{1}{N} \log \la_N(\tri)=0$ for every
open subset $\tri$ of $\mbb{R}$.} This follows from
definition of LDP.

Now let us consider the case, when the Hamiltonian is supported on a
compact set. For each $N$, let $\la_N$ be a compactly supported
symmetric probability with density $\phi_N$ and $\{\xi_N(\si):\,
\si\in\Si_N\}$ be independent random variables having density
$\phi_N$. Consider the Hamiltonian $$H_N(\si)=N\xi_N(\si).$$ Let
$[-\al_N,\al_N]$ be the support of $\phi_N$. Let us assume that
$\al_N\ra \al$ as $N\ra \infty$. Here we allow the possibility that
$\al=\infty$. For $s\geq 0$, denote $a_N(s)= P\{\xi_N(\si) \geq
s\}$. Note that in this setup if $\{\la_N\}$ satisfy LDP with rate
function $\mathcal{I}(x)=0$, then $\frac{1}{N}\log a_N(s) \ra 0$ as
$N\ra \infty$ for $0\leq s<\al$.

The following theorem suggests that if the tail probability does not
decay exponentially fast over $N$, then we can not expect any annealed
phase transition.

\begin{thm}\label{t1.3.1}
Let $[-\al_N,\al_N]$ be the support of $\xi_N$ and for $s\geq 0$,
denote $a_N(s)= P\{\xi_N(\si) \geq s\}$. If $\al_N \ra \al$ and
$\frac{1}{N}\log a_N(s) \ra 0$ as $N\ra \infty$ for $0\leq s<\al$,
then
$$\dsty\lim_{N\ra \infty} \frac{1}{N}E\log Z_N(\be) = \log 2 +
\al\be.$$
\end{thm}

\begin{proof}
As
$$2^Ne^{-\be N\al_N}\leq Z_N(\be) \leq 2^Ne^{\be N\al_N},$$
the proof for $\al=0$ is immediate. Moreover, in this case, for
every sample point $$\dsty\lim_{N\ra \infty} \frac{1}{N}\log
Z_N(\be) = \log 2.$$

So let $\al>0$ (may be $\al=\infty$). Since $\log$ is concave, by
Jensen's inequality
\begin{equation}\label{e1.3.1}
E\log Z_N(\be) \leq \log E Z_N(\be).
\end{equation}

As $H_N$ are bounded by $N\al_N$, $$E Z_N(\be) = 2^N E e^{\be H_N}
< 2^N e^{\be N\al_N}.$$

Hence, by assumption and (\ref{e1.3.1}),
\begin{equation}\label{e1.3.2}
\dsty\limsup_{N \ra \infty}\frac{1}{N}E\log Z_N(\be) \leq \log 2 +
\al\be.
\end{equation}

Now we show,
\begin{equation}\label{e1.3.3}
\dsty\liminf_{N \ra \infty}\frac{1}{N}E\log Z_N(\be) \geq \log 2 +
\al\be.
\end{equation}
For that, with arbitrary but fixed $0\leq s<\al$, let $X_N= \#\{\si:
H_N(\si) \geq s N\}$. Then $E X_N = 2^N a_N(s)$ and $E X_N^2 =
2^N(2^N-1)a_N^2(s) + 2^N a_N(s).$ Hence
\begin{equation}\label{e1.3.4}
E(X_N-E X_N)^2 = E X_N^2 - (E X_N)^2 \leq 2^N a_N(s).
\end{equation}

If $X_N\leq 2^{N-1} a_N(s)$ then $E X_N-X_N\geq 2^{N-1} a_N(s)$ so
that $(X_N-E X_N)^2 \geq 2^{2N-2} a_N^2(s)$. Let $A_N = \{X_N \leq
2^{N-1} a_N(s)\}$. So $A_N \subset \{(X_N-E X_N)^2 \geq 2^{2N-2}
a_N^2(s)\}$. Hence, by Markov inequality and (~\ref{e1.3.4}),
$$P(A_N) \leq \frac{E(X_N - EX_N)^2}{2^{2N-2}a_N^2(s)} \leq
\frac{4}{2^Na_N(s)}.$$

i.e., $P(A_N^c) \geq 1- \frac{4}{2^Na_N(s)}$. But on $A_N^c$,
$$Z_N(\be) \geq X_N e^{\be s N} \geq 2^{N-1} a_N(s) e^{\be s N},$$ and hence
\begin{equation}\label{e1.3.5}
E\left[\log Z_N(\be) 1_{A_N^c}\right] \geq [(N-1)\log2 + \log a_N(s) + \be s
N]\left(1-\frac{4}{2^N a_N(s)}\right).
\end{equation}

Now $A_N= \{X_N=0\}\cup\{1\leq X_N \leq 2^{N-1} a_n(s)\}$. Since
$Z_N(\be) \geq 2^N e^{-\be N\al_N}$ and $PP(X_N=0)=(1-a_N(s))^{2^N}$
we have
\begin{equation}\label{e1.3.6}
E \left[\log Z_N(\be)1_{\{X_N=0\}}\right] \geq (N\log2 -\be
N\al_N)(1-a_N(s))^{2^N}.
\end{equation}

On $\{1 \leq X_N \leq 2^{N-1}a_N(s)\}$, $\log Z_N(\be) \geq \be
\dsty\max_\si H_N(\si) \geq \be s N >0$ and hence
\begin{equation}\label{e1.3.7}
E\left[\log Z_N(\be)1_{\{1 \leq X_N \leq 2^{N-1}a_N(s)\}}\right] \geq 0.
\end{equation}

Thus from (\ref{e1.3.5}), (\ref{e1.3.6}) and (\ref{e1.3.7}) we have
$$\begin{array}{rcl}
\frac{1}{N}E\log Z_N(\be)& \geq& \left[\frac{N-1}{N}\log2+
\frac{\log a_N(s)}{N} +\be s\right]\left(1-\frac{4}{2^Na_N(s)}\right)\vspace{1ex}\\
&& + \left[\log2 -\be \al_N\right](1-a_N(s))^{2^N}.
\end{array}$$

By assumption, $\frac{1}{N}\log a_N(s) \ra 0$ so that $2^Na_N(s) \ra
\infty$ and hence $(1-a_N(s))^{2^N} \ra 0$ as $N\ra \infty$. Thus,
under the assumption,
$$\dsty\liminf_{N\ra \infty}\frac{1}{N}E\log Z_N(\be)
\geq \log2 + \be s.$$

Since $0\leq s<\al$ is arbitrary, we have
$$\dsty\liminf_{N \ra \infty}\frac{1}{N}E\log Z_N(\be) \geq \log 2 +
\al\be$$ which remain true even when $\al=\infty$ with the
understanding that the right side of the above inequality is
$\infty$.

This completes the proof.
\end{proof}

Since $\xi_N$ has density $\phi_N$ with support $[-\al_N, \al_N]$
and $H_N=N\xi_N$, the support of $H_N$ will be $[-T_N, T_N]$ where
$T_N=N\al_N$. If we assume, $\varphi_N$ be the density of $H_N$ with
support $[-T_N, T_N]$, then we can apply the above theorem with
$\al_N=\frac{T_N}{N}$.

The following examples will illustrate the applications of the above
theorem.
\begin{ex}[Uniform Distribution]

Let $\varphi_N(x) =\frac{1}{2T_N}1_{[-T_N,T_N]}$. If
$\frac{T_N}{N}\ra \al>0$, then $a_N(s) \ra \frac{\al-s}{2\al} >0$
for all $0\leq s<\al$. So utilizing the above theorem we get,

a) if $T_N = \sqrt{N}$ then $\dsty\lim_{N\ra \infty} \frac{1}{N}\log
Z_N(\be) = \log 2$ for every sample point,

b) if $T_N = N$ then $\dsty\lim_{N\ra \infty} \frac{1}{N}E\log
Z_N(\be) = \log2 + \be$,

c) if $T_N = N^2$ then $\dsty\lim_{N\ra \infty} \frac{1}{N}E\log
Z_N(\be) = \infty$.

Similar remarks follows for the other examples also.
\end{ex}

\begin{ex}

Let $\de>0$ be fixed and $$\varphi_N(x) =\frac{\de +
1}{2T_N^{\de+1}}(T_N-|x|)^\de, \hspace{5ex}-T_N\leq x \leq T_N.$$
If $\frac{T_N}{N}\ra \al>0$, then $a_N(s) \ra
\frac{1}{2}\left(\frac{\al-s}{\al}\right)^{\de+1}
>0$ for all $0\leq s<\al$.
\end{ex}

\begin{ex}
Let $$\varphi_N(x) =\frac{1}{2T_N}\cos\frac{x}{T_N},
\hspace{5ex}-\frac{\pi T_N}{2}\leq x \leq \frac{\pi T_N}{2}.$$ If
$\frac{T_N}{N}\ra \al>0$, then $a_N(s) \ra
\frac{1}{2}\left(1-\sin\frac{s}{\al}\right)>0$ for all $0\leq
s<\al\frac{\pi}{2}$.
\end{ex}

\begin{ex}
Let $$\varphi_N(x)
=\frac{N}{2(e^{NT_N}-1)}e^{N|x|}1_{[-T_N,T_N]}.$$ If
$\frac{T_N}{N}\ra \al>0$, then $a_N(s) \ra \frac{1}{2}$ as $N\ra
\infty$ for all $0\leq s<\al$.
\end{ex}\vspace{1ex}

In the following examples, above
theorem is not applicable.
\begin{ex}[Truncated Double Exponential]
Let $$\varphi_N(x) =\frac{1}{2(1-e^{-T_N})}e^{-|x|}1_{[-T_N,T_N]}.$$
Let $\frac{T_N}{N}\ra \al (>0)$ as $N\ra \infty$. Then
$a_N(s)=\frac{e^{T_N-sN}-1}{2(e^{T_N}-1)}$. Hence $\frac{\log
a_N(s)}{N} \ra -s \neq 0$ as $N\ra \infty$ for all $s$ with
$0<s<\al$. Thus we can not apply the Theorem~\ref{t1.3.1} any more.
However, if $H_N(\si)$ has density $\varphi_{N}(x)$ and $\la_N$ is
the law of $\frac{1}{N}H_N(\si)$ then if $\frac{T_N}{N}\ra \al
(>0)$, by analysis of subsection~\ref{s1.3.2}, we can easily see
that the sequence $\{\la_N\}$ satisfies large deviation principle
with rate function $\mathcal{I}$ given by,
$$\mathcal{I}(x)=
\begin{cases}
|x| &\text{for $|x|\leq \al$}\\
\infty &\text{otherwise}.
\end{cases}$$
Hence we can use Theorem~\ref{t1.2.3} to conclude that the free
energy will be same as that of exponential REM as long as
$\al\geq\log2$. Where as if $\al<\log2$ then almost surely,
$$\lim_N\frac{1}{N}\log Z_N(\be)=\begin{cases}\log2 &\text{ for
$0\leq \be \leq 1$}\\ \log2-\al+\be\al &\text{for $\be\geq1$}.
\end{cases}$$
\end{ex}

\begin{ex}[Truncated Gaussian]
Let $$\varphi_N(x)
=\frac{1}{C_N}e^{-\frac{x^2}{2N}}1_{[-T_N,T_N]}.$$ Let
$\frac{T_N}{N}\ra \al (>0)$ as $N\ra \infty$. Then $\frac{\log
a_N(s)}{N} \ra -\frac{1}{2}s^2 \neq 0$ as $N\ra \infty$ for all $s$
with $0<\frac{1}{2}s^2<\al$. Thus we can not apply the
Theorem~\ref{t1.3.1} once again. However, if $H_N(\si)$ has density
$\varphi_{N}(x)$ and $\la_N$ is the law of $\frac{1}{N}H_N(\si)$
then if $\frac{T_N}{N}\ra \al (>0)$, by analysis of
subsection~\ref{s1.3.1}, we can easily see that the sequence
$\{\la_N\}$ satisfy large deviation principle with rate function
$\mathcal{I}$ given by,
$$\mathcal{I}(x)=
\begin{cases}
\frac{1}{2}x^2 &\text{for $\frac{1}{2}x^2\leq \al$}\\
\infty &\text{otherwise}.
\end{cases}$$
Hence we can use Theorem~\ref{t1.2.3} to conclude that the free
energy will be same as that of Gaussian REM as long as
$\al\geq\log2$. Where as if $\al<\log2$ then almost surely,
$$\lim_N\frac{1}{N}\log Z_N(\be)=\begin{cases}\log2+\frac{1}{2}\be^2
&\text{ for $0\leq \be \leq \sqrt{2\al}$}\\ \log2-\al+\be\sqrt{2\al}
&\text{for $\be\geq \sqrt{2\al}$}. \end{cases}$$
\end{ex}

\bigskip
\section{Discrete Distributions}
We conclude this chapter by considering the REM driven by some
discrete distributions.
\subsection{Poisson Distribution}
Let us consider the Hamiltonian for the $N$ particle system
$H_N(\si)\thicksim P(N\theta)$ where $P(N\theta)$ is the Poisson
distribution with parameter $N\theta$. Let $\la_N$ be the law
$\frac{1}{N}P(N\theta)$. We also can think of $\la_N$ as the law of
the sample mean for a sample of size $N$ from $P(\theta)$. Then by
Cramer's theorem (Theorem \ref{p0.3.6}), $\{\la_N\}$ satisfies LDP
with convex good rate function $\mathcal{I}$ given by
\begin{equation}
\mathcal{I}(x)=\left\{\begin{array}{ll}\theta-x+x\log\frac{x}{\theta}
& \mbox{for } x\geq0\\ \infty &\mbox{otherwise}\end{array}\right..
\end{equation}
Hence by Theorem \ref{t1.2.3}, if for $\si\in\Si_N$ the random
variables $\xi(\si)$ is distributed like $\la_N$ and the Hamiltonian
$H_N(\si)=N f(\xi(\si))$ with any continuous function $f$ on
$\mbb{R}$, we have the following:

\begin{cor}\label{c1.5.1}
If $\la_N\thicksim \frac{1}{N}P(N\theta)$, then almost surely,
$$\dlim_{N\ra \infty} \frac{1}{N} \log Z_N(\be) = \log2\, -
\dinf_{\mathcal{I}(x)\leq \log 2} \left(\be f(x) +
\mathcal{I}(x)\right).$$
\end{cor}
\noindent{\em Notation:} Note that here $\mathcal{I}$ is a convex
continuous function on $[0,\infty)$ so that $\mathcal{I}(0)=\theta$;
$\mathcal{I}(\theta)=0$ and $\mathcal{I}(x)\ra\infty$ as
$x\ra\infty$. So the set $\{x:\mathcal{I}(x)=\log2\}$ contains only
one point when $\theta<\log2$; contains zero and one non-zero-point
for $\theta=\log2$; contains two positive points for $\theta>\log2$.
As a consequence, the set $\{\mathcal{I}(x)\leq\log2\}$ is an
interval $[x_1,x_2]$; $x_1=0$ in case of $\theta\leq\log2$ where as
$x_1>0$ in case of $\theta>\log2$. In any case, $\theta\in
(x_1,x_2)$.

Hence when $f(x)=x$ the above corollary implies that
$$\begin{array}{ll}
\dlim_N\frac{1}{N}\log Z_N(\be)&=\log2-\dinf_{[x_1,x_2]}\left\{\be
x+\theta-x+x\log\frac{x}{\theta}\right\}\\
&=\log2-\theta-\dinf_{[x_1,x_2]}\left\{(\be-1)x+x\log\frac{x}{\theta}\right\}.
\end{array}$$
To calculate the above infimum, let
$g(x)=(\be-1)x+x\log\frac{x}{\theta}$ on $[0,\infty)$. Clearly, $g$
is a convex function. $g'(x)=\be+\log\frac{x}{\theta}$, so that
$g'(\ol{x})=0$ implies $\ol{x}=\theta e^{-\be}$. Hence $g$ attains
its infimum at $\ol{x}=\theta e^{-\be}$. We consider two cases,
namely, $\theta\leq\log2$ and $\theta>\log2$.\\

\noindent\doublebox{ $\boldsymbol{\theta\leq\log2}$}\vspace{1ex}

For $\be\geq0$, $0<\theta e^{-\be}\leq \theta<x_2$. That is, the
point $\ol{x}=\theta e^{-\be}$, where $g$ attains minimum, belongs
to $\in[x_1,x_2]$ for every $\be\geq0$. Hence
$\dinf_{[x_1,x_2]}g(x)=(\be-1)\theta e^{-\be}-\theta e^{-\be}\log
e^{-\be}=-\theta e^{-\be}$\\

\noindent\doublebox{ $\boldsymbol{\theta>\log2}$}\vspace{1ex}

As $\be$ increases from $0$ to $\infty$, $\theta e^{-\be}$ decreases
from $\theta$ to $0$. Since $0<x_1<\theta$, there exists $\be_0>0$
such that
$$\theta e^{-\be_0}=x_1.$$ Clearly, for $\be\leq \be_0$, $\theta
e^{-\be}\in[x_1,x_2]$ so that $\dinf_{[x_1,x_2]}g(x)=-\theta
e^{-\be}$. Since $g$ attains its infimum at $\ol{x}=\theta
e^{-\be}$, $g$ is increasing (by convexity) on $(\ol{x},\infty)$.
For $\be>\be_0$, $\ol{x}=\theta e^{-\be}<\theta e^{-\be_0}=x_1$.
Thus $g$ is increasing on $[x_1, x_2]$. As a consequence, when
$\be>\be_0$, we have $\dinf_{[x_1,x_2]}g(x)=g(x_1)=\be
x_1+\mathcal{I}(x_1)-\theta=\be x_1+\log2-\theta$.

All this leads to
\begin{thm}
Consider REM where the Hamiltonian $H_N(\si)$ is  Poisson with
parameter $N\theta$.

a) For $\theta\leq\log2$, almost surely,
$$\dlim\frac{1}{N}\log Z_N=\log2-\theta+\theta e^{-\be}\quad\mbox{for } \be\geq0.$$

b) For $\theta>\log2$; let $x_1$ be the least positive solution of
$x(\log\frac{x}{\theta}-1)=\theta-\log2$, and
$\be_0=\log\frac{\theta}{x_1}=\frac{\theta-\log2}{x_1}-1$. Then
almost surely,
$$\begin{array}{lll}
\dlim\frac{1}{N}\log Z_N&=\log2-\theta+\theta e^{-\be}&\mbox{for }
\be\leq\be_0\\&=\be x_1&\mbox{for } \be>\be_0.\end{array}$$
\end{thm}

Now if we take $f(x)=-x$, then by Corollary \ref{c1.5.1}, almost
surely, the limiting free energy is given by
\begin{align*}
\dlim_N\frac{1}{N}\log Z_N(\be)&=\log2-\dinf_{\mathcal{I}(x)\leq\log 2}\{\mathcal{I}(x)-\be x\}\\
&=\log2-\theta-\dinf_{\mathcal{I}(x)\leq \log 2}\left\{x\log\frac{x}{\theta}-(\be+1)x\right\}.
\end{align*}

To calculate the above infimum, let
$g(x)=x\log\frac{x}{\theta}-(\be+1)x$ on $[0,\infty)$ so that
$g'(x)=\log\frac{x}{\theta}-\be$ and $g''(x)=\frac{1}{x}>0$ on
$(0,\infty)$. Hence $g$ attains its infimum at $\ul{x}=\theta
e^{\be}$. Note that, $\ul{x}=\theta$ for $\be=0$ and
$\ul{x}\ra\infty$ as $\be\ra\infty$. So there exist $\be_1>0$ such
that $\mathcal{I}(\theta e^{\be_1})=\log2$, that is, $\theta
e^{\be_1}=x_2$. So the infimum $\dinf_{\mathcal{I}(x)\leq \log
2}\left\{x\log\frac{x}{\theta}-(\be+1)x\right\}$ occurs at $\theta
e^{\be}$ for $\be\leq\be_1$ and at $x_2$ for $\theta>\be_1$ leading
to the following

\begin{thm}
In REM, if the Hamiltonian $H_N(\si)$ negative of Poisson with
parameter $N\theta$, then almost surely
$$\begin{array}{lll}
\dlim\frac{1}{N}\log Z_N&=\log2-\theta+\theta e^{\be}&\mbox{for }
\be\leq\be_1\\&=\be x_2&\mbox{for } \be>\be_1.\end{array}$$
\end{thm}

\subsection{Binomial Distribution}
Let $X_N\thicksim B(N,p)$ where $B(N,p)$ is the Binomial
distribution with parameter $p\,(0<p<1)$. Put $\xi_N=\frac{X_N}{N}$.
Observe that when $p=0$ or $1$, then the Hamiltonian is
deterministic one and uninteresting. Let $\la_N$ be the law of
$\xi_N$. Thus $\xi_N$ is nothing but the proportion of heads in $N$
tosses of a coin (with chance of heads $p$). By Cramer's theorem
(Theorem \ref{p0.3.6}), $\{\la_N\}$ satisfies LDP with convex good
rate function $\mathcal{I}$ given by
\begin{equation}
\mathcal{I}(x)=\begin{cases}x\log\frac{x}{p}+(1-x)\log\frac{1-x}{1-p}
& \mbox{for } x\in[0,1]\\ \infty &\mbox{otherwise}.\end{cases}
\end{equation}
Note that here $\mathcal{I}$ is a strictly convex continuous
function. Now fix a continuous function $f$ on $\mbb{R}$. Consider
$N$ particle system with Hamiltonian $H_N=Nf(\xi_N(\si))$. By
Theorem \ref{t1.2.3}, to calculate the limiting free energy we only
have to solve the optimization problem
$$\dinf_{\mathcal{I}(x)\leq \log 2} \left(\be f(x) +
\mathcal{I}(x)\right).$$

Note that here $\mathcal{I}(0)=-\log(1-p)$; $\mathcal{I}(p)=0$ and
$\mathcal{I}(1)=-\log p$. So the set
$\{x:\mathcal{I}(x)=\log2\}=\{0,1\}$ when $p=\frac{1}{2}$ otherwise
the set $\{x:\mathcal{I}(x)=\log2\}$ is a singleton. Let us denote
the set $\{\mathcal{I}(x)\leq\log2\}$ as $[x_1,x_2]$ where
$0<x_1<p<x_2=1$ for $p>\frac{1}{2}$; $0=x_1<p<x_2<1$ for
$p<\frac{1}{2}$ and $0=x_1<p<x_2=1$ for $p=\frac{1}{2}$.

When $f(x)=x$, by Theorem \ref{t1.2.3}, we have almost surely,
\begin{align*} \dlim_N\frac{1}{N}\log
Z_N(\be)&=\log2-\dinf_{[x_1,x_2]}\left\{\be
x+x\log\frac{x}{p}+(1-x)\log\frac{1-x}{1-p}\right\}.
\end{align*}
To calculate the above infimum, let $g(x)=\be
x+x\log\frac{x}{p}+(1-x)\log\frac{1-x}{1-p}$ on $[0,1]$. Clearly $g$
is a convex function and $g'(x)=\be+\log\frac{x(1-p)}{(1-x)p}$. So
$g$ attains its infimum at $\ol{x}(\be)$ given by
$\ol{x}(\be)=\dfrac{p}{p+(1+p)e^\be}$. We now consider two cases.\\

\noindent\doublebox{ $\boldsymbol{p\leq\frac{1}{2}}$ }\vspace{1ex}

As $[x_1,x_2]=[0,x_2]$ where $x_2>p$ and $\ol{x}(\be)\leq
\frac{p}{2p+1}<p$ for every $\be\geq0$, we have
$\ol{x}(\be)\in[x_1,x_2]=[0,x_2]$. Hence on $[x_1,x_2]$, $g$ attains
its infimum at $\ol{x}(\be)$ and by routine algebraic manipulations,
we get, $\dinf_{[x_1,x_2]}g(x)=\be
\ol{x}(\be)+\mathcal{I}(\ol{x}(\be))=\be-\log(p+(1-p)e^\be)$.\\

\noindent\doublebox{ $\boldsymbol{p>\frac{1}{2}}$ }\vspace{1ex}

Since $\ol{x}(\be)$ decreases from $\frac{p}{2p+1}<p$ to $0$ as
$\be$ increases from $0$ to $\infty$ and $0<x_1<p$, there exists
$\be_0>0$ such that
$$\ol{x}(\be_0)=x_1.$$
Hence as $\ol{x}(\be)\in[x_1,x_2]$ for $\be\leq \be_0$, $g$ attains
its infimum at $\ol{x}(\be)$ on $[x_1,x_2]$ and we get
$\dinf_{[x_1,x_2]}g(x)=\be-\log(p+(1-p)e^\be)$. On the other hand,
for $\be>\be_0$, $g$ being a convex function and attains it infimum
at $\ol{x}<x_1$, $g$ is increasing for $x>\ol{x}$. Hence $g$ attains
its infimum on $[x_1,x_2]$ at $x_1$ leading to
$\dinf_{[x_1,x_2]}g(x)=g(x_1)=\be x_1+\mathcal{I}(x_1)=\be
x_1+\log2$.

All this leads to
\begin{thm}
In REM, if the Hamiltonian $H_N(\si)$ is Binomial with parameter $N$
and $p$, then almost surely,

a) for $p\leq\frac{1}{2}$,
$$\dlim\frac{1}{N}\log Z_N=\log2-\be+\log(p+(1-p)e^\be)\quad\mbox{for } \be\leq0.$$

b)for $p>\frac{1}{2}$,
$$\begin{array}{lll}
\dlim\frac{1}{N}\log Z_N&=\log2-\be+\log(p+(1-p)e^\be)&\mbox{for }
\be\leq\be_0\\&=-\be x_1&\mbox{for } \be>\be_0.\end{array}$$
\end{thm}

On the other hand, if we take $f(x)=-x$, then by Theorem
\ref{t1.2.3}, almost surely,
\begin{align*}
\dlim_N\frac{1}{N}\log Z_N(\be)&=\log2-\dinf_{\mathcal{I}(x)\leq\log 2}\{\mathcal{I}(x)-\be x\}\\
&=\log2-\dinf_{[x_1,x_2]}\left\{x\log\frac{x}{p}+(1-x)\log\frac{1-x}{1-p}-\be
x\right\},
\end{align*}
where we use the same notation for $x_1, x_2$ as in the case for
$f(x)=x$.

To calculate the above infimum, let
$h(x)=x\log\frac{x}{p}+(1-x)\log\frac{1-x}{1-p}-\be x$ on $[0,1]$ so
that $h'(x)=\log\frac{x(1-p)}{(1-x)p}-\be$ and
$g''(x)=\frac{1}{x(1-x)}>0$ on $(0,1)$. Hence $g$ attains its
infimum at $\ul{x}(\be)=\dfrac{p}{p+(1-p)e^{-\be}}$. Note that,
$\ul{x}(0)=p$ and $\ul{x}(\be)\ra 1$ as $\be\ra\infty$. So the
infimum $\dinf_{[x_1,x_2]}h(x)$ is attained at $\ul{x}(\be)$ for
every $\be>0$ for $p\geq\frac{1}{2}$ and for $\be\leq \be_1$ for
$p<\frac{1}{2}$, where $\be_1>0$ is such that
$\mathcal{I}(\ul{x}(\be_1))=\log2$, that is, $\ul{x}(\be_1)=x_2$.
For $p<\frac{1}{2}$ and $\be> \be_1$ the infimum
$\dinf_{[x_1,x_2]}h(x)$ is attained at at $x_2$.  To be more
precise, we get the following:

\begin{thm}
In REM, if the Hamiltonian $H_N(\si)$ is negative Binomial random
variable with parameter $N$ and $p$, then almost surely

a) for $p\geq\frac{1}{2}$,
$$\dlim\frac{1}{N}\log Z_N=\log2+\be+\log(p+(1-p)e^{-\be})\quad\mbox{for } \be\geq0.$$

b)for $p<\frac{1}{2}$,
$$\begin{array}{lll}
\dlim\frac{1}{N}\log Z_N&=\log2+\be+\log(p+(1-p)e^{-\be})&\mbox{for
} \be\leq\be_1\\&=\be x_2&\mbox{for } \be>\be_1.\end{array}$$
\end{thm}

\cleardoublepage

\chapter{The Generalized Random Energy Model}


In the random energy model (REM) {\cite{De80,D1}} of Derrida, the
Hamiltonians in distinct configurations are independent. The idea in
generalized random energy model (GREM) is to bring an amount of
dependence in the Hamiltonians through the structure of
configurations. Of course, very little can be achieved by assuming
an arbitrary covariance matrix. To introduce hierarchy, an $n$-level
tree structure was suggested by Derrida {\cite{D2}}, where the
branches of the tree are in correspondence with the configuration
space. In this chapter we discuss this model with some
modifications. There are two essential differences from what is
usually considered in the literature. First, we provide a general
framework of trees. However, they will be considered in the next
chapter. Second, we split the number of particles $N$ into $n$
groups rather than splitting the number of spins (or `factorizing' 2
as is customary in the literature). This allows us to introduce a
further randomization at the tree level, like Poisson trees and
multinomial trees. These will be consider in chapter 3.

In this chapter, we specialize to the driving distributions having
exponentially decaying tails. The basic inequalities lead to the
large deviation principle (LDP) for the random probabilities as in
the case of REM considered in the previous chapter. This leads to an
explicit formula for the free energy. For the exponential GREM, the
driving distribution does not depend on the number of particles.
This does not make it less interesting. In fact, the Gaussian case
is no more complicated than the exponential case. The present
treatment clearly brings out the similarities between the two cases.
There are dissimilarities too. As expected, for small values of
$\be$ (inverse temperature), the energy function in the exponential
case does not depend on $\be$ where as for the Gaussian it is
quadratic in $\be$. In the Gaussian case, all the weights associated
with all the levels of the tree participate in the expression for
free energy, where as in the exponential case it is not always so.

Even though for any finite number of particles, we have a truly $n$
level tree, in the limit, it may collapse to a lower level tree --
it may even correspond to REM (see remarks \ref{r2.6.1} and
\ref{r2.6.4}). This leads to the notion of reduced GREM. For such
models, the energy function determines all the parameters of the
model. It is also possible to characterize the energy functions. It
is interesting to note that in the SK-model, subject to certain
moment conditions of the underlying distribution, the energy
function is universal {\cite{CH}}, while it is not true here.
\smallskip
\goodbreak
\section{Derrida's Model}

Let us first describe the model in detail.  As a generalization to
his REM\cite{D1}, in GREM\cite{D2} Derrida introduced a tree-like
structure in the energy levels. This is what we now explain. Fix a
positive integer $n\geq 1$. This $n$ will be the level of the tree.
For each level $i= 1, 2, \cdots,n$ of the tree,  fix number $\al_i$
so that $\al_i\in (1,2)$ and $\prod_{i=1}^n \al_i=2$. For fixed $N$,
in the tree, there will be $\al_1^N$ many nodes at the first level.
Below each of the first level nodes, there will be $\al_2^N$ many
nodes in the 2nd level. Hence, there will be a total of
$(\al_1\al_2)^N$ many nodes at the 2nd level. In general, at the
$i$th level, there will be $\al_i^N$ many nodes below each of the
$(i-1)$th level nodes giving a total $(\al_1\al_2\cdots\al_i)^N$
many nodes in the $i$th level. So at the $n$-th, that is, last level
there will be $(\al_1\al_2\cdots\al_n)^N=2^N$ many nodes (leaves).
Derrida associates the configuration space $\Si_N$ with the all
possible branches from root to leaves of the above tree. Since there
are $2^N$ many configurations, he assumes $\prod_{i=1}^n \al_i=2$.
To define the Hamiltonian, he associates an independent random
variable to each edge of the tree.  For $i= 1, 2, \cdots,n$ there
are $(\al_1\al_2\cdots\al_i)^N$ independent Gaussian mean zero
random variables $\xi_j^{(i)}$ with variance $a_i N$ associated to
each of the $i$th level edges. Here $a_1, a_2, \cdots, a_n$ are
positive numbers so that $\sum_{i=1}^n a_i=1$. The Hamiltonian for a
configuration, that is, for a branch from root to a leaf is the sum
of the $n$ random variables associated with the edges constituting
the branch. So the partition function, in this model, reduces to
$$Z_N(\be)=\sum_{i_1=1}^{\al_1^N}\sum_{i_2=(i_1-1)\al_2^N+1}^{i_1\al_2^N}
\cdots\sum_{i_2=(i_{n-1}-1)\al_n^N+1}^{i_{n-1}\al_n^N}
e^{-\be\left(\sum_{k=1}^n \xi_{i_k}^{(k)}\right)}.$$ In the entire
explanation above, we pretended that each $\al_i^N$ is an integer.
But is this possible? -- No. One way out is to consider $[\al_i^N]$.
Being the number of edges, each $\alpha_i^N$ has to be an integer
which divides $2^N$ because $(\al_1\al_2\cdots\al_n)^N=2^N$. By the
fundamental theorem of arithmetics, $\alpha_{i}^{N}=2^{k(i,N)}$ for
some positive integer $k(i,N)$. Moreover, $k(i,N)$ for
$i=1,\cdots,n$ is such that $k(1,N)+ k(2,N) +\cdots+ k(n,N) =N$. In
other words: given  any tree with $2^N$ leaves the construction
allows only for furcations in powers of 2 at each layer. This was
also noted in~\cite{CEGG}.  To eliminate the confusion regarding
whether $\al_i^N$ is an integer or not, we made the natural
modification to the model in the next section.

\goodbreak

\section{A Reformulation}

We formulate GREM as follows. As above, fix an integer $n\geq 1$.
Let $N\geq n$ be the number of particles, each of which can have two
states/spins $+1, -1$; so that the configuration space is $2^N$.
Consider a partition of $N$, into integers $k(i,N)$ for $1\leq i
\leq n$ with each $k(i,N) \geq 1$ and $\dsum_i k(i,N)=N$. The
configuration space $2^N$, naturally splits into product, $\prod
2^{k(i,N)}$ and $\si \in 2^N$ can be written as
$\si_1\si_2\cdots\si_n$ with $\si_i\in 2^{k(i,N)}$. An obvious tree
structure can be brought in the configuration space. As earlier
imagine an $n$-level tree. There are $2^{k(1,N)}$ nodes at the first
level. These will be denoted as $\si_1$, with $\si_1\in2^{k(1,N)}$.
Below each of the first level nodes there are $2^{k(2,N)}$ nodes at
the second level. The second level nodes below $\si_1$ of the first
level will be denoted by $\si_1 \si_2$ with $\si_2 \in 2^{k(2,N)}$.
In general, below a node $\si_1\si_2\cdots\si_{i-1}$ of the
$(i-1)$-th level there are $2^{k(i,N)}$ nodes at the $i$-th level
denoted by $\si_1\si_2\cdots\si_{i-1}\si_i$ for $\si_i\in
2^{k(i,N)}$. Thus a typical branch of the tree reads like
$\si_1\si_2\cdots\si_n$. Obviously the branches are in one one
correspondence with $2^N$, the configuration space. At the node
$\si_1\cdots\si_i$, we place a random variables
$\xi(\si_1\cdots\si_i)$. We assume that all these random variables
are i.i.d. with a symmetric distribution. We associate one weight
for each level, say weight $a_i>0$ for the $i$-th level. These are
not random. In a configuration $\si=\si_1\cdots\si_n$ the
Hamiltonian is
\begin{equation}\label{e2.2.1}
H_N(\si) = \dsum_{i=1}^{n} a_i\xi(\si_1\cdots\si_i).
\end{equation}
For $\be>0$ the partition function is
\begin{equation}\label{e2.2.2}
Z_N(\be) =\dsum_\si e^{-\be H_N(\si)}= 2^N \mbf{E}_\si e^{-\be
H_N(\si)}.
\end{equation}
Here $\mbf{E_\si}$ stands for expectation with respect to $\si$ when
the configuration space $2^ N$ has uniform distribution. In other
words, $\mbf{E_\si}$ is simply the usual average over $\si$.

Since $\xi$'s are random variables both $H_N$ and $Z_N$ are random
variables. We suppress the parameter $\om$. As usual
$\frac{1}{N}\log Z_N(\be)$ is the free energy of the $N$-particle
system. This is the object of study. As $N$ changes, the common
distribution of the $\xi$'s would in general change and so in $H_N$.

\goodbreak

\section{Tree Formulation}\label{s2.3}
We now reformulate the setup as a general tree structure. Though
most of the trees that we consider later are {\em regular} $--$ in
the sense that the number of furcations of a node depend only on its
level, and not on the particular node $--$  the present formulation
is general. It allows randomization of the tree, which we do
consider later in the next chapter.

Let $n\geq 1$ be fixed integer as earlier. For each $N\geq n$, let
$T_N$ be a tree of height $n$ with each branch extending up to the
$n$-th level. $\si_1$ denotes a typical node at the first level and
in general below a node $\si_1\si_2\cdots\si_{i-1}$ of the
$(i-1)$-th level, $\si_1\si_2\cdots\si_{i-1}\si_i$ is a typical node
at the $i$-th level. We shall now define some useful quantities
associated with the tree. Let $\Si_N$ be the set of all branches
$\si_1\si_2\cdots\si_n$ of the tree $T_N$. Let $B_{iN}$ denote the
number of nodes at the $i$-th level. In particular, $B_{nN}$ is the
total number of branches of the tree, which will simply be denoted
by $B_N$. For a node $\si_1\si_2\cdots\si_i$ of the $i$-th level,
let $e(\si_1\si_2\cdots\si_i)$ denote the number of nodes at the
$n$-th level below the node $\si_1\si_2\cdots\si_i$. Equivalently,
$e(\si_1\si_2\cdots\si_i)$ is the total number of branches extending
$\si_1\si_2\cdots\si_i$. Clearly,
$\sum\limits_{\si_1,\cdots,\si_i}e(\si_1\cdots\si_i)= B_N$ for any
$i$. Let
$s_{iN}^2=\sum\limits_{\si_1,\cdots,\si_i}e^2(\si_1\cdots\si_i)$.

Assume that  $\xi(\si_1\cdots\si_i)$ is a symmetric random variable
associated with node $\si_1\si_2\cdots\si_i$. We assume that these
random variables are i.i.d. Strictly speaking we should be using
superscript $N$ for the nodes, random variables etc. But for ease in
reading we suppress the superscript. This should be borne in mind.
We do assume that all our random variables are defined on one
probability space. Consider the map $\Si_N \rightarrow \mathbb{R}^n$
defined by
\[\si \mapsto \xi_\si = (\xi(\si_1), \xi(\si_1\si_2), \cdots,
\xi(\si_1\cdots\si_n)).\]

Let $\mu_N$ be the induced probability on $\mathbb{R}^n$ when
$\Si_N$ has uniform distribution, that is, each $\si \in \Si_N$ has
probability $\frac{1}{B_N}$. In other words, for any Borel set $A
\subset \mathbb{R}^n$,
\[\mu_N(A)=\frac{1}{B_N} \#\{\si:\xi_\si \in A\}.\]

In particular, if $A$ is a box, say
$\tri=\tri_1\times\cdots\times\tri_n$, with each
$\tri_i\subseteq\mathbb{R}$, then
\[\mu_N(\tri) = \frac{1}{B_N} \dsum_{<\si_1\cdots\si_n>}
\prod\limits_{i=1}^{n}
\mbf{1}_{\tri_i}(\xi(\si_1\si_2\cdots\si_i)).\]

Denote $q_{iN} = P(\xi \in \tri_i)$ for $1\leq i \leq n$. Since all
the $\xi(\si_1\cdots\si_i)$ (for fixed $N$) are i.i.d., we did not
use suffix for $\xi$ in defining $q_{iN}$. However since the common
distribution will in general change with $N$, $q_{iN}$ would in
general depend on $N$. Then
\begin{equation}\label{e2.3.0}
\mathbf{E}\mu_N(\triangle)=q_{1N}q_{2N}\cdots q_{nN}.
\end{equation}

Here now is the basic result.
\begin{thm}\label{t2.3.3}
Let $\tri = \tri_1 \times \cdots \times \tri_n \subset
\mathbb{R}^n$. Denote $q_{iN} = P(\xi \in \tri_i)$ for $1\leq i \leq
n$.

a) If  $\sum\limits_{N \geq n}B_{iN} q_{1N}\cdots q_{iN} < \infty$,
for some $i, 1\leq i \leq n$ then a.s. eventually, $$\mu_N(\tri)= 0.$$

b) If for all $i=1,\cdots,n$, $\sum\limits_{N\geq n}
\frac{s_{iN}^2}{B_N^2 q_{1N}\cdots q_{iN}} < \infty$, then for any
$\ep >0$ a.s. eventually,
\[(1-\ep)\mbf{E}\mu_N(\tri) \leq \mu_N(\tri) \leq (1+\ep)
\mbf{E}\mu_N(\tri).\]
\end{thm}

In proving the first part of the theorem we will use the idea of
Dorlas and Dukes \cite{DD}, where as for the last part, we follow
Capocaccia {\em et al} \cite{CCP}.
\begin{proof}
a)  Let $j_0$ be such that
$\sum\limits_{N\geq 1}B_{j_0N}q_{1N}\cdots q_{j_0 N} < \infty.$ Then
\[\begin{array}{lll}
\mu_N(\tri) &=& \dfrac{1}{B_N}\sum\limits_{\si_1 \cdots \si_n}
\prod\limits_{i=1}^n \mbf{1}_{\tri_i}(\xi(\si_1 \cdots \si_i))\\
&\leq&\dfrac{1}{B_N}\sum\limits_{\si_1 \cdots \si_{j_0}}
\prod\limits_{i=1}^{j_0} \mbf{1}_{\tri_i}(\xi(\si_1 \cdots \si_i))e(\si_1 \cdots \si_{j_0})\\
&=& G_N, (\mbox{say}).
  \end{array}\]

Let $A_N$ be the event $\{G_N=0\}.$ Observe that $$A_N^c=
\left\{\sum\limits_{\si_1 \cdots \si_{j_0}} \prod\limits_{i=1}^{j_0}
\mbf{1}_{\tri_i}(\xi(\si_1 \cdots \si_i))\geq 1\right\}.$$ Now by
Chebyshev's inequality,
$$\mbf{P}(A_N^c) \leq\mbf{E}\sum\limits_{\si_1 \cdots \si_{j_0}}
\prod\limits_{i=1}^{j_0} \mbf{1}_{\tri_i}(\xi(\si_1 \cdots \si_i))=B_{j_0 N}q_{1N} \cdots q_{j_0 N}.$$
Thus by assumption and Borel-Cantelli,  $A_N$ will occur a.s.
eventually. i.e. $G_N = 0$ and hence $\mu_N(\tri) = 0$.

b) We first get an estimate for the variance of $\mu_N(\tri)$.
\[\begin{array}{ll}
&var(\mu_N(\tri))\\
=& \mbf{E}(\mu_N(\tri))^2 - (\mbf{E}\mu_N(\tri))^2\\
=& \dfrac{1}{B_N^2}\sum\limits_{\substack{\si_1\cdots \si_n\\
 \tau_1 \cdots \tau_n }} \left[{\mbf{E}}\prod\limits_{i=1}^n
 {\mbf{1}}_{\tri_i}(\xi(\si_1 \cdots \si_i)){\mbf{1}}_{\tri_i}(\xi(\tau_1 \cdots \tau_i))-
 q_{1N}^2 \cdots q_{nN}^2\right]\\
\leq& \dfrac{1}{B_N^2}\sum\limits_{j=1}^{n} \sum\limits_{\si_1
\cdots
\si_j}\sum\limits_{\substack{\si_{j+1}\cdots \si_n\\
\tau_{j+1} \cdots \tau_n\\ \si_{j+1}\neq \tau_{j+1}}} \mbf{E}\prod\limits_{i=1}^j\mbf{1}_{\tri_i}(\xi(\si_1
\cdots \si_i))
\prod\limits_{i=j+1}^n \mbf{1}_{\tri_i}(\xi(\si_1 \cdots \si_i))\mbf{1}_{\tri_i}(\xi(\tau_1 \cdots \tau_i))\\
\leq&\dfrac{1}{B_N^2} \sum\limits_{j=1}^{n} q_{1N}\cdots
q_{jN}q_{(j+1)N}^2\cdots q_{nN}^2\sum\limits_{\si_1 \cdots
\si_j}e^2(\si_1 \cdots \si_j)\\
=&\dfrac{1}{B_N^2} \sum\limits_{j=1}^{n} q_{1N}\cdots
q_{jN}q_{(j+1)N}^2\cdots q_{nN}^2 s_{jN}^2
\end{array}\]

Hence for any $\ep >0$, by Chebyshev's inequality and (\ref{e2.3.0})
\[\mbf{P}(|\mu_N(\tri)-\mbf{E}\mu_N(\tri)|>\ep \mbf{E}\mu_N(\tri))<
\frac{1}{\ep^2 B_N^2}\sum\limits_{j=1}^n\frac{s_{jN}^2}{ q_{1N}
\cdots q_{jN}}.\] But, in view of the assumption, the sum over N of
the right side is finite. So by Borel-Cantelli lemma, a.s.
eventually,
$$(1-\ep)\mbf{E}\mu_N(\tri) \leq \mu_N(\tri) \leq (1+\ep) \mbf{E}\mu_N(\tri).$$
\end{proof}

For GREM type regular trees the condition above will simplify as
follows. This result is in ~\cite{DD} though not explicitly stated.
\begin{cor}\label{c2.3.4}
Let $k(i,N)$, $1\leq i \leq n$ be positive integers with
$\sum\limits_i k(i,N)=N$. Suppose that the tree has $2^{k(i,N)}$
nodes of the i-th level below each node of the (i-1)-th level.

a) If $\dsum_{N\geq n} 2^{k(1,N)+\cdots+k(i,N)} q_{1N}\cdots q_{iN}
< \infty$, for {\bf some} $i, 1\leq i \leq n$, then a.s. eventually,
$\mu_N(\tri) = 0$.

b) If $\sum\limits_{N\geq n}
2^{-(k(1,N)+\cdots+k(i,N))}q_{1N}^{-1}\cdots q_{iN}^{-1} < \infty$,
for {\bf each} $i= 1, \cdots, n$, then for any $\ep>0$, a.s. eventually,
\[(1-\ep)q_{1N}\cdots q_{nN} \leq \mu_N(\tri) \leq (1+\ep)q_{1N}\cdots
q_{nN}.\]

\end{cor}

\section{Exponentially Decaying Driving Distributions}\label{s2.4}

We fix a number $\ga> 0$. In this section we consider an $n$ level
GREM where for the $N$ particle system the random variables
$\xi(\si_1\cdots\si_i)$ are i.i.d. having probability density
$$\phi_{N,\ga}(x) =
\mbox{Const.} e^{-\frac{|x|^\ga}{\ga N^{\ga-1}}} \quad
-\infty<x<\infty,$$ More precisely,
\begin{equation}\label{e2.4.21}
\phi_{N,\ga}(x)
=
\frac{1}{2\Gamma(\frac{1}{\ga})}\left(\frac{\ga}{N}\right)^{\frac{\ga-1}{\ga}}
e^{-\frac{|x|^{\ga}}{\ga N^{\ga-1}}} \quad -\infty<x<\infty.
\end{equation}

Note that $\phi_{N,1}$ is independent of $N$ and is two sided
exponential density with parameter 1. On the other hand,
$\phi_{N,2}$ is Gaussian density with mean 0 and variance $N$. Of
course, $\ga$ can be larger than 2 as well.

If we define the map $\Si_N=\prod_i 2^{k(i,N)} \rightarrow
\mathbb{R}^n$ by
\[\si \mapsto \left(\frac{\xi(\si_1,\om)}{N}, \frac{\xi(\si_1\si_2,\om)}{N}, \cdots,
\frac{\xi(\si_1\cdots\si_n,\om)}{N}\right)\] and transport the
uniform probability of $\Si$ to $\mathbb{R}^n$, we get a probability
$\mu_N(\om)$ on $\mathbb{R}^n$. In evaluating the free energy, we
will be applying Varadhan's lemma (Proposition \ref{p0.3.5}). This
explains the factor $\frac{1}{N}$ in the above map, which was not
present in the general framework of Theorem \ref{t2.3.3}.

Let $\tri = \tri_1 \times \cdots \times \tri_n$ be a non-empty open
rectangle of $\mathbb{R}^n$. For such $\tri$ and $1 \leq i \leq n$
define $m_i = \dsty\inf_{x \in \tri_i} |x|$ and $M_i = \dsty\sup_{x
\in \tri_i}|x|$. Clearly, $m_i <\infty$ for all $i$. Observe that in
case $m_i>0$ then $\tri_i \subseteq (-M_i, -m_i) \cup (m_i, M_i)$
and in case $m_i=0$ then $\tri_i \subseteq (-M_i,M_i)$. In any case
$\tri_i \subseteq (-M_i, -m_i] \cup [m_i, M_i)$ for each $i$. Let
$\tilde{m} = (m_1, \cdots ,m_n)$. Also define $q_{iN}=
P(\frac{\xi}{N} \in \tri_i)$, for $1\leq i \leq n$.

{\em First let us assume that $\ga\geq 1$.} Let $J\subset \mbb{R}$
be an interval. Denote $m = \dsty\inf_{x \in J} |x|$ and $M =
\dsty\sup_{x \in J}|x|$. Denote $q_N=P(\frac{\xi}{N}\in J)$. With
these notations, we have the following two observations:
\begin{equation}\label{e2.4.3}
\begin{split} q_{N}=P\left(\frac{\xi}{N} \in J\right) &\leq
\frac{1}{\Gamma(\frac{1}{\ga})}\left(\frac{\ga}{N}\right)^{\frac{\ga-1}{\ga}}
\dsty\int_{Nm}^{NM}
e^{-\frac{x^\ga}{\ga N^{\ga-1}}} dx\\
&<\frac{1}{\Gamma(\frac{1}{\ga})}\dsty\int_{N\frac{m^\ga}{\ga}}^{\infty}
x^{-\frac{\ga-1}{\ga}}e^{-x} dx \\
&\leq\frac{\ga^{\frac{\ga-1}{\ga}}}{\Gamma(\frac{1}{\ga})(Nm^\ga)^{\frac{\ga-1}{\ga}}}
\dsty\int_{N\frac{m^\ga}{\ga}}^{\infty}
e^{-x} dx\\
&=\frac{\ga^{\frac{\ga-1}{\ga}}}{\Gamma(\frac{1}{\ga})(Nm^\ga)^{\frac{\ga-1}{\ga}}}
e^{-N\frac{m^\ga}{\ga}},
\end{split}
\end{equation}
with the understanding that when $m=0$, the last expression is 1 and

\begin{equation}\label{e2.4.4}
\begin{split}
q_{N}=P\left(\frac{\xi}{N} \in J\right) &\geq
\frac{1}{2\Gamma(\frac{1}{\ga})}\left(\frac{\ga}{N}\right)^{\frac{\ga-1}{\ga}}
\dsty\int_{Nm}^{NM}
e^{-\frac{x^\ga}{\ga N^{\ga-1}}} dx\\
&>\frac{1}{2\Gamma(\frac{1}{\ga})}\dsty\int_{N\frac{m^\ga}{\ga}}^{N\frac{m^\ga}{\ga}+\de}
x^{-\frac{\ga-1}{\ga}}e^{-x} dx \\
&>\frac{\de}{2\Gamma(\frac{1}{\ga})(N\frac{m^\ga}{\ga}+\de)^{\frac{\ga-1}{\ga}}}
e^{-(N\frac{m^\ga}{\ga}+\de)},
\end{split}
\end{equation}
for any $0<\de < \frac{1}{\ga}(M^\ga - m^\ga)$.

{\em Now let $\ga< 1$.} $J, m, M$ as above except that $J$ is now
assumed to be bounded interval of $\mbb{R}$ so that $0\leq
m<M<\infty$. With $q_N$ as earlier, we have
\begin{equation}\label{e2.4.5}
\begin{split} q_{N}=P\left(\frac{\xi}{N} \in J\right) &\leq
\frac{1}{\Gamma(\frac{1}{\ga})}\left(\frac{\ga}{N}\right)^{\frac{\ga-1}{\ga}}
\dsty\int_{Nm}^{NM}
e^{-\frac{x^\ga}{\ga N^{\ga-1}}} dx\\
&=\frac{1}{\Gamma(\frac{1}{\ga})}\dsty\int_{N\frac{m^\ga}{\ga}}^{N\frac{M^\ga}{\ga}}
x^{-\frac{\ga-1}{\ga}}e^{-x} dx \\
&\leq\frac{\ga^{\frac{\ga-1}{\ga}}}{\Gamma(\frac{1}{\ga})(NM^\ga)^{\frac{\ga-1}{\ga}}}
\dsty\int_{N\frac{m^\ga}{\ga}}^{\infty}
e^{-x} dx\\
&=\frac{\ga^{\frac{\ga-1}{\ga}}}{\Gamma(\frac{1}{\ga})(NM^\ga)^{\frac{\ga-1}{\ga}}}
e^{-N\frac{m^\ga}{\ga}},
\end{split}
\end{equation}
with the understanding that when $m=0$, the last expression is 1.
The difference between (\ref{e2.4.5}) and (\ref{e2.4.3}) is just
that in the penultimate inequality the lower bound of the  integral
appeared in (\ref{e2.4.3}) where as in (\ref{e2.4.5}), the upper
bound of the integral appeared.

\begin{equation}\label{e2.4.6}
\begin{split}
q_{N}=P\left(\frac{\xi}{N} \in J\right) &\geq
\frac{1}{2\Gamma(\frac{1}{\ga})}\left(\frac{\ga}{N}\right)^{\frac{\ga-1}{\ga}}
\dsty\int_{Nm}^{NM}
e^{-\frac{x^\ga}{\ga N^{\ga-1}}} dx\\
&>\frac{1}{2\Gamma(\frac{1}{\ga})}\dsty\int_{N\frac{m^\ga}{\ga}}^{N\frac{m^\ga}{\ga}+\de}
x^{-\frac{\ga-1}{\ga}}e^{-x} dx \\
&>\frac{\de\ga^{\frac{\ga-1}{\ga}}}{2\Gamma(\frac{1}{\ga})(Nm^\ga)^{\frac{\ga-1}{\ga}}}
e^{-(N\frac{m^\ga}{\ga}+\de)},
\end{split}
\end{equation}
for any $0<\de < \frac{1}{\ga}(M^\ga - m^\ga)$. The difference
between (\ref{e2.4.6}) and (\ref{e2.4.4}) is just that in the
penultimate inequality the upper limit of the  integral appeared in
(\ref{e2.4.4}) where as in (\ref{e2.4.6}), the upper limit of the
integral to bound $e^{-x}$ and lower limit to bound
$x^{-\frac{\ga-1}{\ga}}$ is used. Moreover, when $m=0$ the lower
bound for $q_N$ can be given by
$\frac{1}{2\Gamma(\frac{1}{\ga})}\dsty\int_{\de_1}^{\de_2}x^{-\frac{\ga-1}{\ga}}e^{-x}
dx$ for $0<\de_1<\de_2<\frac{M^\ga}{\ga}$. As earlier, this bound
does not depend on $N$.

From now on we assume that $\frac{k(i,N)}{N} \ra p_i$ for $1\leq i
\leq n$ with $p_1>0$. Clearly, $\sum p_i = 1$. Let
\begin{equation}\Psi = \{\widetilde{x}\in \mathbb{R}^n :
\dsum_{i=1}^k \frac{|x_i|^{\ga}}{\ga} \leq \dsum_{i=1}^k p_i \log2,
\hspace{1ex} 1\leq k\leq n\}.\end{equation}

\begin{prop}\label{p2.4.1}
$\mu_N \Ra \de_0$ a.s. as $N\ra \infty$.
\end{prop}
\begin{proof}
For any $\ep>0$, define $\tri(\ep) = [-\ep,\ep]\times \cdots \times
[-\ep,\ep] \subseteq \mathbb{R}^n.$ By Markov inequality,
$$\mbf{P}(\mu_N(\tri^c(\ep))>\ep) < \frac{1}{\ep}\mbf{E}\mu_N(\tri^c(\ep) )<
\frac{n}{\ep}\mbf{P}(|\xi|>\ep N)<\frac{2n}{\ep}\mbf{P}(\xi>\ep N) <
\frac{2n}{\ep}C_N e^{-N\frac{\ep^\ga}{\ga}},$$ where $C_N$ can be
obtained from (\ref{e2.4.3}) for $\ga\geq 1$ and from (\ref{e2.4.5})
for $0<\ga<1$. Since $\frac{1}{N}\log C_N \ra0$ as $N\ra \infty$,
the proposition follows from the Borel-Cantelli lemma.
\end{proof}

\begin{prop}\label{p2.4.2}
If $\bar{\tri} \cap \Psi = \phi$, then a.s. eventually $\mu_N (\tri)
= 0$. Moreover, the sequence $\{\mu_N\}$ is supported on a compact
set.
\end{prop}

\begin{proof}
$\bar{\tri} \cap \Psi = \phi$ implies $\tilde{m}\notin \Psi$. This
is seen as follows. By definition of $m_i$, either $m_i$ or $-m_i$
is in $\bar\tri_i$. Thus for each $i$, there is an $\ep_i=\pm 1$
such that $\ep_i m_i\in \bar\tri_i$. Thus the vector $(\ep_1
m_1,\cdots,\ep_n m_n)\in \bar\tri$ and hence $\notin \Psi$. By the
symmetry of $\Psi$, $\tilde{m}\notin \Psi$ as well. As a
consequence, for some $j$, $1\leq j \leq n$,
\begin{equation}\label{e2.4.7}
\sum\limits_{i=1}^j \frac{m_i^\ga}{\ga} > \sum\limits_{i=1}^j p_i
\log 2.
\end{equation}
For $\ga \geq 1$ using (\ref{e2.4.3}) and for $0<\ga<1$ using
(\ref{e2.4.5}) we can say that $q_{iN}<C_{iN}
e^{-N\frac{m_i^\ga}{\ga}}$ where $\frac{1}{N}\log C_{iN}\ra 0$ as
$N\ra \infty$ for $1\leq i\leq j$. Hence as a consequence of
(\ref{e2.4.7}) and the fact $\frac{k(i,N)}{N}\ra p_i$, we have
$$\sum\limits_{N\geq 1}2^{k(1,N)+\cdots+k(j,N)}q_{1N}\cdots q_{jN}
< \sum\limits_{N\geq 1}e^{-N\sum\limits_{i=1}^j
\left(\frac{m_i^\ga}{\ga}-\frac{k(i,N)}{N}\log 2-\frac{1}{N}\log C_{iN}\right)} < \infty.$$

Thus by Corollary \ref{c2.3.4}, a.s. eventually $\mu_N (\tri) = 0$.

To see the last statement of the Proposition, fix any $\de>0$. Let
$J$ be the compact set $[-\log2-\de, +  \log2+\de]^n$. Since the
complement of this set is union of $2^n$ open rectangles of
$\mathbb{R}^n$, each of whose closures are disjoint with $\Psi$, the
earlier part implies that eventually $\mu_N(J)=1$.
\end{proof}
\begin{prop}\label{p2.4.3}
If $(\bar{\tri} \cap \Psi)^0 \neq \phi$, then for any $\ep
> 0$ a.s. eventually
\[(1-\ep)q_{1N} \cdots q_{nN}\leq \mu_N(\tri) \leq (1+\ep)q_{1N} \cdots q_{nN}.\]
\end{prop}

\begin{proof}

The assumption $(\bar{\tri} \cap \Psi)^0 \neq \phi$ implies
$\tilde{m} \in \Psi^0$. Indeed, since $(\bar{\tri} \cap \Psi)^0 \neq
\phi$, pick $(x_1,\cdots,x_n)\in(\bar{\tri} \cap \Psi)^0$. By
symmetry of $\Psi$, $(|x_1|, \cdots, |x_n|)\in \Psi^0$ as well, and
now $0\leq m_i \leq |x_i|$ for all $i$ yields $(m_1,\cdots,m_n) \in
\Psi^0$.

We are going to show that the hypothesis of Corollary
\ref{c2.3.4}(b) holds. Fix $i$, $1\leq i\leq n$. Using
(\ref{e2.4.4}) for $\ga\geq 1$ and using (\ref{e2.4.6}) for $\ga<1$,
we can say that $q_{jN}>C_{jN} e^{-N\frac{m_j^\ga}{\ga}+\de}$ for
$1\leq j\leq n$ with sufficiently small $\de>0$. Thus
$$2^{-(k(1,N)+\cdots+k(i,N))}q_{1N}^{-1}\cdots q_{iN}^{-1} <
e^{-N\left[ \sum\limits_{j=1}^i\left(\frac{k(j,N)}{N}\log2 -
\frac{m_j^\ga}{\ga}+\frac{1}{N}\log C_{iN}\right) - i\de\right]}.$$
Since $\tilde{m}$ is an interior point of $\Psi$, there is an
$\al>0$ such that $\sum\limits_{j=1}^i p_j\log2 -
\sum\limits_{j=1}^i \frac{m_j^\ga}{\ga}>\al$. Now use the fact that
$\frac{k(j,N)}{N} \ra p_j$ and $\frac{1}{N}\log C_{jN}\ra 0$ as
$N\ra \infty$ to deduce that eventually
$\sum\limits_{j=1}^i\left(\frac{k(j,N)}{N}\log2-\frac{m_j^\ga}{\ga}+
\frac{1}{N}\log C_{iN}\right)>\al$. Making $\de>0$ smaller, if
necessary, assume that eventually
$\sum\limits_{j=1}^i\left(\frac{k(j,N)}{N}\log2-\frac{m_j^\ga}{\ga}+
\frac{1}{N}\log C_{iN}\right)-i\de>\al$. Hence, eventually
$e^{-N\left[ \sum\limits_{j=1}^i\left(\frac{k(j,N)}{N}\log2 -
\frac{m_j^\ga}{\ga}+\frac{1}{N}\log C_{iN}\right) -
i\de\right]}<e^{-N\al}$. As a consequence,
$$\sum\limits_{N\geq 1}2^{-(k(1,N)+\cdots+k(i,N))}q_{1N}^{-1}\cdots q_{iN}^{-1} < \infty.$$
Hence by Corollary \ref{c2.3.4}, the proposition follows.
\end{proof}

\smallskip
\begin{rem}\label{r2.4.1}
{\em $(\bar{\tri} \cap \Psi)^0 \neq \phi$ implies in particular,
that $p_1>0$. In fact, $\Psi^0 \neq \phi$ iff $p_1>0$. }
\end{rem}
Now, we have the following,

\begin{prop}\label{p2.4.4}
For a.e. sample point $\om$,
\[\begin{array}{llll}
\dsty\lim_{N\ra \infty} \frac{1}{N} \log\mu_N(\tri) &=&
-\dsum_{i=1}^n \frac{m_i^\ga}{\ga} & \mbox{ if } (\bar{\tri} \cap
\Psi)^0 \neq
\phi\\
&=& -\infty & \mbox{ if } \bar{\tri}\cap \Psi = \phi.
\end{array}\]
\end{prop}

\begin{proof}
When $\bar{\tri}\cap \Psi = \phi$, the result is immediate from
Proposition \ref{p2.4.2}.

Assume that $(\bar{\tri} \cap \Psi)^0 \neq \phi$. Fix any $\ep$,
$0<\ep<1$. Let $\ga\geq 1$. By (\ref{e2.4.3}), $\frac{1}{N}\log
q_{iN}< \frac{1}{N} \log C_{iN}-\frac{m_i^\ga}{\ga}$ where
$\frac{1}{N} \log C_{iN}\ra 0$ as $N\ra \infty$. Hence $\limsup
\frac{1}{N}\log q_{iN}\leq \frac{m_i^\ga}{\ga}$. Similarly, by using
(\ref{e2.4.4}) we get $\liminf \frac{1}{N}\log q_{iN}\geq
\frac{m_i^\ga}{\ga}$. Thus $\lim \frac{1}{N}\log q_{iN}$ exists and
equals to $\frac{m_i^\ga}{\ga}$ for each $i$. The same holds even if
$0<\ga<1$, where we need to use (\ref{e2.4.5}) and (\ref{e2.4.6}).
Then by proposition \ref{p2.4.3} we have a.s eventually, a.s.
eventually
\[(1-\ep)q_{1N} \cdots q_{nN}\leq \mu_N(\tri) \leq (1+\ep)q_{1N} \cdots q_{nN}.\]
So by taking logarithms and using $\dsty\lim_{N\ra
\infty}\frac{1}{N} \log q_{iN} = -\frac{m_i^\ga}{\ga},$ for each $i$
we get the proposition.
\end{proof}

Let us consider the map $I: \mathbb{R}^n \ra \mathbb{R}$, defined as
follows,
\begin{equation}\label{e2.4.8}
\begin{array}{llll}
I(\widetilde{x}) &=& \frac{1}{\ga}\dsum_{i=1}^n |x_i|^\ga &\mbox{ if
} \widetilde{x} \in
\Psi\\
&=& \infty &\mbox{ otherwise}.
\end{array}
\end{equation}
\begin{thm}\label{t2.4.5}
Almost surely, the sequence $\{\mu_N\}$ satisfies LDP with the rate
function $I$.
\end{thm}

\begin{proof}
Let $\mathcal{A}$ be the collection of all rectangles $\tri =
\tri_1\times\cdots\times\tri_n \subseteq \mathbb{R}^n$ such that
each $\tri_i$ is a bounded interval with rational endpoints and
either $\bar{\tri}\cap \Psi = \phi$ or $(\bar{\tri} \cap \Psi)^0
\neq \phi.$

It is easy to check that $\mathcal{A}$ forms a base for the usual
topology of $\mathbb{R}^n$. For $\tri \in\mathcal{A}$, by
Proposition \ref{p2.4.4}, the limit, $ - \dsty\lim_{N\ra \infty}
\frac{1}{N} \log \mu_N(\tri)$ exists almost surely. Denote this
limit by $L_{\tri}$. Since $\mathcal{A}$ is a countable family, out
side a null set, these limits are well defined for all $\tri
\in\mathcal{A}$.

In view of Proposition \ref{p0.3.3}, to complete the proof, we show
that for $\tilde{x}\in \mathbb{R}^n$,
\begin{equation}\label{e2.4.9}
I(\tilde{x}) = \dsup_{\tilde{x} \in \tri \in \mathcal{A}}L_{\tri}.
\end{equation}

If $\tilde{x} \notin \Psi$, clearly $\sup\limits_{\tilde{x} \in \tri
\in \mathcal{A}} L_\tri = \infty = I(\tilde{x}).$

Now consider, $\tilde{x}=(x_1,\cdots,x_n) \in \Psi$. Suppose
$\tilde{x}\in\tri\in\mathcal{A}$. If
$\tri=\tri_1\times\cdots\times\tri_n$ with
$m_i=\inf\limits_{y\in\tri_i}|y|$, then $m_i\leq |x_i|$ and hence
$\frac{m_i^\ga}{\ga}\leq \frac{|x_i|^\ga}{\ga}$. Therefore, by
Proposition \ref{p2.4.4}, $L_{\tri}=\sum\limits_{i=1}^n
\frac{m_i^\ga}{\ga}\leq\sum\limits_{i=1}^n \frac{|x_i|^\ga}{\ga}$.
Thus
\begin{equation}\label{e2.4.10}
\dsup_{\tilde{x} \in \tri \in \mathcal{A}}L_{\tri}\leq I(\tilde{x}).
\end{equation}
On the other hand, consider $\ep>0$ so that $\ep<|x_i|$ for any $i$
with $x_i\neq 0$. Let $\tri$ be the box with sides
$\tri_i=(x_i-\ep,x_i+\ep)$. By choice of $\ep$,
$m_i=\dinf_{y\in\tri_i} |y|$ equals $|x_i\pm \ep|$ depending on the
sign of $x_i$. Of course, if $x_i=0$ then $m_i=0$. Thus for the
$\tri$ so constructed, we have, $L_{\tri}=\dsum_{\{i:x_i\neq0\}}
\frac{|x_i\pm \ep|^\ga}{\ga}$. This being true for all sufficiently
small $\ep$, we conclude that
\begin{equation}\label{e2.4.11}
\dsup_{\tilde{x} \in \tri \in \mathcal{A}}L_{\tri}\geq \dsum_{i=1}^n
\frac{|x_i|^\ga}{\ga}=I(\tilde{x})
\end{equation}
(\ref{e2.4.10}) and (\ref{e2.4.11}) complete the proof of
(\ref{e2.4.9}) thus completing the proof of the theorem.
\end{proof}
\smallskip

We shall now proceed towards an expression for the free energy.
Denoting $f(\tilde{x})= \dsum_{i=1}^n\be a_i x_i$,
$$\begin{array}{lll}
\dlim_N \frac{1}{N} \log Z_N(\be) &=&\log2+\dlim_N \frac{1}{N} \log
\mbf{E}_N e^{-Nf}\\
&=&\log2-\dinf_{\tilde{x}\in\Psi}\{\dsum_{i=1}^n\be a_i
x_i+\dsum_{i=1}^n \frac{|x_i|^\ga}{\ga}\}.
\end{array}$$
by Proposition \ref{p0.3.5}. This last infimum equals
$\dinf_{\tilde{x}\in\Psi}\dsum_{i=1}^n\left(
\frac{|x_i|^\ga}{\ga}-\be a_i x_i\right)$. Since $\be>0$, $a_i>0$ it
is easy to see that the above infimum is attained when all the $x_i$
are negative. In other words, by symmetry of $\Psi$, the infimum is
attained at a point $-\tilde{x}$ for some $\tilde{x}\in \Psi^+ =
\Psi \cap \{\tilde{x} : x_i \geq 0 \mbox{ for }1\leq i\leq n\}.$
Thus
$$\dlim_N \frac{1}{N} \log
Z_N(\be)=\log2-\dinf_{\tilde{x}\in\Psi^+}\dsum_{i=1}^n \left(
\frac{x_i^\ga}{\ga}-\be a_i x_i\right).$$

In this way, for the above mentioned class of driving distributions,
the free energy exists almost surely and is a constant. Not only
that, finding an explicit formula for the free energy reduces to
calculating the above infimum.

\begin{rem}
It is also worth noting that the LDP holds good even when the
driving distributions at various levels are different. To be more
specific, let us fix $n$ numbers $\ga_1,\cdots,\ga_n$; each greater
than zero and consider an $n$ level GREM where the driving
distribution at the $i$-th level is $\phi_{N,\ga_i}$. More
precisely, for any node $\si_1\cdots\si_i$ at the $i$-th level
$\xi(\si_1\cdots\si_i)$ has density $\phi_{N,\ga_i}$. Of course, all
the random variables are independent. Define as earlier, the map
$\Si_N \ra \mathbb{R}^n$ by
\[\si \mapsto \left(\frac{\xi(\si_1)}{N}, \frac{\xi(\si_1\si_2)}{N}, \cdots,
\frac{\xi(\si_1\cdots\si_n)}{N}\right).\] Let $\mu_N$ be the induced
probability on $\mathbb{R}^n$ when $\Si_N$ is equipped with uniform
probability. The same arguments as above, with
$q_{iN}=P\left(\frac{\xi(\si_1\cdots\si_i)}{N}\in\tri_i\right)$,
will show that almost surely, the sequence of probabilities
$\{\mu_N, N\geq n\}$ on $\mathbb{R}^n$ satisfies LDP. In this case,
with rate function $I$ will be given by
\end{rem}

\begin{equation}\label{e2.4.13}
\begin{array}{llll} I(\widetilde{x}) & = & \dsum_{i=1}^n
\frac{|x_i|^{\ga_i}}{\ga_i} & \mbox{ if } \widetilde{x} \in\Psi \\
                      & =& \infty& \mbox{ otherwise},
\end{array}
\end{equation}
where
\begin{equation}\label{e2.4.14}
\Psi = \{\widetilde{x}\in \mathbb{R}^n :
\dsum_{i=1}^k \frac{|x_i|^{\ga_i}}{\ga_i} \leq \dsum_{i=1}^k p_i
\log2, \hspace{1ex} 1\leq k\leq n\},
\end{equation}
with $p_i=\dlim_{N\ra\infty}
\frac{k(i,N)}{N}$. Let, as earlier, $\Psi^+$ be the part of $\Psi$
in the positive orthant of $\mathbb{R}^n$. As a consequence of all
this, we have the following:
\begin{thm}\label{t2.4.6}
If the driving distribution has density $\phi_{N,\ga_i}$ at the
$i$-th level, we have almost surely,
\[\dsty\lim_{N} \frac{1}{N} \log
Z_N(\be)=\log2-\dinf_{\widetilde{x}\in\Psi^+}\left\{\dsum_{i=1}^n\left(\frac{x_i^{\ga_i}}{\ga_i}-\be
a_i x_i\right)\right\}.\]
\end{thm}

\section{Inside Out}\label{s2.5}
A close observation of the above discussion reveals the following
cute idea. Though the identification, at first glance, will look
like very simple, its implication in GREM will be understood through
the rest of this chapter.

Let for each $j$, $1\leq j\leq n$, we have a sequence of
probabilities $\{\la_N^j, N\geq n\}$ on $\mbb{R}$ which obey LDP
with a strictly quasi-convex continuous good rate function
$\mathcal{I}_j$. That is, $\mathcal{I}_j$ has compact level sets and
for any two distinct points $x$ and $y$ in
$\{0<\mathcal{I}_j<\infty\}$ we have $\mathcal{I}_j(\theta
x+(1-\theta)y)<\max\{\mathcal{I}_j(x),\mathcal{I}_j(y)\}$ for any
$\theta$ with $0<\theta<1$. For the sake of simplicity, we will also
assume that $\mathcal{I}_j(0)=0$. The assumption of strict
quasi-convexity is purely a technical assumption and this can be
replaced by similar other conditions also. For example, one can
replace this by requiring that the set $\{x:\,
\mathcal{I}_j(x)=\al\}$ is a nowhere dense set for every $\al>0$. We
mentioned this condition in Remark \ref{r1.2.2}, but there we
demanded this only for $\al=\log 2$. Now, let us denote
$\{+1,-1\}^N$ by $\Si_N$. For each $N$, let $k(1,N),\ldots, k(n,N)$
be non-negative integers adding to $N$ and put
$\Si_{jN}=\{+1,-1\}^{k(j,N)}$. Clearly,
$\Si_N=\Si_{1N}\times\Si_{2N}\times\cdots\times\Si_{nN}$ and we
express $\sigma\in \Si_N$ as $\si_1\si_2\cdots\si_n$ with
$\si_i\in\Si_{iN}$, in an obvious way. Suppose for fixed $N$, we
have a bunch of independent random variables as follows:
$\{\xi(\si_1):\, \si_1\in\Si_{1N}\}$ having distributions $\la_N^1$,
$\{\xi(\si_1\si_2):\,\si_2\in\Si_{2N}, \si_1\in\Si_{1N}\}$ having
distributions $\la_N^2$ and in general
$\{\xi(\si_1\si_2\cdots\si_{j-1}\si_j):\, \si_j\in\Si_{jN},
\cdots,\si_1\in\Si_{1N}\}$ having distribution $\la_N^j$.

Define for each $\omega$, $\mu_N(\omega)$ to be the empirical
measure on $\mbb{R}^n$, namely,
$$\mu_N(\omega)=\frac{1}{2^N}\sum_\si
\delta\left<\xi(\si_1,\omega),\xi(\si_1\si_2,\omega),\cdots,\xi(\si_1\cdots\si_n,\omega)\right>$$
where $\delta\left<x\right>$ denotes the point mass at
$x\in\mbb{R}^n$.
\begin{thm}\label{t2.5.1}
Suppose $\frac{k(j,N)}{N}\ra p_j>0$ for $1\leq j\leq n$. Then for
a.e. $\omega$, the sequence $\{\mu_N(\omega), N\geq n\}$ satisfies
LDP with rate function $\mathcal{J}$ given as follows:

Supp$(\mathcal{J})=\{(x_1,\cdots,x_n): \dsum_{k=1}^j
\mathcal{I}_k(x_k) \leq \dsum_{k=1}^j p_k\log2\; \mbox{ for }\,
1\leq j\leq n \}$

and $$\begin{array}{rll} \mathcal{J}(x)& = \dsum_{k=1}^n
\mathcal{I}_k(x_k)& \mbox{if }
x\in \mbox{Supp}(\mathcal{J})\\
&=\infty&\mbox{otherwise}.
\end{array}$$
\end{thm}
\begin{proof} In what follows $\tri$ denotes a box in $\mbb{R}^n$ with
sides $\tri_j;\; 1\leq j\leq n$ where each $\tri_j$ is an interval.
The proof consists of the following steps. The steps are executed
one by one as in Propositions  \ref{p2.4.2} to \ref{p2.4.4}, so will
not be repeated here.

{\bf Step 1:} {\em If $\bar{\tri} \cap \mbox{Supp}(\mathcal{J}) =
\phi$, then a.s. eventually $\mu_N (\tri)
= 0$.}\\

{\bf Step 2:} {\em If $(\bar{\tri} \cap \mbox{Supp}(\mathcal{J}))^0
\neq \phi$, then for any $\ep>0$ a.s. eventually
$$(1-\ep)\prod_{i=1}^n\la_N^i(\tri_i) \leq \mu_N(\tri)\leq
(1+\ep)\prod_{i=1}^n\la_N^i(\tri_i).$$}

{\bf Step 3:} {\em For a.e. sample point $\om$,
\[\begin{array}{llll}
\dsty\lim_{N\ra \infty} \frac{1}{N} \log\mu_N(\tri) &=&
-\dsum_{i=1}^n \mathcal{I}_i(\tri_i) & \mbox{ if } (\bar{\tri} \cap
\mbox{Supp}(\mathcal{J}))^0 \neq
\phi\\
&=& -\infty & \mbox{ if } \bar{\tri}\cap \mbox{Supp}(\mathcal{J}) =
\phi,
\end{array}\]
where $\mathcal{I}_i(\tri_i)=\inf\{\mathcal{I}(x):\, x\in \tri_i\}$.}

To conclude the proof we use the idea of Theorem \ref{t2.4.5}.
\end{proof}

We note that continuity of the rate functions $\mathcal{I}_j$ is not
necessary, but then one needs to go through $\limsup$ and $\liminf$
of $\frac{1}{N}\log\mu_N(\tri)$ an in Theorem \ref{t1.2.1}, instead
of limits which we used above.

The implications of the above theorem for GREM~\cite{D2} are clear.
For fixed $N$, and $\si\in\Si_N$ one defines the Hamiltonian
$$H_N(\si)=N\sum_{i=1}^n a_i \xi(\si_1\cdots\si_i).$$ Here $a_i,\;
1\leq i\leq n$ are positive numbers called weights. In the Gaussian
case, it is customary to take $\sum a_i^2=1$, though it is not a
mathematical necessity. As earlier, $Z_N(\be)=\sum_\si e^{-\be
H_N(\si)}$. Special choices of $\la_N^i$ lead to all the known
models considered. Centered Gaussian were consider in
\cite{D2,CCP,DD,JR1}. More general distributions as well as the
cases when some $p_j$ are zero were considered in \cite{JR1}.
Moreover one could take different distributions for different values
of $j$, see \S \ref{s2.7} for some interesting consequences. Thus
the main problem of GREM is reduced to a variational problem. Note
that, if $n=1$, GREM reduces to REM.

\section{The Variational Problem}

In this section, we derive explicit formulae for the free energy. We
return back to the driving distribution given by (\ref{e2.4.21}),
namely, having density
\begin{equation}\label{e2.6.a12}
\phi_{N,\ga}(x) =
\frac{1}{2\Gamma(\frac{1}{\ga})}\left(\frac{\ga}{N}\right)^{\frac{\ga-1}{\ga}}
e^{-\frac{|x|^{\ga}}{\ga N^{\ga-1}}} \quad -\infty<x<\infty.
\end{equation}
 We now consider the model with same driving distributions at different
levels. In this setup, we need to calculate the infimum
\begin{equation}\label{e2.6.a12a}
 \dinf_{\tilde{x}\in\Psi^+}\dsum_{i=1}^n \left(
\frac{x_i^\ga}{\ga}-\be a_i x_i\right)
\end{equation}
in order to get an explicit
formula for the limiting free energy. Note that, putting $\ga_i=\ga$
got all $i$ in (\ref{e2.4.14}) we get
\begin{equation}\label{e2.6.a12b}
 \Psi=\{\widetilde{x}\in\mbb{R}^n:\, \sum_{i=1}^k|x_i|^\ga \leq
\sum_{i=1}^k\ga p_i\log2,\quad 1\leq k\leq n\}
\end{equation}
and as usual
$\Psi^+$ is the part of $\Psi$ in the positive orthant of
$\mbb{R}^n$.
\smallskip
\subsection{\texorpdfstring{$\ga>1$}{gamma>1}}\label{s2.6.1}
First let us assume $\ga>1$. To evaluate the infimum let us put, for
$1 \leq j \leq k \leq n$,
\begin{equation}\label{e2.5.a12}
B(j,k) =
\left(\frac{(p_j+\cdots+p_k)\ga\log2}{a_j^\frac{\ga}{\ga-1}+
\cdots+a_k^\frac{\ga}{\ga-1}}\right)^\frac{\ga-1}{\ga}.
\end{equation}
Set $r_0=0$ and for $l\geq0$ (integer),
\begin{equation}\label{e2.5.b12}
\be_{l+1} = \dmin_{k>r_l}B(r_l+1,k)\quad r_{l+1} = \max\{i>r_l:
B(r_l+1,i)=\be_{l+1}\}.
\end{equation}

Clearly, for some $K$ with $1\leq K\leq n$, we have $r_K=n$. Put
$\be_0=0$ and $\be_{K+1}=\infty$, so  that $0=\be_0 <\be_1 <
\be_2\cdots< \be_K <\be_{K+1}=\infty$.

Fix $j\leq K$ and let  $\be\in( \be_j, \be_{j+1}].$  Define
$\widetilde{\ol{x}}\in\Psi^+$ as follows:
\begin{equation}\label{e2.5.12a}
\begin{array}{llll}
\ol{x}_i &=& (\be_l a_i)^\frac{1}{\ga-1} &\mbox{ if } i\in
\{r_{l-1}+1,\cdots,r_l\}
\mbox{ for some }l, 1\leq l \leq j\\
&=& (\be a_i)^\frac{1}{\ga-1} &\mbox{ if } i\geq r_j+1.
\end{array}\end{equation}
\noindent{\bf Claim:} $\dsty\inf_{\widetilde{x}\in \Psi^+}
\dsum_{i=1}^n\left(\frac{x_i^\ga}{\ga}-\be a_ix_i\right)$ occurs at
$\widetilde{\ol{x}}.$

In order to prove the claim, fix any  $\widetilde{x}\in \Psi^+$. For
$k\leq j$ (recall that $j\leq K$ was fixed above), first note that,
by Holder's inequality,
$$\dsum_{i=1}^{r_k} x_i \ol{x}_i^{\ga-1}
\leq \left(\dsum_{i=1}^{r_k} x_i^\ga\right)^{\frac{1}{\ga}}
\left(\dsum_{i=1}^{r_k}\ol{x}_i^\ga\right)^\frac{\ga-1}{\ga} \leq
\dsum_{i=1}^{r_k}\ol{x}_i^\ga.$$ where the last inequality follows
from the facts $\widetilde{x}\in\Psi^+$ and
$\dsum_{i=1}^{r_k}\ol{x}_i^\ga=\dsum_{i=1}^{r_k}\ga p_i\log2$ so
that $\dsum_{i=1}^{r_k}x_i^\ga\leq \dsum_{i=1}^{r_k}\ga
p_i\log2=\dsum_{i=1}^{r_k}\ol{x}_i^\ga$.

Hence, $\dsum_{i=1}^{r_k}
\ol{x}_i^{\ga-1}(\ol{x}_i-x_i)\geq 0$.\\

Since $\be>\be_j$, we have $(\frac{\be}{\be_l}-1)>0$ for $1\leq
l\leq j$. Moreover since $\be_l$ are increasing with $l$, these
number $(\frac{\be}{\be_l}-1)$ are decreasing. It follows that,
$$\dsum_{l=1}^{j}\left(\frac{\be}{\be_l}-1\right)\dsum_{i=r_{l-1}+1}^{r_l}
\ol{x}_i^{\ga-1}(\ol{x}_i-x_i)\geq 0$$ In other words, using the
definition of $\ol{x}_i$,
\begin{equation}\label{e2.5.12}
\dsum_{i=1}^{r_j}\be a_i(\ol{x}_i-x_i)\geq
\dsum_{i=1}^{r_j}\ol{x}_i^{\ga-1}(\ol{x}_i-x_i).
\end{equation}
Now,
\begin{equation}\label{e2.5.13}
\begin{split}
&\dsum_{i=1}^{r_j}\left(\frac{x_i^\ga}{\ga}-\be a_i x_i\right)-
\dsum_{i=1}^{r_j}\left(\frac{\ol{x}_i^\ga}{\ga}-\be a_i \ol{x}_i\right)\\
= &\dsum_{i=1}^{r_j}\left(\frac{x_i^\ga}{\ga}+\be
a_i(\ol{x}_i-x_i)-\frac{\ol{x}_i^\ga}{\ga} \right)\\
\geq
&\dsum_{i=1}^{r_j}\left(\frac{x_i^\ga}{\ga}+\ol{x}_i^{\ga-1}(\ol{x}_i-x_i)-
\frac{\ol{x}_i^\ga}{\ga}\right)\hspace{10ex}\mbox{by (\ref{e2.5.12})}\\
=&\dsum_{i=1}^{r_j}\left(\frac{x_i^\ga}{\ga}+\frac{\ga-1}{\ga}\ol{x}_i^\ga-x_i
\ol{x}_i^{\ga-1}\right)\\
\geq& 0.
\end{split}
\end{equation}
where in the last inequality we used $x_i\ol{x}_i^{\ga-1}\leq
\frac{1}{\ga}x_i^\ga +\frac{\ga-1}{\ga}\ol{x}_i^\ga$.

On the other hand, utilizing the definition of $\ol{x}_i$ and  the
inequality $\be a_ix_i \leq
\frac{x_i^\ga}{\ga}+\frac{\ga-1}{\ga}(\be a_i)^{\frac{\ga}{\ga-1}}$
we have,
\begin{equation}\label{e2.5.14}
\begin{split}
&\dsum_{i=r_j+1}^{n}\left(\frac{x_i^\ga}{\ga}-\be a_i x_i\right)-
\dsum_{i=r_j+1}^{n}\left(\frac{\ol{x}_i^\ga}{\ga}-\be a_i \ol{x}_i\right)\\
= &\dsum_{i=r_j+1}^{n}\left(\frac{x_i^\ga}{\ga}+\frac{\ga-1}{\ga}(\be a_i)^{\frac{\ga}{\ga-1}}-\be a_i x_i\right)\\
\geq& 0.
\end{split}
\end{equation}
Clearly, (\ref{e2.5.13}) and (\ref{e2.5.14}) complete proof of the
claim. This argument is in fact a generalization of Dorlas \&
Dukes\cite{DD}, Capocaccia et. al.\cite{CCP}.

All this leads to the following explicit formula for the free
energy.
\begin{thm}\label{t2.6.1}
For GREM with driving distribution having density $\phi_{N,\ga}$ as
defined in (\ref{e2.6.a12}), almost surely,
\[\begin{array}{ll}
\dsty\lim_{N} \frac{1}{N} \log
Z_N(\be)\hspace{-1ex}&=\dsum_{i=r_j+1}^n p_i\log2+
\frac{\ga-1}{\ga}\dsum_{i=r_j+1}^n (\be a_i)^{\frac{\ga}{\ga-1}}
+\be \dsum_{l=1}^j\be_l^{\frac{1}{\ga-1}}\dsum_{i=r_{l-1}+1}^{r_l}
a_i^{\frac{\ga}{\ga-1}}\vspace{1ex}\\
&\hspace{30ex}\mbox{if }\quad \be_j<\be \leq\be_{j+1},\;0\leq j\leq K-1\vspace{1ex}\\
&=\be\dsum_{l=1}^K \be_l^{\frac{1}{\ga-1}}\dsum_{i=r_{l-1}+1}^{r_l}
a_i^{\frac{\ga}{\ga-1}}\hspace{4ex}\mbox{if }\quad\be>\be_K
\end{array}\]
\end{thm}

Observe that for $\ga=2$, that is when the driving distribution is
Normal, with proper identification of parameters this is essentially
the same formula as in \cite{CCP,DD}. In defining the $\be_i$,
Capocaccia {\em et. al.} use a variant in~\cite{CCP} ($\S$3.2). In
defining $r_i$, Dorlas and Dukes \cite{DD} consider the least index,
where as recall that, we define $r_{l+1}$ as $\max\{i>r_l:
B(r_l+1,i)=\be_{l+1}\}$. This makes no difference because {\em
`nothing happens'} in between these two indices. This follows from
the fact that if $a_i>0$, $b_i>0$ for $i=1,2,3$ and
$\frac{b_1}{a_1}=\frac{b_1+b_2+b_3}{a_1+a_2+a_3}<\frac{b_1+b_2}{a_1+a_2}$,
then this will imply that
$\frac{b_1}{a_1}=\frac{b_2+b_3}{a_2+a_3}<\frac{b_2}{a_2}$. So
defining $r_{l+1}$ as $\min\{i>r_l: B(r_l+1,i)=\be_{l+1}\}$, when
$\{i>r_l: B(r_l+1,i)=\be_{l+1}\}$ is not a singleton set,
$\be_{l+2}$ will be same as $\be_{l+1}$. And this will continue
until the maximum index of the set $\{i>r_l: B(r_l+1,i)=\be_{l+1}\}$
is attained.

Moreover, the weights $a_i$ in Dorlas and Dukes \cite{DD} are
incorporated in the density, there was no need to assume $\sum
a_i=1$, their parameter $J$ can be incorporated in the weights. In
fact, there is one benefit of putting the weights in the density.
The large deviation technique will easily allow us to consider
variable weights $a_{iN}$ depending on $N$ at the $i$-th level
instead of a constant weights $a_i$. For instance, let $a_{iN}>0$
for all $1\leq i\leq n$ and $N\geq 1$ be the weights of the $i$-th
level for the $N$ particle system with $a_{iN} \ra a_i$ as
$N\ra\infty$. When the weights $a_i$ did not depend on $N$, they
were not brought in the large deviation argument. The free energy
was
\begin{equation}\label{e2.5.14b}
\log2 - \dinf_{\tilde{x}\in\Psi^+}\dsum_{i=1}^n \left(
\frac{x_i^2}{2}-\be a_i x_i\right),
\end{equation}
where $$\Psi^+=\{\widetilde{x}\in\mathbb{R}^n:\,
\dsum_{i=1}^k\frac{1}{2}x_i^2\leq \dsum_{i=1}^k p_i\log
2\}\cap\{\widetilde{x}\in\mathbb{R}^n:\, x_i\geq 0, 1\leq i\leq
n\}.$$ If we consider variable weights $a_{iN}$ as above then they
must enter in the large deviation arguments. If $\xi(i,N)\sim
\mathcal{N}(0,\frac{1}{N})$, then it is not hard to show that the
distribution of $a_{iN}\xi(i,N)$ satisfies LDP with rate function
$\frac{x^2}{2a_i^2}$. Accordingly, we get the limiting free energy
as
\begin{equation}\label{e2.5.14c}
\log2 - \dinf_{\tilde{x}\in\Upsilon^+}\dsum_{i=1}^n \left(
\frac{x_i^2}{2a_i^2}-\be x_i\right),
\end{equation}
where
$$\Upsilon^+=\{\widetilde{x}\in\mathbb{R}^n:\,
\dsum_{i=1}^k\frac{1}{2a_i^2}x_i^2\leq \dsum_{i=1}^k p_i\log 2\}
\cap\{\widetilde{x}\in\mathbb{R}^n:\, x_i\geq 0, 1\leq i\leq n\}.$$
Though we seem to have two different optimization problems in
(\ref{e2.5.14b}) and (\ref{e2.5.14c}), they produce the same
result as the former can be transform to the later using affine
transforms $y_i=\frac{x_i}{a_i}$. So this will lead to the expected
result that the limiting free energy of the Gaussian GREM is
continuous with respect to its weights. Not only in the Gaussian case
this can be made precise in all the models $(\ga>1)$ discussed in
this subsection by the same way and for other models with some extra
efforts.

In the Gaussian case, that is when $\ga=2$, two simple cases are
worth mentioning. The numbers $\be_j$ mentioned below are same as
the above, in these particular cases.

\begin{cor}\label{c2.5.2}{(Gaussian Case)}

i) Let
$0<\frac{p_1}{a_1^2}<\frac{p_2}{a_2^2}<\cdots<\frac{p_n}{a_n^2}$.
Put $\be_j=\frac{\sqrt{2p_j\log2}}{a_j}$ for $j=1,\cdots,n$. Then
a.s.
$$\begin{array}{lll}
\dlim_N \frac{1}{N}\log
Z_N(\be)&\hspace{-1ex}=\log2+\frac{\be^2}{2}\dsum_1^na_i^2 &\mbox{if
}
\be<\be_1,\\
&\hspace{-1ex}=\dsum_{j+1}^n p_i\log2+\dsum_{1}^j\be
a_i\sqrt{2p_i\log2}+\frac{\be^2}{2}\dsum_{j+1}^n a_i^2&\\
&&\hspace{-22ex}\mbox{if }\be_j\leq\be<\be_{j+1} \mbox{ for } 1\leq j<n,\\
&\hspace{-1ex}=\be \dsum_1^n a_i \sqrt{2p_i \log2}&\mbox{if
}\be\geq\be_n.
\end{array}$$
ii) Let
$\frac{p_1}{a_1^2}=\frac{p_2}{a_2^2}=\cdots=\frac{p_n}{a_n^2}>0$.
Then a.s.
$$\begin{array}{llll}
\dlim_N \frac{1}{N}\log
Z_N(\be)&=&\log2+\frac{\be^2}{2}\dsum_1^na_i^2 &\mbox{ if } \be<
\sqrt{\frac{2\log2}{\sum a_i^2}}\\
&=&\be\sqrt{2\log2 \sum a_i^2}&\mbox{ if } \be\geq
\sqrt{\frac{2\log2}{\sum a_i^2}}.
\end{array}$$
\end{cor}

\begin{rem}\label{r2.6.1}
We say that an $n$ level GREM with some particular driving
distribution is in {\em reduced form}, if the limiting free energy
of the model can not be obtained from any $k$ level GREM with same
driving distribution where $k<n$.

For a Gaussian $n$-level GREM, as the above analysis shows, if it
can not be obtained as a $k$-level GREM then $\be_i$s are defined
for $1\leq i\leq n$. On the other hand if it can be obtained as the
energy function of a $k$-level Gaussian GREM for some $k<n$, then
the $\be_i$s of the construction are only for $1\leq i\leq k$.

If a GREM is in reduced form, according to this definition, we do
not know whether its energy function can be obtain as that of a
$k$-level GREM for some $k<n$ with, of course, different driving
distributions. Along with the setup of the model in this subsection, we
are lucky enough to get the explicit expression of the limiting free
energy. Moreover we know the explicit expression of the $\be_i$s
where the expression of the free energy are changing. We observed in
this case that there may be at most $n$ many $\be_i$s. But we do not
know, whether this the intrinsic property of the model or there are
some driving distributions so that for an $n$-level GREM, we can get
more than $n$ many $\be_i$s.
\end{rem}

\begin{rem}
It can be shown that the energy function determines the parameters
of the model for every $\ga>1$ and one could characterize functions
those arise as energy functions for GREM. As observed in in the
above Corollary, an $n$ level GREM may reduce to a $k$ level GREM
for some $k<n$ or even to a REM. In such a case, some weights $a_i$
occur in groups and get added up. Of course, in such a case when the
model is not in reduced form, clearly it is not possible to recover
the weights from the formula for energy. But it is interesting to
note that when the GREM is in reduced form, we can recover the
parameters from the energy function. To make the statement precise
first of all note that, in this set up, GREM is in reduced from if
and only if all the $p_i$, $a_i$ are non zero and
$\frac{p_1^\frac{\ga-1}{\ga}}{a_1}<\frac{p_2^\frac{\ga-1}{\ga}}{a_2}
<\cdots<\frac{p_n^\frac{\ga-1}{\ga}}{a_n}$, let us assume this to be
the case. This is similar to that of the Gaussian case. Note that,
in this case $\be_i=\frac{(\ga p_i\log2)^\frac{\ga-1}{\ga}}{a_i}$
for $1\leq i\leq n$. From Theorem \ref{t2.6.1}, it follows that the
limiting free energy $\mathcal{E}(\be)$ is a continuous function
with $\mathcal{E}(0)=\log2$. It has a continuous derivative
$\mathcal{E}'(\be)$ with
\begin{equation*}
\mathcal{E}'(\be)=
\begin{cases}
0&\text{if $\be=0$}\\
\be^\frac{1}{\ga-1}\dsum_{i=k+1}^n a_i^{\frac{\ga}{\ga-1}} +
\dsum_{i=1}^k\be_i^{\frac{1}{\ga-1}} a_i^{\frac{\ga}{\ga-1}}
&\text{if $\be_k<\be \leq\be_{k+1},\;0\leq k\leq n-1$}\\
\dsum_{i=1}^n \be_i^{\frac{1}{\ga-1}}a_i^{\frac{\ga}{\ga-1}}&
\text{if $\be>\be_n$.}
\end{cases}
\end{equation*}
Further,
\begin{equation*}
\mathcal{E}''(\be) =\begin{cases}
\frac{1}{\ga-1}\be^\frac{2-\ga}{\ga-1}\dsum_1^n
a_i^\frac{\ga}{\ga-1},& \text{for $0< \be <\be_1$}\\
\frac{1}{\ga-1}\be^\frac{2-\ga}{\ga-1}\dsum_{k+1}^n
a_i^\frac{\ga}{\ga-1},& \text{for $\be_k<\be<\be_{k+1}$, for $1\leq
k<
n $}\\
0,& \text{for $\be>\be_n$}.\end{cases}\end{equation*}

The energy function can be characterized in this case. To start
with, observe that the above energy function has the following
properties:

i) $\mathcal{E}(0)=\log2$ and $\mathcal{E}'(0)=0$,

ii) $\mathcal{E}$ is a continuously differentiable function,

iii) denote $x_k=\frac{(\ga p_k\log2)^\frac{\ga-1}{\ga}}{a_k}$;
$c_k=\frac{1}{\ga-1}\be^\frac{2-\ga}{\ga-1}\dsum_k^n a_i^\frac{\ga}{\ga-1}$
and $\theta=\frac{2-\ga}{\ga-1}$ then
$0\equiv x_0<x_1<\cdots<x_n<x_{n+1}\equiv \infty$;
$c_1>c_2>\cdots>c_n>c_{n+1}\equiv 0$ and $\theta>-1$ with
$\mathcal{E}''(\be)=(1+\theta)c_i\be^\theta$ in $(x_{i-1},x_i)$ for
$1\leq i\leq n+1$.

Conversely, let $f$ be a function on $[0,\infty)$ such that

i)$f(0)=\log2$ and $f'(0)=0$,

ii) $f$ has continuous first derivative,

iii) there are finitely many points $0<x_1<\cdots<x_n$ and
$c_1>\cdots>c_n>c_{n+1}=0$ so that the left and right derivatives of $f'$ are
unequal at $x_i$ for $1\leq i\leq n$ and $f''(x)=x^\theta c_i$ in
$(x_{i-1},x_i)$ for $1\leq i\leq n+1$ with $x_0=0$ and $x_{n+1}=\infty$. Then $f$ is the
energy function for $\ga$-GREM with driving distribution having parameter
$\ga=\frac{\theta+2}{\theta+1}$, $p_i=\frac{\theta+1}{(\theta+2)\log2}(x_i(c_i- c_{i+1}))^{\theta+2}$
and $a_i=(c_i-c_{i+1})^\frac{1}{\theta+2}$ for $1\leq i\leq n$ if
\begin{equation}\label{e14a}
\dsum_{i=1}^n
x_i^{\theta+2} (c_i-c_{i+1}) =\frac{\theta+2}{\theta+1} \log2.
\end{equation}

In particular for a Gaussian GREM, it is in reduced from if and only
if all the $p_i$, $a_i$ are non zero and
$\frac{p_1}{a_1^2}<\frac{p_2}{a_2^2}<\cdots<\frac{p_n}{a_n^2}$. When
that is the case, from Theorem \ref{t2.6.1} with $\ga=2$, it follows
that the limiting free energy $\mathcal{E}(\be)$ is piecewise
quadratic continuous function with $\mathcal{E}(0)=\log2$. It has a
continuous derivative $\mathcal{E}'(\be)$ with $\mathcal{E}'(0)=0$
and
\begin{equation*}\mathcal{E}''(\be) =
\begin{cases}
\dsum_1^n a_i^2,& \text{in $(0,\frac{\sqrt{2p_1\log2}}{a_1})$,} \\
\dsum_{k+1}^n a_i^2,& \text{in
$(\frac{\sqrt{2p_k\log2}}{a_k},\frac{\sqrt{2p_{k+1}\log2}}{a_{k+1}})$, for $1\leq k \leq n-1$,}\\
0,& \text{if $\be>\frac{\sqrt{2p_n\log2}}{a_n}$.}
\end{cases}
\end{equation*}

Moreover, if $f$ is a $C^1$ function on $[0,\infty)$ with
$f(0)=\log2$ and $f'(0)=0$ so that there are finitely many points
$0<x_1<\cdots<x_n$ where the left and right derivatives of $f'$ are
unequal and $f''$ is a positive constant, say, $c_i$ in
$(x_{i-1},x_i)$ with $x_0=0$ and $x_{n+1}=\infty$. Then $f$ is the
energy function for some Gaussian GREM iff
\begin{equation}\label{e14}
c_1>\cdots>c_n>c_{n+1}=0 \quad \mbox{ and } \quad
\dsum_{i=1}^n x_i^2 (c_i-c_{i+1}) = 2\log2. \end{equation}
\end{rem}

\subsection{\texorpdfstring{$\ga=1$}{gamma=1}}

Now, let us assume $\ga=1$. Note that $\ga=1$ represents the two
sided exponential distribution with mean 0 and parameter 1. In this
case, we can not use the above argument directly as the ratios
$B(j,k)$ defined in (\ref{e2.5.a12}), the constants $a_i$ appear
with exponent $\frac{\ga}{\ga-1}$. However, to get the expression
for the free energy, we can directly proceed to evaluate
$$\dinf_{\tilde{x}\in\Psi^+}\dsum_{i=1}^n (1-\be a_i) x_i,$$ where
\[\Psi^+ = \left\{\widetilde{x}\in \mathbb{R}^n : \dsum_{i=1}^k
|x_i| \leq \dsum_{i=1}^k p_i \log2, \hspace{1ex} 1\leq k\leq
n\right\}\cap \{\tilde{x} : x_i \geq 0 \mbox{ for }1\leq i\leq
n\}.\] This is what we will do now.  To calculate the infimum, let
us set $r_0=0$ and for $k=1,2,\cdots$, let us define $\be_k$, $r_k$
as follows:
\begin{align*}
\be_1&=\min\{\frac{1}{a_i}: 1\leq i\leq n\}\\
r_1&=\max\{i: \frac{1}{a_i}=\be_1\}\\
\intertext{and in general, for $k>1$,}
\be_k &= \min\{\frac{1}{a_i}: r_{k-1} < i \leq n\}\\
r_k &=\max\{i:r_{k-1} < i \leq n,\; \frac{1}{a_i} = \be_k\}.
\end{align*}
Obviously this process stops at a finite stage say at $K$, so that
$\be_K=\frac{1}{a_n}$ and $r_K=n$. We put $\be_{K+1} = \infty$. For
example, if $a_1> a_2> \cdots >a_n$ then $\be_k = \frac{1}{a_k}$,
$r_k=k$ for $k= 1,2, \cdots,n$, and $K=n$. On the other hand if
$a_1< a_2< \cdots <a_n$ then $\be_1 = \frac{1}{a_n}$, $r_1=n$ and
$K=1$.

Clearly, $0=\be_0 <\be_1 < \be_2\cdots< \be_K <\be_{K+1}=\infty$.

\begin{rem}\label{r2.6.2}
The case $\ga=1$ can also be recovered as a limiting case
from the previous section. We can proceed by defining $\be_k$ as
done in the last subsection. But now we have to take limit
$\lim_{\ga\downarrow
1}\left(\frac{(p_j+\cdots+p_k)\ga\log2}{a_j^\frac{\ga}{\ga-1}+
\cdots+a_k^\frac{\ga}{\ga-1}}\right)^\frac{\ga-1}{\ga}$ to define
$B(j,k)$. A simple calculation shows that,
$B(j,k)=\frac{1}{\max_{j\leq i\leq k} a_i}$. Hence the $\be_k$s
defined in the earlier subsection lead to the same formula as above
when $\be\downarrow 1$.
\end{rem}

Now, fix $j\leq K$ and let  $\be\in[ \be_j, \be_{j+1}).$  Define
$\widetilde{\ol{x}}\in\Psi^+$ as follows:
\begin{equation}\label{e2.5.15a}
\begin{array}{llll}
\ol{x}_i &=& \dsum_{j=r_{(l-1)}+1}^{r_l}p_j\log2 &\mbox{ if } i=r_l
\mbox{ for some }l, 1\leq l \leq j\\
&=& 0 &\mbox{ otherwise. }
\end{array}\end{equation}

\noindent{\bf Claim:} $\dsty\inf_{\widetilde{x}\in \Psi^+}
\dsum_{i=1}^n(1-\be a_i)x_i$ occurs at $\widetilde{\ol{x}}.$\\

In case $j=0$ that is $\be_j=\be_0=0$, the claim is obvious. Indeed,
for $\be<\be_1$, $(1-\be a_i)$ is positive for all $i$, the infimum
occurs at $\widetilde{\ol{x}}$ with $\ol{x}_i=0$ for all $i$. So let
us assume that $j\geq 1$. First note that, since $\be\geq\be_j$, we
have $\be\geq\be_l$ for $1\leq l\leq j$ and
$(1-\frac{\be}{\be_l})=(1-\be a_{r_l})\leq 0$ . Moreover since
$\be_l$ are strictly increasing with $l$, the numbers $a_{r_l}$ are
strictly decreasing, that is, $a_{r_1}> a_{r_2}>\cdots>a_{r_K}$. Now
to prove the claim, fix any $\widetilde{x}\in \Psi^+$.

\begin{align*}\label{e2.5} &\dsum_{i=1}^{n}(1-\be a_i)
x_i-\dsum_{i=1}^{n}(1-\be a_i)
\ol{x}_i\\
\geq &\dsum_{i=1}^{r_j}(1-\be a_i) x_i-\dsum_{i=1}^{r_j}(1-\be a_i) \ol{x}_i\\
\intertext{(Since $(1-\be a_i)\geq 0$ and $\ol{x}_i=0$ for $i>r_j$)}
=&\dsum_{l=1}^j\left(\dsum_{i=r_{(l-1)}+1}^{r_l}(1-\be a_i) x_i-(1-\be a_{r_l})\ol{x}_{r_l}\right)\\
\intertext{(By definition of $\ol{x}_i$)}
\geq &\dsum_{l=1}^j(1-\be
a_{r_l})\left(\dsum_{i=r_{(l-1)}+1}^{r_l} x_i-\ol{x}_{r_l}\right)\\
\intertext{(Since by definition $a_{r_l}\geq a_i$ for
$r_{(l-1)}+1\leq i\leq r_l$) } = &(1-\be
a_{r_j})\left(\dsum_1^{r_j}x_i-\dsum_1^j
\ol{x}_{r_l}\right)+\dsum_{l=1}^{j-1}\be(a_{r_{(l+1)}}-
a_{r_l})\left(\dsum_1^{r_l} x_i-\dsum_1^l
\ol{x}_{r_i}\right)\\
= &(1-\be
a_{r_j})\dsum_1^{r_j}(x_i-p_i\log2)+\dsum_{l=1}^{j-1}\be(a_{r_{(l+1)}}-
a_{r_l})\dsum_1^{r_l}( x_i-p_i\log2)\\
\geq&0
\end{align*}

The last inequality follows from the facts that (i) by definition
$\dsum_1^{r_l}( x_i-p_i\log2)\leq 0$, (ii) $(a_{r_{(l+1)}}-
a_{r_l})< 0$ for $1\leq l\leq j$ and (iii) $(1-\be a_{r_j})\leq 0$.
Hence, the proof of the claim.

Here then is the formula for the free energy.
\begin{thm}\label{t2.5.3}
For two sided exponential GREM, almost surely,
\begin{equation*}
\dsty\lim_{N} \frac{1}{N} \log Z_N(\be)=
\begin{cases} \log2  & \text{if  $\be < \be_1$}\\
\log2+\dsum_{l=1}^j (\be a_{r_l} - 1)\dsum_{r_{l-1}+1}^{r_l}p_i
\log2 & \text{if $\be_j\leq \be < \be_{j+1}$.}
\end{cases}
\end{equation*}
\end{thm}

\begin{rem}\label{r2.6.3}
Once again, the free energy for the case $\ga=1$ can be recovered
from that of $\ga>1$ as a limiting case. It is quite easy to check
that with notation of $\be_l$ and $r_l$ as in subsection 2.6.1,
$\lim_{\ga\downarrow 1}(\be_l
a_i)^\frac{1}{\ga-1}=\frac{1}{k}\dsum_{j=r_{(l-1)}+1}^{r_l}p_j\log2$
where $k=\#\{i:\, a_i=\max_{r_{(l-1)}+1\leq i\leq r_l} a_i\, \&\,
r_{(l-1)}+1\leq i\leq r_l\}$. Moreover, for $\be a_i<1$, we have
$\lim_{\ga\downarrow 1}(\be a_i)^\frac{1}{\ga-1}=0$. So
$y_i=\lim_{\ga\downarrow 1}\ol{x}_i(\ga)$, where $\ol{x}_i$ as given
by (\ref{e2.5.12a}), may not give $\ol{x}_i$ as defined in
(\ref{e2.5.15a}). The only difference will be that $\ol{x}_i=0$ for
$r_{(l-1)}+1\leq i< r_l$ and
$\ol{x}_{r_l}=\dsum_{j=r_{(l-1)}+1}^{r_l}p_j\log2$ whereas, with the
same notation of $k$ as above,
$y_i=\frac{1}{k}\dsum_{j=r_{(l-1)}+1}^{r_l}p_j\log2$ for those $i$
where $r_{(l-1)}+1\leq i\leq r_l$ and $a_i=\max_{r_{(l-1)}+1\leq
i\leq r_l} a_i$. But it is easy to see that if
$\dsty\inf_{\widetilde{x}\in \Psi^+} \dsum_{i=1}^n(1-\be a_i)x_i$
occurs at $\widetilde{\ol{x}}$, then it will occur also at
$\widetilde{y}=(y_i)$. Thus the limiting free energy in the case of
$\ga=1$ is nothing but the limiting ($\ga\ra 1$) case of the
limiting free energy of the GREM where $\ga>1$.
\end{rem}
Thus Remarks \ref{r2.6.2} and \ref{r2.6.3} will lead to the
following
\begin{thm}
If $\mathcal{E}_\ga(\be)$ and $\mathcal{E}(\be)$ denote the limiting
free energy as given by Theorems \ref{t2.6.1} and \ref{t2.5.3}
respectively, then for all $\be\geq 0$, almost surely,
$$\lim_{\ga\downarrow 1}\mathcal{E}_\ga(\be)=\mathcal{E}(\be).$$
\end{thm}

As in the case of $\ga=2$, for the case $\ga=1$ also two special situations are worth mentioning
\begin{cor}
i) Let $a_1> a_2> \cdots >a_n$.  Then a.s.
\begin{equation*}
\dsty\lim_{N \ra \infty} \frac{1}{N} \log Z_N(\be)  =
\begin{cases} \log2& \text{if $\be < \frac{1}{a_1}$}\\
\log2+\dsum_{i=1}^k (\be a_i - 1)p_i \log2 & \text{if $\frac{1}{a_k}\leq \be < \frac{1}{a_{k+1}}$}\\
\be \dsum_{i=1}^n a_i p_i\log2             & \text{if $\be
\geq\frac{1}{a_n}$.}
\end{cases}
\end{equation*}

ii) Let $a_1\leq a_2 \leq \cdots \leq a_n$. Then a.s.
\begin{equation*}
\dsty\lim_{N \ra \infty} \frac{1}{N} \log Z_N(\be) =
\begin{cases}
\log2 &\text{if $\be < \frac{1}{a_n}$}\\
\be a_n\log2 &\text{if $\be \geq \frac{1}{a_n}$.}
\end{cases}
\end{equation*}
\end{cor}

\begin{rem}\label{r2.6.4}
Returning to Theorem \ref{t2.5.3}, it is interesting to note that
exponential GREM with parameters $(p_1,\cdots,p_n,a_1,\cdots,a_n)$
is equivalent to GREM with parameters
$(p_1',\cdots,p_K',a_1',\cdots,a_K')$ where $p_1'=\dsum_1^{r_1}
p_j,\; p_2'=\dsum_{r_1+1}^{r_2} p_j, \cdots,
p_K'=\dsum_{r_{(K-1)}+1}^{n} p_j$ and $a_1'=a_{r_1}, \;
a_2'=a_{r_2}, \cdots, a_K'=a_{r_K}.$ This is evident from Theorem
\ref{t2.5.3}. Here `equivalent' is used in the sense that for every
$\be$, both systems have the same free energy. Thus, in order that
an $n$-level GREM does not collapse to a lower level GREM it is
necessary and sufficient that the weights $a_i$ be strictly
decreasing. One should keep in mind that we are using the same
distribution at all levels of the GREM.
\end{rem}

The purpose of the following remark is to show that the energy
function determines the parameters of the model. One could
characterize functions that arise as energy functions for exponential GREM.
\begin{rem}
 As observed in the previous Remark, an $n$ level GREM may reduce
to a $K$ level GREM for some $K<n$. In the exponential GREM, some
weights $a_i$ do not appear in the formula for free energy. When
such a thing happens it is clearly not possible to recover the
weights from the formula for energy. It is interesting to note that
when the GREM is in reduced form, we can recover the parameters from
the energy function. Here is the precise statement.

Since an exponential GREM is in reduced form if and only
if $a_1>\cdots>a_n>0$ and $p_i\neq 0$ for $1\leq i\leq n$, let us
assume this to be the case. Let $\mathcal{E}(\be)$ be the energy
function, that is $\mathcal{E}(\be)=\dlim_N\frac{1}{N} \log
Z_N(\be)$. From Theorem \ref{t2.5.3}, it is easy to see that
$\mathcal{E}(\be)$ is a piecewise linear continuous function of
$\be$ taking value $\log2$ near zero. Further, its derivative
$\mathcal{E}'(\be)=\dsum_{i=1}^k a_ip_i\log2$ in
$(\frac{1}{a_k},\frac{1}{a_{k+1}})$. These properties are good
enough to show the following: $\mathcal{E}(\be)$ uniquely determines
all the quantities $p_i$ and $a_i$. In other words, the energy
function identifies the parameters.

If $0<x_1<\cdots<x_n$ be the points where the left and right
derivatives of $\mathcal{E}(\be)$ are unequal, then
$a_i=\frac{1}{x_i}$. Further, if $\mathcal{E}'(\be)=c_i$ in
$(x_i,x_{i+1})$ then $p_i= \frac{x_i(c_i-c_{i-1})}{\log2}$ for
$1\leq i\leq n$. Here $x_0=0$ and $x_{n+1}=\infty$.

In fact the above considerations lead to a characterization of
energy functions for exponential GREM. Suppose $f$ is a continuous
function on $[0,\infty)$ with $f(0)=\log2$. Further suppose that
there are finitely many points $0<x_1<\cdots<x_n$ where the left and
right derivatives are unequal and $f'$ is a constant, say, $c_i$ in
$(x_i,x_{i+1})$. Here $x_0=0$ and $x_{n+1}=\infty$. Then $f$ is the
energy function for some exponential GREM iff
\begin{equation}
0=c_0<c_1<\cdots<c_n \quad \mbox{ and } \quad \dsum_{i=1}^n
x_i(c_i-c_{i-1})=\log2.
\end{equation}
\end{rem}

\subsection{\texorpdfstring{$0<\ga<1$}{0<gamma<1}}

Now we come to the case $\ga<1$. Unlike in the above two
subsections, here we have not been able to derive the closed form
expression of the free energy for general $n$ level trees. For
$\ga<1$, the function $\dsum_{i=1}^n \left( \frac{x_i^\ga}{\ga}-\be
a_i x_i\right)$ is not a convex function, rather a concave function.
Moreover the domain $\Psi=\{\widetilde{x}\in\mbb{R}^n:\,
\sum_{i=1}^k|x_i|^\ga \leq \sum_{i=1}^k\ga p_i\log2,\quad 1\leq
k\leq n\}$ is also a non-convex set. Hence in order to calculate
\begin{equation}
 \dinf_{\tilde{x}\in\Psi^+}\dsum_{i=1}^n \left(
\frac{x_i^\ga}{\ga}-\be a_i x_i\right)
\end{equation}
with
\begin{equation}
 \Psi^+=\{\widetilde{x}\in\mbb{R}^n:\, \sum_{i=1}^k|x_i|^\ga \leq
\sum_{i=1}^k\ga p_i\log2, \&\; x_k\geq0 \text{ for } 1\leq k\leq n\},
\end{equation}
we can not use the convex analysis, as we did for the case $\ga>1$.
However, by change of variables, the problem can be brought back to
optimizing a convex function over a convex set. Now we specialize to
the case $n=2$. We shall calculate
\begin{equation}
\dinf_{\tilde{x}\in\Psi^+}\left\{(\frac{1}{\ga}x_1^\ga-\be a_1 x_1)+(\frac{1}{\ga}x_2^\ga-\be a_2 x_2)\right\}
\end{equation}
with
\begin{equation}
\Psi^+=\{(x_1,x_2)\geq 0:\, x_1^\ga \leq \ga p_1\log2, x_1^\ga+x_2^\ga \leq \ga (p_1+p_2)\log2\}.
\end{equation}
To do so, we transform the problem by denoting
$\frac{1}{\ga}x_1^\ga=x$, $\frac{1}{\ga}x_2^\ga=y$, $a_1
\ga^\frac{1}{\ga}=a$, $a_2 \ga^\frac{1}{\ga}=b$, $p_1\log2=c$,
$p_2\log2=d$ and $\al=\frac{1}{\ga}$ so that, $\al>1$ and we need to
calculate
\begin{equation}
-\sup_{\tilde{x}\in\Psi^+}\left\{(\be a x^\al-x)+(\be b y^\al-y)\right\}
\end{equation}
with
\begin{equation}
\Psi^+=\{(x,y)\geq 0:\, x \leq c, x+y \leq c+d\}.
\end{equation}
Let $f(x,y)=(\be a x^\al-x)+(\be b y^\al-y)$. Since $f(x,y)$ is a
convex function and we are looking for supremum over a convex set,
the supremum occurs at the boundary points. Note that $\Psi^+$ is a
polygon. Where as for any $\ep\in \mathbb{R}$, the set
$\{f(x,y)=\ep\}$ is either empty set or a smooth curve. Hence the
above supremum occurs at one of the corner points, $A\equiv(0,0)$,
$B\equiv(c,0)$, $C\equiv(c,d)$ and $D\equiv(0,c+d)$, of $\Psi^+$.
Now
\begin{align}
f(A)\gtreqqless f(B)\quad& \text{iff $\be \lesseqqgtr \frac{1}{ac^{\al-1}}$},\\
f(A)\gtreqqless f(C)\quad& \text{iff $\be \lesseqqgtr \frac{c+d}{ac^\al+bd^\al}$},\\
f(A)\gtreqqless f(D)\quad& \text{iff $\be \lesseqqgtr \frac{1}{b(c+d)^{\al-1}}$},\\
f(B)\gtreqqless f(C)\quad& \text{iff $\be \lesseqqgtr \frac{1}{bd^{\al-1}}$},\\
f(B)\gtreqqless f(D)\quad& \text{iff $b(c+d)^\al >ac^\al$  and $\be \lesseqqgtr \frac{d}{b(c+d)^\al-ac^\al}$},\\
f(B)> f(D)\quad& \text{if $b(c+d)^\al \leq ac^\al$},\\
f(C)\gtreqqless f(D)\quad& \text{iff $ac^\al+bd^\al\gtreqqless b(c+d)^\al.$}
\end{align}
Note that the last two relations do not depend on $\be$. Now comparing all
the possibilities, we obtain the following three scenarios:\\

\noindent\shadowbox{$\boldsymbol{b(c+d)^\al\leq ac^\al+bd^\al}$}

Let us assume $b(c+d)^\al\leq ac^\al+bd^\al$. Then it is easy to see
that
\begin{equation*}
b(c+d)^\al\leq ac^\al+bd^\al \Rightarrow
\begin{cases}
f(C)\geq f(D),\\
\frac{c+d}{ac^\al+bd^\al}\leq\frac{1}{b(c+d)^{\al-1}}<\frac{1}{bd^{\al-1}}.
\end{cases}
\end{equation*}
Now
$$\frac{c+d}{ac^\al+bd^\al}<\frac{1}{bd^{\al-1}}\Ra \frac{1}{ac^{\al-1}}<\frac{1}{bd^{\al-1}},$$
and $$\frac{1}{ac^{\al-1}}\lesseqqgtr
\frac{c+d}{ac^\al+bd^\al}\Leftrightarrow
\frac{1}{ac^{\al-1}}\lesseqqgtr\frac{1}{bd^{\al-1}}$$ implies
$$\frac{1}{ac^{\al-1}}<\frac{c+d}{ac^\al+bd^\al}\leq\frac{1}{b(c+d)^{\al-1}}<\frac{1}{bd^{\al-1}}.$$
Hence we get
\begin{equation}
\sup_{\Psi^+}f(x,y)=\begin{cases}
f(A)\quad \text{if $0\leq \be\leq \frac{1}{ac^{\al-1}}$}\\
f(B)\quad \text{if $\frac{1}{ac^{\al-1}}\leq \be \leq
\frac{1}{bd^{\al-1}}$}\\
f(C)\quad \text{if $\be \geq \frac{1}{bd^{\al-1}}$}.
\end{cases}
\end{equation}

\noindent\shadowbox{$\boldsymbol{b(c+d)^\al>
ac^\al+bd^\al\;\&\;b(c+d)^{\al-1}\leq ac^{\al-1}}$}

Let us
assume $b(c+d)^\al> ac^\al+bd^\al\;\&\;b(c+d)^{\al-1}\leq
ac^{\al-1}$. Then it is easy to see that
\begin{equation*}
b(c+d)^\al> ac^\al+bd^\al \Rightarrow
\begin{cases}
f(D)> f(C),\\
\frac{1}{b(c+d)^{\al-1}}<\frac{c+d}{ac^\al+bd^\al},\\
\frac{d}{b(c+d)^\al-ac^\al}<\frac{1}{bd^{\al-1}}.
\end{cases}
\end{equation*}
Moreover,
\begin{equation*}
b(c+d)^{\al-1}\leq ac^{\al-1} \Rightarrow
\frac{1}{ac^{\al-1}}\leq\frac{1}{b(c+d)^{\al-1}}\leq\frac{d}{b(c+d)^\al-ac^\al}.
\end{equation*}
Thus
$$\frac{1}{ac^{\al-1}}\leq\frac{1}{b(c+d)^{\al-1}}\leq\frac{d}{b(c+d)^\al-ac^\al}<\frac{1}{bd^{\al-1}}.$$
Hence we get
\begin{equation}
\sup_{\Psi^+}f(x,y)=\begin{cases}
f(A)\quad \text{if $0\leq \be\leq \frac{1}{ac^{\al-1}}$}\\
f(B)\quad \text{if $\frac{1}{ac^{\al-1}}\leq \be \leq
\frac{d}{b(c+d)^\al-ac^\al}$}\\
f(D)\quad \text{if $\be \geq \frac{d}{b(c+d)^\al-ac^\al}$}.
\end{cases}
\end{equation}

\noindent\shadowbox{$\boldsymbol{b(c+d)^\al>
ac^\al+bd^\al\;\&\;b(c+d)^{\al-1}> ac^{\al-1}}$}

Let us assume
$b(c+d)^\al> ac^\al+bd^\al\;\&\;b(c+d)^{\al-1}> ac^{\al-1}$. Then it
is easy to see that
\begin{equation*}
b(c+d)^\al> ac^\al+bd^\al \Rightarrow
\begin{cases}
f(D)> f(C),\\
\frac{1}{b(c+d)^{\al-1}}<\frac{c+d}{ac^\al+bd^\al}.
\end{cases}
\end{equation*}
Moreover,
\begin{equation*}
b(c+d)^{\al-1}> ac^{\al-1} \Rightarrow
\frac{1}{ac^{\al-1}}>\frac{1}{b(c+d)^{\al-1}}>\frac{d}{b(c+d)^\al-ac^\al}.
\end{equation*}
Thus
$$\frac{d}{b(c+d)^\al-ac^\al}<\frac{1}{b(c+d)^{\al-1}}<\frac{c+d}{ac^\al+bd^\al}.$$
Hence we get
\begin{equation}
\sup_{\Psi^+}f(x,y)=\begin{cases}
f(A)\quad \text{if $0\leq \be\leq \frac{1}{b(c+d)^{\al-1}}$}\\
f(D)\quad \text{if $\be \geq \frac{1}{b(c+d)^{\al-1}}$}.
\end{cases}
\end{equation}

We can conclude the above three cases in the following:

\begin{thm}
For two level GREM with driving distribution having density
$\phi_{N,\ga}$ as defined in (\ref{e2.6.a12}) with $0<\ga<1$, we
have, almost surely,
\begin{enumerate}
\item{if $a_2\leq a_1(p_1\log2)^\frac{1}{\ga}+a_2(p_2\log2)^\frac{1}{\ga}$,
then the limiting free energy is

$\begin{cases} \log2&\text{for
$0\leq \be\leq
\frac{1}{a_1\ga^\frac{1}{\ga}(p_1\log2)^\frac{1-\ga}{\ga}}$,}\\
p_2\log2+\be a_1 (\ga p_1\log2)^\frac{1}{\ga}& \text{for
$\frac{1}{a_1\ga^\frac{1}{\ga}(p_1\log2)^\frac{1-\ga}{\ga}}\leq
\be\leq
\frac{1}{a_2\ga^\frac{1}{\ga}(p_2\log2)^\frac{1-\ga}{\ga}}$,}\\
\be (a_1 (\ga p_1\log2)^\frac{1}{\ga}+a_2 (\ga
p_2\log2)^\frac{1}{\ga}) &\text{for $\be \geq
\frac{1}{a_2\ga^\frac{1}{\ga}(p_2\log2)^\frac{1-\ga}{\ga}}.$}
\end{cases}$}
\item{if $a_2> a_1(p_1\log2)^\frac{1}{\ga}+a_2(p_2\log2)^\frac{1}{\ga}$
and $a_2\leq a_1(p_1\log2)^\frac{1-\ga}{\ga}$, then the limiting
free energy is

$\begin{cases}
\log2&\text{for $0\leq \be\leq
\frac{1}{a_1\ga^\frac{1}{\ga}(p_1\log2)^\frac{1-\ga}{\ga}}$,}\\
p_2\log2+\be a_1 (\ga p_1\log2)^\frac{1}{\ga}& \text{for
$\frac{1}{a_1\ga^\frac{1}{\ga}(p_1\log2)^\frac{1-\ga}{\ga}}\leq
\be\leq
\frac{p_2\log2}{a_2(\ga\log2)^\frac{1}{\ga}-a_1(\ga p_1\log2)^\frac{1}{\ga}}$,}\\
\be a_2 (\ga\log2)^\frac{1}{\ga}) &\text{for $\be \geq
\frac{p_2\log2}{a_2(\ga\log2)^\frac{1}{\ga}-a_1(\ga
p_1\log2)^\frac{1}{\ga}}.$}
\end{cases}$}
\item{if $a_2> a_1(p_1\log2)^\frac{1}{\ga}+a_2(p_2\log2)^\frac{1}{\ga}$
and $a_2> a_1(p_1\log2)^\frac{1-\ga}{\ga}$, then the limiting free
energy is

$\begin{cases}
\log2&\text{for $0\leq \be\leq
\frac{1}{a_2\ga^\frac{1}{\ga}(\log2)^\frac{1-\ga}{\ga}}$,}\\
\be a_2 (\ga\log2)^\frac{1}{\ga})  &\text{for $\be \geq
\frac{1}{a_2\ga^\frac{1}{\ga}(\log2)^\frac{1-\ga}{\ga}}.$}
\end{cases}$}
\end{enumerate}
\end{thm}

\begin{rem}
Note that in a $2$ level double exponential GREM (in the earlier subsection) with weights $a_1$
and $a_2$, we had at most two cases, namely, $\frac{1}{a_1}\leq
\frac{1}{a_2}$ and $\frac{1}{a_1}> \frac{1}{a_2}$. Where as for
$\ga<1$, we have three cases.

Though there are three cases, we can think of them as two cases like
the double exponential GREM, namely,
$\frac{1}{a_1(p_1\log2)^\frac{1-\ga}{\ga}}\leq
\frac{1}{a_2(\log2)^\frac{1-\ga}{\ga}}$ and
$\frac{1}{a_1(p_1\log2)^\frac{1-\ga}{\ga}}>
\frac{1}{a_2(\log2)^\frac{1-\ga}{\ga}}$, where the first case has
two more subcases, namely,
$\frac{1}{a_2(p_2\log2)^\frac{1-\ga}{\ga}}\leq\frac{p_2\log2}{a_2(\log2)^\frac{1}{\ga}-a_1(p_1\log2)^\frac{1}{\ga}}$
and
$\frac{1}{a_2(p_2\log2)^\frac{1-\ga}{\ga}}>\frac{p_2\log2}{a_2(\log2)^\frac{1}{\ga}-a_1(p_1\log2)^\frac{1}{\ga}}$.
\end{rem}

\medskip
\section{Level-dependant Distributions}\label{s2.7}

We already mentioned that the LDP holds good even when the driving
distributions at various levels are different. To be precise, fix
numbers $\ga_1,\cdots,\ga_n$; each greater than 0. Consider an $n$
level GREM with the driving distribution at the $i$-th level being
$\phi_{N,\ga_i}$ given by (\ref{e2.6.a12}). That is, at the first
level for each edge $\si_1$ the associated random variable
$\xi(\si_1)$ has density $\phi_{N,\ga_1}$. In general, for any edge
$\si_1\cdots\si_i$ at the $i$-th level, the associated random
variables $\xi(\si_1\cdots\si_i)$ has density $\phi_{N,\ga_i}$. Then
the map $\Si_N \ra \mathbb{R}^n$ by
\[\si \mapsto \left(\frac{\xi(\si_1,\om)}{N}, \frac{\xi(\si_1\si_2,\om)}{N}, \cdots,
\frac{\xi(\si_1\cdots\si_n,\om)}{N}\right)\] induce random
probability $\mu_N(\om)$ on $\mathbb{R}^n$ by transporting the
uniform probability on $\Si_N$. Theorem \ref{t2.4.6} suggest that in
this case the free energy of the system will be
\begin{equation}\label{e2.7.16}
\dsty\lim_{N} \frac{1}{N} \log
Z_N(\be)=\log2-\dinf_{\widetilde{x}\in\Psi^+}
\left\{\dsum_{i=1}^n\left(\frac{x_i^{\ga_i}}{\ga_i}-\be
a_i x_i\right)\right\},
\end{equation}

where $\Psi^+$ is the intersection of
\[\Psi = \{\widetilde{x}\in \mathbb{R}^n : \dsum_{i=1}^k
\frac{|x_i|^{\ga_i}}{\ga_i} \leq \dsum_{i=1}^k p_i \log2,
\hspace{1ex} 1\leq k\leq n\},\] with the positive orthant of
$\mathbb{R}^n$  and $p_i=\dlim_{N\ra\infty} \frac{k(i,N)}{N}$.

In its generality, it is very difficult to have a closed form
expression for the above infimum. May be there is no general closed
form expression, for the infimum and hence for the free energy of
the system. To make a beginning and to see what one can expect, we
now specialize to the case $n=2$. The limiting frequencies $\dlim_N
\frac{k(i,N)}{N}$ are $p_i$ for $i=1, 2$. The weights for the two
level are $a_1$ and $a_2$ respectively. We assume $p_1, p_2, a_1,
a_2$ are strictly
positive.\\

\subsection{Exponential - Gaussian GREM}

In this case we consider the distributions at the first level to be
$\phi_{N,1}$ and at the second level to be $\phi_{N,2}$ --- that is,
exponential and Gaussian respectively. So from (\ref{e2.7.16}), the
expression for the free energy for this case will read as follows:
\begin{align}
\mathcal{E}(\be)&=\dsty\lim_{N} \frac{1}{N} \log Z_N(\be)\notag\\
&=\log2-\inf \{f(x,y): x,y\geq 0;\; x\leq p_1\log2;\;
x+\frac{1}{2}y^2\leq \log2\}\\
\intertext{where}
f(x,y)&= x(1-\be a_1)+\frac{1}{2}y^2-\be a_2 y.\label{e2.7.a17}
\end{align}

To calculate $\mathcal{E}(\be)$ explicitly we proceed as follows.
First we discuss the case $\be\leq\frac{1}{a_1}$. Then we discuss
$\be>\frac{1}{a_1}$. This later case leads to three subcases. In
each subcase combining the conclusion along with the case
$\be\leq\frac{1}{a_1}$, we give a full picture of $\mathcal{E}(\be)$
for all values of $\be$.\\

\noindent\doublebox{\bf I.\; $\boldsymbol{\be\leq\frac{1}{a_1}}$}\vspace{1ex}

On the interval $[0,\infty)$, the function $\frac{1}{2}y^2-\be a_2 y$
decreases up to $\be a_2$ and then increases. So when $\be\leq
\frac{1}{a_1}$, that is when $1-\be a_1\geq 0$, the above function
attains its minimum at the point, $(0,\be a_2\wedge \sqrt{2\log2})$.\\

\noindent\doublebox{\bf II.\; $\boldsymbol{\be>\frac{1}{a_1}}$}\vspace{1ex}

If $\be>\frac{1}{a_1}$, here is how to calculate the infimum. The
function $g(y)= \dinf_x f(x,y)$ is given by
\begin{equation*}
g(y) =\begin{cases}p_1(1-\be a_1)\log2+ \frac{1}{2}y^2 - \be a_2 y & \text{for $0\leq y\leq\sqrt{2p_2\log2}$}\\
(1-\be a_1)\log2+ \frac{1}{2}\be a_1 y^2 - \be a_2 y & \text{for $\sqrt{2p_2\log 2}\leq y \leq\sqrt{2\log2}$.}
\end{cases}
\end{equation*}
This is because, when $0\leq y\leq \sqrt{2p_2\log2}$, $\dinf_x
f(x,y)$ is attained at $x=p_1\log2$, whereas in the other case the
infimum is attained at $x=\log2-\frac{1}{2}y^2$.

Since the required infimum of $f$ is just the infimum of $g(y)$, one
has to calculate $\dinf_{0\leq y\leq \sqrt{2\log 2}} g(y)$  by
analyzing $g$ in the two intervals separately. This is what we do below. First note that
the function $g$ is continuous. Now we have the following three
scenarios.\\

\noindent\shadowbox{\bf A1: $\boldsymbol{\frac{a_2}{a_1} < \sqrt{2p_2\log2} }$}\vspace{1ex}

Let us assume $\frac{a_2}{a_1} < \sqrt{2p_2\log2}$. First let us
consider $\be$ such that $\frac{1}{a_1}<\be\leq
\frac{\sqrt{2p_2\log2}}{a_2}$. In particular, $\be
a_2\leq\sqrt{2p_2\log2}$ where as $\frac{a_2}{a_1}
<\sqrt{2p_2\log2}$. So the function $\frac{1}{2}y^2 - \be a_2 y$ is
decreasing up to $\be a_2$ in $[0, \sqrt{2p_2\log2}]$ and then
increasing. Thus in $[0, \sqrt{2p_2\log2}]$, $g$ attains its minimum
at $\be a_2$. On $[\sqrt{2p_2\log 2}, \sqrt{2\log2}]$ the function
$\frac{1}{2}\be a_1 y^2 - \be a_2 y=\be a_1(\frac{1}{2}y^2 -
\frac{a_2}{a_1} y)$ is increasing. Hence, $g$ being continuous, for
the values of $\be$ under consideration, the infimum will occur at
$y=\be a_2$.

Now $\be$ be such that $\be>\frac{\sqrt{2p_2\log2}}{a_2}$ so that
$\be a_2> \sqrt{2p_2\log2}$. Since
$\frac{a_2}{a_1}<\sqrt{2p_2\log2}$ the function $\frac{1}{2}\be a_1
y^2 - \be a_2 y$ is increasing on $[\sqrt{2p_2\log 2},
\sqrt{2\log2}]$. The function $\frac{1}{2}y^2 - \be a_2 y$ is
decreasing on $[0, \sqrt{2p_2\log2}]$ attaining infimum at
$y=\sqrt{2p_2\log2}$. As a consequence, for
$\be>\frac{\sqrt{2p_2\log2}}{a_2}$, the infimum of $g(y)$ is occurs
at $\sqrt{2p_2\log2}$.

Thus combining {\bf I.} and above para we conclude that if
$\frac{a_2}{a_1} < \sqrt{2p_2\log2}$ then phase transitions take
place at $\be = \frac{1}{a_1}$ and
$\be=\frac{\sqrt{2p_2\log2}}{a_2}$. So, substituting this
corresponding arguments where minimum is attained in
(\ref{e2.7.a17}), we have the following

\begin{thm}
In the Exponential-Gaussian GREM, if $\frac{a_2}{a_1} <
\sqrt{2p_2\log2}$ then almost surely,
\begin{equation*}
\dlim_{N\ra\infty} \frac{1}{N} \log Z_N(\be) =
\begin{cases}\log2+\frac{1}{2}\be^2a_2^2  & \text{if $\be\leq\frac{1}{a_1}$} \\
p_2\log2+\frac{1}{2}\be^2a_2^2+\be p_1a_1\log2 & \text{if $\frac{1}{a_1}<\be\leq \frac{\sqrt{2p_2\log2}}{a_2}$}\\
\be(a_2\sqrt{2p_2\log2}+a_1p_1\log2) & \text{if
$\be>\frac{\sqrt{2p_2\log2}}{a_2}$.} \end{cases}
\end{equation*}
\end{thm}
We can picture the value of $\mathcal{E}(\be)$ against $\be$ as
given below. The values of $\be$ are given under the line and values
of $\mathcal{E}(\be)$ are given above the line. The phase
transitions occur at the dark lines.

\vspace{3ex}
\begin{texdraw}
\htext (2 .1){\bf Subcase A1}

\drawdim in

\linewd 0.01

\move(0 -.5) \avec(5.4 -.5)

\move(0 -0.6) \lvec(0 -0.4)

\htext (0 -.75){0}

\htext (-.1 -.2){$\boldsymbol{\mathcal{E}(\be) \rightarrow}$}

\htext (-.1 -.95){$\boldsymbol{\beta \rightarrow}$}

\htext (.05 -.4){$\frac{1}{2}\be^2a_2^2+\log2$}

\linewd 0.02 \move(3.09 -.65) \lvec(3.09 -.3)

\htext (2.8 -.9){$\frac{\sqrt{2p_2\log2}}{a_2}$}

\htext (3.31 -.4){$\be(a_2\sqrt{2p_2\log2}+a_1p_1\log2)$}

\move(1.02 -.65) \lvec(1.02 -.3)

\htext (.95 -.9){$\frac{1}{a_1}$}

\htext (1.22 -.4){$\frac{1}{2}\be^2a_2^2+(\be p_1a_1+p_2)\log2$}

\linewd 0.01 \move(4.5 -.6) \lvec(4.5 -.5)

\htext (4.3 -.9){$\frac{\sqrt{2\log2}}{a_2}$}
\end{texdraw}

This case seems rather peculiar. This is indeed a sum of two REMs,
as follows. Imagine placing exponential random variables
$\xi_{\si_1}$ at the first level and one i.i.d bunch
$\{\xi_{\si_1\si_2}\}$ is placed below each first level node. In
other words, consider $\{\eta_{\si_2}:\si_2\in 2^{k(2,N)}\}$ i.i.d
$\mathcal{N}(0,N)$ and set $\xi_{\si_1\si_2}=\eta_{\si_2}$ for all
$\si_1, \si_2$. Consider the corresponding Hamiltonian
$H_N(\si)=a_1\xi_{\si_1}+a_2\xi_{\si_1\si_2}=a_1\xi_{\si_1}+a_2\eta_{\si_2}$.
Let us set $Z_N^1=\dsum_{\si_1}e^{\be a_1\xi_{\si_1}}$, the
partition function for the $k(1,N)$-particles system consisting of
exponential Hamiltonian with weight $a_1$. Let
$Z_N^2=\dsum_{\si_2}e^{a_2\eta_{\si_2}}$, the partition function for
$k(2,N)$ particle system consisting of Gaussian, $\mathcal{N}(0,N)$
Hamiltonian with weight $a_2$. Clearly, $Z_N=Z_N^1\cdot Z_N^2$. If,
for $i=1,2;\; \mathcal{E}_i=\dlim_N \frac{1}{N} \log Z_N^i$ then the
exponential REM formula \cite{N1,JR1} yields, a.s.,
\begin{equation}\label{e2.7.17}
\mathcal{E}_1(\be)=
\begin{cases}p_1\log2 &\text{ if $\be\leq \frac{1}{a_1}$}\\
\be p_1 a_1\log2 &\text{ if $\be > \frac{1}{a_1}$.}
\end{cases}
\end{equation}

The Gaussian REM formula (keeping in mind that for $N$ fixed, the
$k(2,N)$ particle system has $\mathcal{N}(0,N)$ Hamiltonians as
opposed to $\mathcal{N}(0,k(2,N))$ yields, a.s,
\begin{equation}\label{e2.7.18}
\mathcal{E}_2(\be)=
\begin{cases}
p_2\log2 +\frac{1}{2} a_2^2 \be^2 &\text{ if $\be\leq \frac{\sqrt{2p_2\log2}}{a_2}$}\\
\be a_2\sqrt{2p_2\log2} &\text{ if $ \be >
\frac{\sqrt{2p_2\log2}}{a_2}$.}
\end{cases}
\end{equation}

One can now verify that, a.s.
$$\mathcal{E}(\be)=\mathcal{E}_1(\be)+\mathcal{E}_2(\be).$$
In other words the GREM behaves like sum of two independent REMs,
one exponential and other Gaussian. The word independent is used
here in the sense that there is no interaction between these two
REMs -- that is, there is no interaction between the $k(1,N)$
particles and the $k(2,N)$ particles, as if there is a barrier
between these two sets of particles. Of course, this is so as long
as $\frac{a_2}{a_1} < \sqrt{2p_2\log2}$.\\

\noindent\shadowbox{\bf A2: $\boldsymbol{\sqrt{2p_2\log2}\leq \frac{a_2}{a_1} < \sqrt{2\log2} }$}\vspace{1ex}

Let us assume $\sqrt{2p_2\log2}\leq \frac{a_2}{a_1} <
\sqrt{2\log2}$. Then $\be> \frac{1}{a_1}$ means $\be
a_2>\frac{a_2}{a_1}\geq\sqrt{2p_2\log2}$ where as $\frac{a_2}{a_1}
<\sqrt{2\log2}$. So the function $\frac{1}{2}y^2 - \be a_2 y$ is
decreasing on $[0, \sqrt{2p_2\log2}]$ and the other function
$\frac{1}{2}\be a_1 y^2 - \be a_2 y=\be a_1(\frac{1}{2}y^2 -
\frac{a_2}{a_1} y)$ is decreasing up to $\frac{a_2}{a_1}$ in
$[\sqrt{2p_2\log 2}, \leq\sqrt{2\log2}]$ and then increasing. Hence,
as $g$ is continuous, the infimum will occur at $y=\frac{a_2}{a_1}$.
Thus, the phase transition takes place at $\be = \frac{1}{a_1}$. So
we have the following

\begin{thm}
In the Exponential-Gaussian GREM, if $\sqrt{2p_2\log2}\leq
\frac{a_2}{a_1} < \sqrt{2\log2}$ then almost surely,
\begin{equation*}
\dlim_{N\ra\infty} \frac{1}{N} \log Z_N(\be) =
\begin{cases}
\log2+\frac{1}{2}\be^2a_2^2 & \mbox{if }\;\be\leq\frac{1}{a_1} \\
\be \left(\frac{1}{2}\frac{a_2^2}{a_1}+a_1\log2\right) & \mbox{if
}\;\be>\frac{1}{a_1}
\end{cases}
\end{equation*}
\end{thm}
As earlier, we can picture the value of $\mathcal{E}(\be)$ against
$\be$ as given below. The values of $\be$ are given under the line
and values of $\mathcal{E}(\be)$ are given above the line. The phase
transitions occur at the dark lines.

\vspace{10ex}
\begin{texdraw}
\htext (2 .1){\bf Subcase A2} \drawdim in

\linewd 0.01

\move(0 -.5) \avec(5.3 -.5)

\move(0 -0.6) \lvec(0 -0.4)

\htext (0 -.75){0}

\htext (-.1 -.2){$\boldsymbol{\mathcal{E}(\be) \rightarrow}$}

\htext (-.1 -.95){$\boldsymbol{\beta \rightarrow}$}

\htext (.5 -.4){$\frac{1}{2}\be^2a_2^2+\log2$}

\linewd 0.02 \move(2.5 -.65) \lvec(2.5 -.35)

\htext (2.44 -.9){$\frac{1}{a_1}$}

\htext (2.75 -.4){$\be(\frac{1}{2} \frac{a_2^2}{a_1}+a_1\log2)$}
\linewd 0.01 \move(3.5 -.6) \lvec(3.5 -.5)

\htext (3.3 -.9){$\frac{\sqrt{2\log2}}{a_2}$}

\move(1.75 -.6) \lvec(1.75 -.5)

\htext (1.5 -.9){$\frac{\sqrt{2p_2\log2}}{a_2}$}
\end{texdraw}

In this case, we observe that the free energy for inverse
temperature up to $\frac{1}{a_1}$ is given by $\log2 +\frac{1}{2}
\be^2a_2^2$. This can be thought of as the Gaussian REM energy but
not going all the way up to $\be\leq \frac{\sqrt{2\log2}}{a_2}$ but
cut short at $\frac{1}{a_1}$. This can also be thought of as the sum
of the two energies $\mathcal{E}_1$ and $\mathcal{E}_2$ as in
(\ref{e2.7.17}) and (\ref{e2.7.18}), but then the Gaussian effect is
prolonged up to $\be\leq \frac{1}{a_1}$ instead of stopping at
$\frac{\sqrt{2p_2\log2}}{a_2}$. We do not know which is the correct
interpretation. For $\be>\frac{1}{a_1}$, the system exhibits a new
phenomenon which we are unable to explain. The term $\be a_1\log2$
is reminiscent of the exponential REM energy. The other term
$\frac{1}{2}\be\frac{a_2^2}{a_1}$ appears to be new.\\

\noindent\shadowbox{\bf A3: $\boldsymbol{\sqrt{2\log2} \leq \frac{a_2}{a_1}}$}\vspace{1ex}

Let us assume $\sqrt{2\log2} \leq \frac{a_2}{a_1}$. Then $\be>
\frac{1}{a_1}$ means $\be a_2> \frac{a_2}{a_1}\leq \sqrt{2\log2}$.
So both the functions $\frac{1}{2}y^2 - \be a_2 y$ and
$\frac{1}{2}\be a_1 y^2 - \be a_2 y=\be a_1(\frac{1}{2}y^2 -
\frac{a_2}{a_1} y)$ are decreasing on $[0, \sqrt{2p_2\log2}]$ and
$[\sqrt{2p_2\log 2}, \leq\sqrt{2\log2}]$ respectively. Hence the
infimum will occur at $y= \sqrt{2\log2}$. Being $\be a_2>
\sqrt{2\log2}$, the phase transition takes place at $\be =
\frac{\sqrt{2\log2}}{a_2}$. Hence we have the following

\begin{thm}
In the Exponential-Gaussian GREM, if $\frac{a_2}{a_1}\geq
\sqrt{2\log2}$ then almost surely,
\begin{equation*}
\dlim_{N\ra\infty} \frac{1}{N} \log Z_N(\be)=
\begin{cases}\log2+\frac{1}{2}\be^2a_2^2&\text{if $\be\leq\frac{\sqrt{2\log2}}{a_2}$}\\
\be a_2\sqrt{2\log2}&\text{if $\be>\frac{\sqrt{2\log2}}{a_2}$}
\end{cases}
\end{equation*}
\end{thm}

As earlier, we can picture the value of $\mathcal{E}(\be)$ against
$\be$ as given below. The values of $\be$ are given under the line
and values of $\mathcal{E}(\be)$ are given above the line. The phase
transitions occur at the dark lines.

\vspace{3ex}
\begin{texdraw}
\htext (2 .1){\bf Subcase A3} \drawdim in

\linewd 0.01

\move(0 -.5) \avec(5.3 -.5)

\move(0 -0.6) \lvec(0 -0.4)

\htext (1 -.4){$\frac{1}{2}\be^2a_2^2+\log2$}

\htext (0 -.75){0}

\htext (-.1 -.25){$\boldsymbol{\mathcal{E}(\be) \rightarrow}$}

\htext (-.1 -.95){$\boldsymbol{\beta \rightarrow}$}

\linewd 0.02 \move(3.5 -.65) \lvec(3.5 -.35)

\htext (3.3 -.9){$\frac{\sqrt{2\log2}}{a_2}$}

\linewd 0.01 \move(4.5 -.6) \lvec(4.5 -.5)

\htext (4.44 -.9){$\frac{1}{a_1}$}

\htext (3.6 -.4){$\be a_2\sqrt{2\log2}$}

\end{texdraw}

Thus in subcase A3, the system behaves like a REM with Gaussian
distributions \cite{D1} having weight $a_2$, that is, as if
$H_N(\si)$ are i.i.d centered Gaussian with variance $a_2^2 N$. For
example, when $a_1=a_2$ then this is just the standard Gaussian REM.
It does not depend on the quantities $p_1$ and $p_2$. Even when
$p_2=0.0001$ (very small) the first level exponentials do not show
up in the limit. Further the GREM reduces to a REM. Of course, this
is so as long as $\sqrt{2\log2}<\frac{a_2}{a_1}$. This should be
contrasted with subcase A1 where the entire system behaves like sum
of two
independent REM, one Gaussian and other exponential.\\

\subsection{Gaussian - Exponential GREM}
Let us consider the situation where the driving distributions at the
first level are Gaussian, $\phi_{N,2}$ and at the second level they
are exponential, $\phi_{N,1}$. Moreover, as earlier $a_1$ and $a_2$
are the weights at the first and second level respectively. We will
use the same notation for $k(1,N), k(2,N)$ and for $p_1, p_2$. In
this case the general formula of Theorem \ref{t2.4.6} reduces to the
following:
$$\lim_N \frac{1}{N}\log Z_N(\be)= \log2-\inf\{\widetilde{f}(x,y):
x,y\geq 0;\; x\leq \sqrt{2p_1\log2};\;\frac{1}{2}x^2+y\leq
\log2\}$$ almost surely, where $$\widetilde{f}(x,y)=\frac{1}{2}x^2-\be
a_1x+y(1-\be a_2).$$

In this case to calculate the infimum we proceed as follows. Put
$g(x)=\dinf_y \widetilde{f}(x,y)$. Since $(1\be a_2)\geq 0$ for
$\be\leq\frac{1}{a_2}$, we have
\begin{align*}
g(x)&=\frac{1}{2}x^2-\be a_1 x&\text{if $\be\leq\frac{1}{a_2}$}\\
&=\frac{1}{2}x^2-\be a_1 x+(1-\be
a_2)(\log2-\frac{1}{2}x^2)&\text{if $\be>\frac{1}{a_2}$}\\
\intertext{ that is}
g(x)&=\frac{1}{2}x^2-\be a_1 x&\text{if $\be\leq\frac{1}{a_2}$}\\
&=\be a_2\left(\frac{1}{2}x^2-\frac{a_1}{a_2} x\right)+(1-\be
a_2)\log2&\text{if $\be>\frac{1}{a_2}$.}
\end{align*}
Since infimum of $f(x,y)$ is same as that of infimum of $g$ over
$x$, one has to calculate $\dinf_{0\leq x\leq \sqrt{2p_1\log 2}}
g(x)$. Here, we will have the following two scenarios.\\

\noindent\shadowbox{\bf B1: $\boldsymbol{\frac{a_1}{a_2} \leq \sqrt{2p_1\log2} }$}\vspace{1ex}

Let us assume $\frac{a_1}{a_2} \leq \sqrt{2p_1\log2}$. If
$\be\leq\frac{1}{a_2}$, then $\be a_1\leq\sqrt{2p_1\log2}$. The
function $\frac{1}{2}x^2-\be a_1 x$ decreases up to $\be a_1$ an
then increases. Hence when $\be\leq\frac{1}{a_2}$ the infimum occurs
at $x=\be a_1$. For $\be>\frac{1}{a_2}$ as $\frac{a_1}{a_2} \leq
\sqrt{2p_1\log2}$, the infimum will occur at $x=\frac{a_1}{a_2}$. So
we have the following

\begin{thm}
In the Gaussian-Exponential GREM, if $\frac{a_1}{a_2} \leq
\sqrt{2p_1\log2}$ then almost surely,
\begin{equation*}
\dlim_{N\ra\infty} \frac{1}{N} \log Z_N(\be)  =
\begin{cases}\log2+\frac{1}{2}\be^2a_1^2 & \mbox{if }\;\be\leq\frac{1}{a_2} \\
\be\left(\frac{1}{2}\frac{a_1^2}{a_2}+a_2\log2\right) & \mbox{if }\;
\be >\frac{1}{a_2}
\end{cases}
\end{equation*}
\end{thm}
As earlier, we can picture the value of $\mathcal{E}(\be)$ against
$\be$ as given below. The values of $\be$ are given under the line
and values of $\mathcal{E}(\be)$ are given above the line. The phase
transitions occur at the dark lines.

\vspace{5ex}
\begin{texdraw}
\htext (2 .1){\bf Subcase B1} \drawdim in

\linewd 0.01

\move(0 -.5) \avec(5.3 -.5)

\move(0 -0.6) \lvec(0 -0.4)

\htext (0 -.75){0}

\htext (-.1 -.2){$\boldsymbol{\mathcal{E}(\be) \rightarrow}$}

\htext (-.1 -.95){$\boldsymbol{\beta \rightarrow}$}

\htext (.5 -.4){$\frac{1}{2}\be^2a_1^2+\log2$}

\linewd 0.02 \move(2.5 -.65) \lvec(2.5 -.35)

\htext (2.44 -.9){$\frac{1}{a_2}$}

\htext (2.75 -.4){$\be(\frac{1}{2} \frac{a_1^2}{a_2}+a_2\log2)$}
\linewd 0.01 \move(3.5 -.6) \lvec(3.5 -.5)

\htext (3.3 -.9){$\frac{\sqrt{2p_1\log2}}{a_1}$}

\end{texdraw}
\vspace{3ex}

\noindent\shadowbox{\bf B2: $\boldsymbol{\sqrt{2p_1\log2} <\frac{a_1}{a_2} }$}\vspace{1ex}

Let us assume $\frac{a_1}{a_2} > \sqrt{2p_1\log2}$. If
$\be\leq\frac{1}{a_2}$ we have $\be a_1\leq\frac{a_1}{a_2}$. So the
quantity $\be a_1$ will be in $[0, \sqrt{2p_1\log2}]$ as long as
$\be\leq \frac{\sqrt{2p_1\log2}}{a_1}$. As the function
$\frac{1}{2}x^2-\be a_1 x$ decreases up to $\be a_1$ an then
increases, for $\be\leq \frac{\sqrt{2p_1\log2}}{a_1}$ the infimum
occurs at $x=\be a_1$. But for and for
$\frac{\sqrt{2p_1\log2}}{a_1}<\be\leq\frac{1}{a_2}$, the infimum
occurs at $x=\sqrt{2p_1\log2}$. For $\be>\frac{1}{a_2}$ consider the
function $\frac{1}{2}x^2-\frac{a_1}{a_2} x$ which decreases up to
$\frac{a_1}{a_2}$ and then increases. As $\frac{a_1}{a_2} >
\sqrt{2p_1\log2}$, the infimum will occur at $x=\sqrt{2p_1\log2}$.
Thus we have the following

\begin{thm}
In the Gaussian-Exponential GREM, if $\frac{a_1}{a_2} >
\sqrt{2p_1\log2}$ then almost surely,
\begin{equation*}
\dlim_{N\ra\infty} \frac{1}{N} \log Z_N(\be)  =
\begin{cases}
\log2+\frac{1}{2}\be^2a_1^2 & \mbox{if }\;\be\leq\frac{\sqrt{2p_1\log2}}{a_2} \\
p_2\log2 + \be a_1\sqrt{2p_1\log2}&\mbox{if }\;\frac{\sqrt{2p_1\log2}}{a_2}<\be\leq \frac{1}{a_2}\\
\be\left(a_1\sqrt{2p_1\log2}+a_2p_2\log2\right) & \mbox{if }\; \be
>\frac{1}{a_2}
\end{cases}
\end{equation*}
\end{thm}

As earlier, we can picture the value of $\mathcal{E}(\be)$ against
$\be$ as given below. The values of $\be$ are given under the line
and values of $\mathcal{E}(\be)$ are given above the line. The phase
transitions occur at the dark lines.

\vspace{5ex}
\begin{texdraw}
\htext (2 .1){\bf Subcase B2}

\drawdim in

\linewd 0.01

\move(0 -.5) \avec(5.4 -.5)

\move(0 -0.6) \lvec(0 -0.4)

\htext (0 -.75){0}

\htext (-.1 -.2){$\boldsymbol{\mathcal{E}(\be) \rightarrow}$}

\htext (-.1 -.95){$\boldsymbol{\beta \rightarrow}$}

\htext (.05 -.4){$\frac{1}{2}\be^2a_1^2+\log2$}

\linewd 0.02 \move(3 -.65) \lvec(3 -.3)

\htext (2.95 -.9){$\frac{1}{a_2}$}

\htext (3.2 -.4){$\be(a_1\sqrt{2p_1\log2}+a_2p_2\log2)$}

\move(1.05 -.65) \lvec(1.05 -.3)

\htext (.8 -.9){$\frac{\sqrt{2p_1\log2}}{a_1}$}

\htext (1.26 -.4){$\be a_1\sqrt{2p_1\log2}+p_2\log2$}

\linewd 0.01 \move(4.5 -.6) \lvec(4.5 -.5)

\htext (4.3 -.9){$\frac{\sqrt{2\log2}}{a_1}$}
\end{texdraw}

Remarks similar to Exponential-Gaussian GREM apply here as well.
Subcase B1 is similar to subcase A2. Here also the term
$\frac{1}{2}\be \frac{a_1^2}{a_2}$ is not reminiscent of anything we
know.

Subcase B2 is similar to that of subcase A1. That is in subcase B2,
the limiting free energy is sum of two REM free energies -- one is
of Gaussian REM and other is of exponential REM. To be precise, the
Gaussian REM  limiting free energy (keeping in mind that for $N$
fixed, the $k(1,N)$ particle system has $\mathcal{N}(0,a^2N)$
Hamiltonian as opposed to $\mathcal{N}(0,a^2k(1,N))$ yields, a.s,
\begin{equation}
\widetilde{\mathcal{E}}_1(\be)=
\begin{cases}
p_1\log2 +\frac{1}{2} a_1^2 \be^2 &\text{ if $\be\leq \frac{\sqrt{2p_1\log2}}{a_1}$}\\
\be a_1\sqrt{2p_1\log2} &\text{ if $ \be >
\frac{\sqrt{2p_1\log2}}{a_1}$.}
\end{cases}
\end{equation}
On the other hand, for fixed $N$, if we have configurations
$2^{k(1,N)}$  then the exponential REM limiting free energy, with
Hamiltonian as $a_2$ times double exponential random variable
yields, a.s.,
\begin{equation}
\widetilde{\mathcal{E}}_2(\be)=
\begin{cases}p_2\log2 &\text{ if $\be\leq \frac{1}{a_2}$}\\
\be p_2 a_2\log2 &\text{ if $\be > \frac{1}{a_2}$.}
\end{cases}
\end{equation}
Now it is easy verify that, in subcase B2, a.s.
$$\mathcal{E}(\be)=\widetilde{\mathcal{E}}_1(\be)+\widetilde{\mathcal{E}}_2(\be).$$
The reader should note that to compare subcase B2 with subcase A1,
we interchange $a_2$ with $a_1$ and $p_2$ with $p_1$ (to maintain
the same weights and proportions for the exponential and Gaussian
levels).

The last interesting note is that in Gaussian-Exponential GREM, the system never
reduces completely to a Gaussian REM as happened in subcase A3.

Thus the large deviation technique allows the use of different
distributions at different levels leading to some interesting
phenomenons. The conclusions of Exponential-Gaussian GREM differ
from those of Exponential-Gaussian. The system may reduce to a
Gaussian REM even with a very small weight is associated to that
level. Even the system may appear as a system of two independent
REMs separated by a big wall preventing them to interact between
each other. Moreover, there are situations where we could not
explain the terms present in the expression for energy.

\cleardoublepage

\chapter{More Tree Structures including Randomness}


In this chapter, we will consider several models similar to that of
Generalized Random Energy Model. In the previous chapter, we
formulated GREM in general tree set up with out giving any examples
of general tree structures. The set up also allows us to randomize
the tree structure. First we consider regular trees but the trees
are random, driven by Poisson random variables. Then we consider
non-regular random trees again driven by Poisson random variables.
We prove that in both the cases the free energy exists for almost
every tree sequences and they are same as that of usual
deterministic tree GREMs for almost every sample point. Also we
consider Multinomial trees. These will be explained later.

The usual GREM has hierarchical structure, and it is so in all the
above mentioned models. In 2006, Bolthausen and Kistler\cite{BK}
defined a model which is a generalization of the GREM where the
model is no longer hierarchical. They called the model as
non-hierarchical version of GREM and prove the existence of the free
energy by using second moment method. Surprisingly, the energy
expression is again the same as that of the usual GREM. So the
non-hierarchy does not play a role in the limiting free energy. We
produce an alternative proof of their result through large deviation
techniques and show that the free energy of this model is minimum of
certain hidden GREMs. Then we introduce another model, block tree
GREM where the energy is maximum over certain GREM energies. We
present further generalization in a model, in the next chapter,
through which we can get all the models REM, GREM,
Bolthausen-Kistler model and their versions with the external field.

\smallskip

\goodbreak
\section{Regular Poisson GREM}
In generalized random energy model, we have randomness coming from
the driving distributions. The reformulation of GREM in general tree
structure allows us to introduce another randomness at the tree
level which is independent of the randomness of the Hamiltonians. As
usual for $N$ particles system, let $\{k(i,N),\,1\leq i\leq n\}$ be
a partition of $N$ into $n$ (the level of the tree) positive
integers. Consider, for each $N$, independent random variables
$L_{1N}, \cdots,L_{nN}$ where $L_{iN} \sim P(2^{k(i,N)})$, i.e. a
Poisson random variable with parameter $2^{k(i,N)}$ for $1\leq i\leq
n$. Let us construct a random tree with $(1+L_{iN})$ nodes at the
i-th level below each node of the (i-1)-th level. That is, at the
first level there will be $1+L_{1N}$ many edges and at the second
level there will be total $(1+L_{1N})(1+L_{2N})$ many edges. Here we
are considering $1+L_{iN}$ instead of $L_{iN}$ itself, to take care
of the situation $L_{iN}=0$ so that each branch in the tree is of
length $n$. Once again we denote the edges at the first level by
$\si_1$ and the second level edges below $\si_1$ as $\si_1\si_2$ and
so on. The weight of the $i$-th level is $a_i>0$. Similarly we will
associate independent random variable $\xi(\si_1\cdots\si_i)$ with
the edge $\si_1\cdots\si_i$. In this case for $N$ particle system,
instead of $2^N$ configurations we will have
$(1+L_{1N})(1+L_{2N})\cdots (1+L_{nN})$ many configurations. Of
course this is also a regular tree, but random, and could be called
regular Poisson tree. The corresponding GREM model, where the
Hamiltonian for the configuration $\si=(\si_1,\cdots,\si_n)$ is
defined as $\dsum_{i=1}^n a_i\xi(\si_1\cdots\si_i)$, can be called a
{\em regular Poisson tree GREM} with parameter
$\tilde{k}=(k(1,N),\cdots,k(n,N))$. The next result says that if the
same conditions as in Corollary \ref{c2.3.4} hold then even with
randomization of tree, the conclusion holds for almost every tree
sequence.
\begin{prop}\label{p3.1.1}
Consider a regular Poisson tree GREM with parameter $\tilde{k}$. The
following is true:

a) If $\sum\limits_{N\geq 1} 2^{k(1,N)+\cdots+k(i,N)} q_{1N}\cdots
q_{iN} < \infty$, for some $i, 1\leq i \leq n$ then for a.e. tree
sequence, a.s. eventually, $\mu_N(\tri) = 0.$

b) If $\sum\limits_{N\geq1}
2^{-(k(1,N)+\cdots+k(i,N))}q_{1N}^{-1}\cdots q_{iN}^{-1} < \infty$,
for each $i= 1, \cdots, n$, then for a.e. tree sequence the
following is true: for any $\ep>0$, a.s. eventually,
\[(1-\ep)q_{1N}\cdots q_{nN} \leq \mu_N(\tri) \leq (1+\ep)q_{1N}\cdots
q_{nN}.\]
\end{prop}
\begin{proof}
a) It suffices to verify the hypothesis of Theorem
\ref{t2.3.3}(a) holds for almost every tree sequence, that is, $\sum\limits_{N\geq n} B_{iN}
q_{1N}\cdots q_{iN}<\infty$ for some $i$. Recall that $B_{iN}$ is the number of branches at the
$i$-th level.

But we could
prove a stronger statement, namely, if for some $i$ with $1\leq
i\leq n$, $\dsum_{N\geq n} 2^{k(1,N)+\cdots+k(i,N)} q_{1N}\cdots
q_{iN} < \infty$, then $\mbf{E}_T\sum\limits_{N\geq n} B_{iN}
q_{1N}\cdots q_{iN}<\infty$ for that $i$ where $\mbf{E}_T$ is the
tree expectation. Since the tree randomness is independent of the
Hamiltonian randomness, in view of the hypothesis, it suffices to
show
\begin{equation}\label{e4.1.0}
\mbf{E}_T B_{iN}\leq 2^{k(1,N)+\cdots+k(i,N)+i}.
\end{equation}
Using independence of the random variables $(L_{jN}, \; 1\leq j \leq
i)$, we get
\begin{equation}\label{e4.1.00}
\mbf{E}_T B_{iN} = \mbf{E}\prod\limits_{j=1}^i
(1+L_{jN})=\prod\limits_{j=1}^i \mbf{E}
(1+L_{jN})=\prod\limits_{j=1}^i (1+2^{k(j,N)})\leq
2^i\prod\limits_{j=1}^i 2^{k(j,N)}.
\end{equation}

b) \quad It is enough to show that for fixed $\ep>0$, almost every
tree sequence satisfies the stated conclusion. This is achieved by
verifying that the hypothesis of Theorem \ref{t2.3.3}(b) holds for
almost every tree sequence, that is, $\sum\limits_{N\geq n} \frac{s_{iN}^2}{B_N^2 q_{1N}\cdots
q_{iN}}<\infty$. Recall that,
$s_{iN}^2=\sum\limits_{\si_1,\cdots,\si_i}e^2(\si_1,\cdots,\si_i)$
where  $e(\si_1\si_2\cdots\si_i)$ denotes the number of nodes at the
$n$-th level below the node $\si_1\si_2\cdots\si_i$ and $B_N$ is the
number of leaves or the the total number of branches in the tree.

Here also we prove a stronger statement, namely,
$\mbf{E}_T\sum\limits_{N\geq n} \frac{s_{iN}^2}{B_N^2 q_{1N}\cdots
q_{iN}}<\infty$ for each $i$ where $\mbf{E}_T$ is the tree
expectation. Again, since the tree randomness is independent of the
Hamiltonian randomness, in view of the hypothesis, it suffices to
show
\begin{equation}\label{e4.1.1}
\mbf{E}_T \left(\frac{s_{iN}^2}{B_N^2}\right)\leq
2^{-(k(1,N)+\cdots+k(i,N))}.
\end{equation}
But due to regularity of the tree $s_{iN}^2 = \prod\limits_{j=1}^i
(1+L_{jN})\prod\limits_{j=i+1}^n (1+L_{jN})^2$ and $B_N^2 =
\prod\limits_{j=1}^n (1+L_{jN})^2$. Hence $$\frac{s_{iN}^2}{B_N^2} = \prod\limits_{j=1}^i
\frac{1}{1+L_{jN}}.$$
Thus using the independence of the random variables $(L_{jN}, \; 1\leq j \leq
n)$, we get
\begin{equation}\label{e4.1.2}
\mbf{E}_T \left(\frac{s_{iN}^2}{B_N^2}\right)= \prod\limits_{j=1}^i
\mbf{E} \left(\frac{1}{1+L_{jN}}\right).
\end{equation}

Since for a Poisson random variables $X$ with parameter $\la$, $\mbf{E}\frac{1}{1+X}=\frac{1}{\la}\left(1-e^{-\la}\right)$ and since $L_{jN}\sim P(2^{k(i,N)})$ we have,
\begin{equation}\label{e4.1.3}
\mbf{E}
\left(\frac{1}{1+L_{jN}}\right)=2^{-k(j,N)}\left(1-e^{-2^{k(j,N)}}\right)\leq
2^{-k(j,N)}.
\end{equation}

Substituting (\ref{e4.1.3}) in (\ref{e4.1.2}) we get (\ref{e4.1.1}).
\end{proof}
Now further if we assume that $\frac{k(i,N)}{N} \ra p_i \,(>0)$ for $1\leq i\leq n$ and the random variables $\xi(\si_1\cdots\si_i)$ are distributed like $\phi_{N,\ga_i}$ as defined in (\ref{e2.4.21}), the rest of the proof for existence of the free energy is the same as that of Theorem \ref{t2.4.6}. Thus if the sequence $\left\{\frac{\xi(\si_1\cdots\si_i)}{N}\right\}$ satisfies LDP with good rate function $\mathcal{I}_i$ for each $i$ with scale parameter $N$, then we have the following.
\begin{thm}\label{t3.1.2}
Assume the setup as in the above paragraph. For regular Poisson tree GREM, for almost every tree sequences, almost
surely,
\[\dsty\lim_{N} \frac{1}{N} \log
Z_N(\be)=\log2-\dinf_{\widetilde{x}\in\Psi}\left\{\dsum_{i=1}^n\left(\mathcal{I}_i(x_i)-\be
a_i x_i\right)\right\},\]
where
\[\Psi = \left\{\widetilde{x}\in \mathbb{R}^n :
\dsum_{i=1}^k \mathcal{I}_i(x_i) \leq \dsum_{i=1}^k p_i
\log2, \hspace{1ex} 1\leq k\leq n\right\}.\]
\end{thm}

Thus, though we have another randomness in
the setup of the model, the limiting free energy remains the same.
That is why, in the original setup of GREM by Derrida, though ($\al_i^N$) may not be an
integer, with out any loss of generality one can consider the number
of branches at the $i$-th level to be $[\al_i^N]$. We make this more precise in Remark
\ref{r3.3.1}.

We also note that, in this model the number of configurations in the
configuration space may not be of the form $\al^N$, where $\al$ ia a
natural number. Instead, it is of the form $l_1l_2\cdots l_n$.

In the above model, there are $n$ random variables controlling the
number of nodes at the $n$ level of the GREM. Since $n$ is fixed and
$N\ra\infty$, one may get the impression that this extra randomness
is not showing up in the final conclusion, namely, in the expression
for the free energy. The next model shows that such an impression is
not correct.

\goodbreak

\section{Poisson GREM}
In the above model we randomized the tree sequences so that the
underlying trees once again remain regular. But the general
formulation allow us to consider a non-regular tree sequence. In the
previous model, for each $N$, we randomized the tree using $n$ many
Poisson random variables corresponding to the levels of the tree.
But it is conceivable to use independent Poisson variables at each
of the nodes to construct the tree as well as the configuration
space. This is what we do now. As in the previous model, let
$\{k(1,N),\cdots, k(n,N)\}$ be a partition of $N$.  Unlike in that
model, now let us consider an $n$-level tree with $P(2^{k(i,N)})+1$,
many nodes below each of the nodes at the $(i-1)$-th level for
$1\leq i\leq n$. Here $P(2^{k(i,N)})$ denotes a Poisson random
variable with parameter $2^{k(i,N)}$. In other words, instead of
fixing one random variable and taking so many nodes below each of
the $(i-1)$-th level nodes, we now fix one random variable for each
node of the $(i-1)$-th level and take so many nodes below that. Let
us assume all these Poisson random variables are  independent. As in
the previous model, we denote a typical edge at the $i$-th level by
$\si_1\cdots\si_i$ below the edge $\si_1\cdots\si_{i-1}$ and we
associate independent random variables $\xi(\si_1\cdots\si_i)$ to
it.  We assume that this family $\{\xi(\si_1\cdots\si_i)\}$ is
independent of the above Poisson family. For $1\leq i\leq n$, we
have a positive number $a_i$ denoting the  weights for the $i$-th
level of the tree. Now we define the Hamiltonian for the
configuration $\si=(\si_1,\cdots,\si_n)$ as
\begin{equation*}
H_N(\si)=\sum_{i=1}^n a_i\xi(\si_1\cdots\si_i).
\end{equation*}

This model can be called a true {\em Poisson tree GREM} with
parameter $\tilde{k} = (k(1,N),\cdots, k(n,N))$. Here we randomize,
rather Poissonize, the tree in its full form. Even in this case also
the model behaves the same way and we get the same conclusions as
that of the above model, that is, (a) and (b) of Proposition
\ref{p3.1.1} remain true. This is the content of the next
proposition.

\begin{prop}\label{p3.2.1}
Consider a  Poisson tree GREM with parameter $\tilde{k}$. The
following is true:

a) If $\sum\limits_{N\geq 1} 2^{k(1,N)+\cdots+k(i,N)} q_{1N}\cdots
q_{iN} < \infty$, for some $i, 1\leq i \leq n$ then for a.e. tree
sequence, a.s. eventually, $\mu_N(\tri) = 0.$

b) If $\sum\limits_{N\geq1}
2^{-(k(1,N)+\cdots+k(i,N))}q_{1N}^{-1}\cdots q_{iN}^{-1} < \infty$,
for each $i= 1, \cdots, n$, then for a.e. tree sequence the
following is true: for any $\ep>0$, a.s. eventually,
\[(1-\ep)q_{1N}\cdots q_{nN} \leq \mu_N(\tri) \leq (1+\ep)q_{1N}\cdots
q_{nN}.\]
\end{prop}
We need the following two inequalities to prove the proposition.

{\em Let $a\geq 1$,  $b\geq 1$  and $\lambda >0$. Suppose that $X
\sim P(a\lambda)$ and $Y \sim P(b\lambda)$ are independent random
variables. Then

\begin{equation}\label{e3.2.6}\mbf{E}\left(\frac{X+a}{X+Y+a+b}\right)^2 \leq
2\left(\frac{a}{a+b}\right)^2, \end{equation}

and

\begin{equation}\label{e3.2.7}\mbf{E}
\frac{X+a}{(X+Y+a+b)^2} \leq
\frac{a}{(a+b)^2}\,\frac{1}{\lambda}.\end{equation}}

Both these rely on conditioning. Since $X$ and $Y$ are independent
Poisson random variables, given $X+Y=l$, the conditional
distribution of $X$ is binomial $\left(l, \frac{a}{a+b}\right)$. So
\begin{align*}
&\mbf{E}\left(\frac{X+a}{X+Y+a+b}\right)^2\\
 =& \mbf{E}\left\{\frac{1}{(X+Y+a+b)^2}\mbf{E}\left[(X+a)^2\;\boldsymbol{|}\;  X+Y\right]\right\}\\
 =& \mbf{E}\left\{\frac{1}{(X+Y+a+b)^2}\frac{a^2(X+Y+a+b)^2+(X+Y)ab}{(a+b)^2}\right\}\\
 \leq& \frac{a^2}{(a+b)^2} \quad\quad \mbox{since } a\geq 1,
\end{align*} and
\begin{align*}
&\mbf{E}\left(\frac{X+a}{(X+Y+a+b)^2}\right)\\
 =& \mbf{E}\left\{\frac{1}{(X+Y+a+b)^2}\mbf{E}[(X+a)\;\boldsymbol{|}\; X+Y]\right\}\\
 =& \mbf{E}\left\{\frac{1}{(X+Y+a+b)^2}\frac{a(X+Y+a+b)}{(a+b)}\right\}\\
 =& \frac{a}{(a+b)}\mbf{E}\frac{1}{(X+Y+a+b)}\\
\intertext{}
 <&\frac{a}{(a+b)}\mbf{E}\frac{1}{(X+Y+1)}\quad\quad \mbox{since } a+b> 1\\
 <&\frac{a}{\la(a+b)^2}.
\end{align*}

\begin{proof}[Proof of Proposition~\ref{p3.2.1}]
The proof is routine but involves rather cumbersome notation. To
describe the random tree for the $N$-particle system, let $L_0\sim
P(2^{k(1,N)})$, the number of edges at the first level. For $1\leq
\si_1 \leq L_0+1$, let $L_{\si_1}\sim P(2^{k(2,N)})$, the number of
edges at the second level below the first level edge $\si_1$. In
general, for $\si_1\si_2\cdots\si_i$, with $1\leq \si_1 \leq L_0+1,
\, 1\leq \si_2 \leq L_{\si_1}+1, \, \cdots, 1\leq \si_i \leq
L_{\si_1\cdots\si_{i-1}}+1$, let $L_{\si_1\cdots\si_i}\sim
P(2^{k(i+1,N)})$, the number of edges at the $(i+1)$-th level
below the edge $\si_1\cdots\si_{i-1}$ at the $i$-th level.

To prove part (a), it suffices to
show, as in Proposition \ref{p3.1.1},  that
\begin{equation*}
\mbf{E}_T B_{iN}\leq 2^{k(1,N)+\cdots+k(i,N)+i}.
\end{equation*}
Since, in this model $$B_{iN}=\sum_{\si_1}\cdots
\sum_{\si_{i-1}}(L_{\si_1\cdots\si_{i-1}}+1),$$ the proof is immediate.

To prove (b), as in Proposition \ref{p3.1.1}, it suffices to show
that for each $i$,
\[\mbf{E}_T \left(\frac{s_{iN}^2}{B_N^2}\right)\leq 2^n
2^{-(k(1,N)+\cdots+k(i,N))}.\]  But in this model,
$$\frac{s_{iN}^2}{B_N^2} = \sum_{\si_1}\cdots\sum_{\si_i}
\left(\frac{\sum_{\si_{i+1}}\cdots\sum_{\si_{n-1}}
(L_{\si_1\cdots\si_{n-1}}+1)}{\sum_{\si_1}\cdots
\sum_{\si_{n-1}}(L_{\si_1\cdots\si_{n-1}}+1)}\right)^2.$$ To
calculate the expectation we proceed as follows. Let $\mathcal{F}_o$
be the $\si$-field generated by $L_0$, $\mathcal{F}_1$ be the
$\si$-field generated by $\{L_0,L_{\si_1}:\, 1\leq \si_1\leq
L_0+1\}$ and in general $\mathcal{F}_i$ be the $\si$-field generated
by $\{L_0, L_{\si_1},\cdots,L_{\si_1\cdots\si_i}:\, 1\leq \si_1\leq
L_0+1,1\leq \si_2\leq L_{\si_1}+1,\cdots,1\leq\si_i\leq
L_{\si_1,\cdots,\si_{i-1}}\}$ for $i=0,\cdots,n-2$. Let $\mbf{E}_i$
be the conditional expectation given $\mathcal{F}_i$. Then
(\ref{e3.2.6}) suggests that,
\[\begin{array}{rl}
\mbf{E}_{n-2} \left(\dfrac{s_{iN}^2}{B_N^2}\right)&=\mbf{E}_{n-2} \dsum_{\si_1}\cdots\dsum_{\si_i}
\left(\dfrac{\sum_{\si_{i+1}}\cdots\sum_{\si_{n-1}}
L_{\si_1\cdots\si_{n-1}}+\sum_{\si_{i+1}}\cdots\sum_{\si_{n-1}}1}{\sum_{\si_1}\cdots
\sum_{\si_{n-1}}L_{\si_1\cdots\si_{n-1}}+\sum_{\si_1}\cdots
\sum_{\si_{n-1}}1}\right)^2 \vspace{2ex}\\
& \leq 2 \dsum_{\si_1}\cdots\dsum_{\si_i}
\left(\dfrac{\sum_{\si_{i+1}}\cdots\sum_{\si_{n-2}}
(L_{\si_1\cdots\si_{n-2}}+1)}{\sum_{\si_1}\cdots
\sum_{\si_{n-2}}(L_{\si_1\cdots\si_{n-2}}+1)}\right)^2.
\end{array}\]
Similarly, $$\mbf{E}_{n-3}\mbf{E}_{n-2} \left(\frac{s_{iN}^2}{B_N^2}\right)\leq 2^2 \sum_{\si_1}\cdots\sum_{\si_i}
\left(\frac{\sum_{\si_{i+1}}\cdots\sum_{\si_{n-3}}
(L_{\si_1\cdots\si_{n-3}}+1)}{\sum_{\si_1}\cdots
\sum_{\si_{n-3}}(L_{\si_1\cdots\si_{n-3}}+1)}\right)^2,$$ and thus
$$\mbf{E}_{i}\cdots\mbf{E}_{n-2} \left(\frac{s_{iN}^2}{B_N^2}\right)\leq 2^{n-i-1} \sum_{\si_1}\cdots\sum_{\si_i}
\left(\frac{1}{\sum_{\si_1}\cdots
\sum_{\si_{i}}(L_{\si_1\cdots\si_{i}}+1)}\right)^2.$$ Now we can use
(\ref{e3.2.7}) to calculate further conditional expectations so that
$$\mbf{E}_{i-1}\mbf{E}_{i}\cdots\mbf{E}_{n-2} \left(\frac{s_{iN}^2}{B_N^2}\right)
\leq 2^{n-i-1}\frac{1}{2^{k(i,N)}} \sum_{\si_1}\cdots\sum_{\si_{i-1}}
\left(\frac{1}{\sum_{\si_1}\cdots
\sum_{\si_{i}}(L_{\si_1\cdots\si_{i-1}}+1)}\right)^2,$$ and so on to get
$$\mbf{E}_{0}\mbf{E}_{1}\cdots\mbf{E}_{n-2} \left(\frac{s_{iN}^2}{B_N^2}\right)
\leq 2^{n-i-1}\frac{1}{2^{k(i,N)+\cdots+k(1,N)}}.$$
Since
$\mbf{E}_T\left(\frac{s_{iN}^2}{B_N^2}\right)=\mbf{E}_{0}\mbf{E}_{1}\cdots\mbf{E}_{n-2}
\left(\frac{s_{iN}^2}{B_N^2}\right)$, the proof is complete.
\end{proof}

Once again to verify the existence of limiting free energy, one has
to verify all the steps involved in section \ref{s2.4}. To be precise,
if we assume that $\frac{k(i,N)}{N} \ra p_i \,(>0)$ for $1\leq i\leq
n$, the sequence $\left\{\frac{\xi(\si_1\cdots\si_i)}{N}\right\}$
satisfies LDP with good rate function $\mathcal{I}_i$ for each $i$, then we get
\begin{thm}
For almost every tree sequences, almost surely,
\[\dsty\lim_{N} \frac{1}{N} \log
Z_N(\be)=\log2-\dinf_{\widetilde{x}\in\Psi}\left\{\dsum_{i=1}^n\left(\mathcal{I}_i(x_i)-\be
a_i x_i\right)\right\},\]
where
\[\Psi = \left\{\widetilde{x}\in \mathbb{R}^n :
\dsum_{i=1}^k \mathcal{I}_i(x_i) \leq \dsum_{i=1}^k p_i
\log2, \hspace{1ex} 1\leq k\leq n\right\}.\]
\end{thm}

\begin{rem}
Recall that a tree is regular if for any $i$, the number of nodes
below a $(i-1)$-th level node depends only on $i$ and not on the
specific nodes. We say that a tree sequence is regular if after some
stage each tree in the sequence is regular. Under suitable
conditions $--$ for instance, when $\sum\limits_N
e^{-k(i,N)}<\infty$  for some $i$ $--$ it is possible to show  that
almost every tree sequence ceases to be regular. Though in this case
the probability that the tree sequence will consist of regular trees
is very small, we did not get any further new result except that,
now the limiting free energy is constant for almost every tree
sequences as well as almost every sample points.
\end{rem}

\section{Multinomial tree GREM}

In the above two models, we randomized the number of nodes at each
level keeping the average fixed. It is also possible to randomize
the vector $\tilde{k}$ suitably. To do that, we fix $p_i>0$ for
$1\leq i \leq n$ with $\sum\limits_1^n p_i= 1$. Now consider an
$n$-faced die with $p_i$ being the chance of face $i$ appearing in a
throw. Now we can consider two experiments with this die. {\em
Firstly}, we can consider an indefinite throws of the die and for
$N$ particle system let $K(i,N)$ be the number of times face $i$
appears in first $N$ throws. In the {\em second}, experiment for $N$
particle system we will throw the die independently $N$ times and
observe the outcomes. With the same notation, let $K(i,N)$ be the
number of times face $i$ appears. Clearly, in both the cases
$K(i,N)\geq 0$ and $\sum\limits_{i=1}^n K(i,N) =N$. We can consider
GREM with parameter $\tilde{K}$, that is where the under lying tree
has $2^{K(i,N)}$ many edges below each of the $(i-1)$ level node.
These can be called {\em multinomial tree GREM} of first kind and
{\em multinomial tree GREM} of second kind with parameter
$\tilde{p}=(p_1,\cdots,p_n)$ respectively. With the same notation as
in section \ref{s2.3}, in this case we have
\begin{equation*}
\mbf{E}_T B_{iN} = \mbf{E}\prod\limits_{j=1}^i
2^{K(j,N)}=\dsum_{k=0}^N 2^k \left(\begin{array}{c}N\\k
\end{array}\right)\left(\dsum_1^i p_j\right)^k \left(1-\dsum_1^i p_j\right)^{N-k}\hspace{-2ex}
=\left(1+\dsum_1^i p_j\right)^N\hspace{-2ex}.
\end{equation*}
But for any $x$, as $(1+x)^N\leq e^{Nx}$, we have
\begin{equation}\label{e3.3.8}
\mbf{E}_T B_{iN}\leq e^{N\dsum_1^i
p_j}=2^{\frac{N}{\log2}(p_1+\cdots+p_i)}.
\end{equation}
On the other hand, in this case,
$s_{iN}^2=2^{K(1,N)+\cdots+K(i,N)}2^{2(K((i+1),N)+\cdots+K(n,N))}$
and $B_N^2=2^{2\dsum_{j=1}^n K(j,N)}$ so that
$\left(\frac{s_{iN}^2}{B_N^2}\right)=2^{-\dsum_1^i K(j,N)}$. Once
again using the fact that $\dsum_1^iK(j,N)$ is binomial with
parameters $N$ and $\dsum_1^i p_j$, we see that
$$\mbf{E}2^{-\dsum_1^i K(j,N)}=\left(1-\frac{1}{2}\dsum_1^i
p_j\right)^N .$$ Hence by same inequality as earlier, we have
\begin{equation}\label{e3.3.9}
\mbf{E}2^{-\dsum_1^i K(j,N)}\leq e^{-\frac{N}{2}\dsum_1^i
p_j}=2^{-\frac{N}{2\log2}\dsum_1^i p_j}.
\end{equation}

Combining the above observations (\ref{e3.3.8}) and (\ref{e3.3.9}),
we have the following
\begin{cor}
Consider a multinomial tree GREM either of first kind or of second
kind with parameter $\tilde{p}$. Let
$\tri=\tri_1\times\cdots\times\tri_n$ be a box in $\mathbb{R}^n$ and
$q_{iN}=P(\frac{\xi(\si_1\cdots\si_i}{N})\in\tri_i$.

a) If $\sum\limits_{N\geq1}
2^{\frac{N}{\log2}(p_1+\cdots+p_i)}q_{1N}\cdots q_{iN} < \infty$,
for some $i, 1\leq i \leq n$ then for a.e. tree sequence, a.s.
eventually, $\mu_N(\tri) = 0.$

b) If $\sum\limits_{N\geq1}
2^{-\frac{N}{2\log2}(p_1+\cdots+p_i)}q_{1N}^{-1}\cdots q_{iN}^{-1} <
\infty$, for each $i= 1, \cdots, n$, then for a.e. tree sequence the
following is true: for any $\ep>0$, a.s. eventually,
\[(1-\ep)q_{1N}\cdots q_{nN} \leq \mu_N(\tri) \leq (1+\ep)q_{1N}\cdots
q_{nN}.\]
\end{cor}

Notice the difference in the hypothesis of (a) and (b) in the above
corollary. More specifically, there is a factor $\frac{1}{2}$ extra
in the exponent of 2 in part (b). This difference will not help us
to obtain the exact support of the empirical measure $\mu_N$. For
the variational problem this support is very essential. But, of
course, we have strong law of large number in our hand. By SLLN,
almost surely, $\frac{1}{N}(K(1,N),\cdots,K(n,N)) \ra
(p_1,\cdots,p_n)$ in both the models, first or second kind. In the
first case, we use SLLN for a sequence of i. i. d. random variables
and in the second case we use SLLN for array of rowwise independent
random variables\cite{HMT89}. According to Hu et al in \cite{HMT89}:
{\em If $\{X_{nk}\}$ be an array of rowwise independent random
variables such that $E X_{nk}=0$ and there exists a random variable
$X$ with $E X^2<\infty$ so that for all $n$ and $k$ and for all
$t>0$, $P(|X_{nk}|>t)\leq P(|X|>t)$, then $\frac{1}{n}\dsum_{k=1}^n
X_{nk} \ra 0$ almost surely.} Since we can write $K(i,N)$ as
$\dsum_{k=1}^N X_{Nk}^{(i)}$ where $X_{Nk}^{(i)}$ are Bernoulli with
success probability $p_i$, the above result is applicable for the
second kind model.  Thus, in either of the cases, for every $\ep>0$
almost every tree sequences after some stage
$$N(p_1+\cdots+p_i -\ep)<K(1,N)+\cdots+K(i,N)<N(p_1+\cdots+p_i
+\ep)$$ for $i=1, \cdots, n$. That is for almost every tree
sequences and for any arbitrary $\ep>0$, $B_{iN}\leq
2^{N(p_1+\cdots+p_i +\ep)}$ and
$\frac{s_{iN}^2}{B_N^2}<2^{-N(p_1+\cdots+p_i -\ep)}$ for $i=1,
\cdots, n$. As a consequence, we can restate Theorem \ref{t2.3.3} as
follows

\begin{cor}
Consider a multinomial tree GREM with parameter $\tilde{p}$. Let
$\tri=\tri_1\times\cdots\times\tri_n$ be a box in $\mathbb{R}^n$ and
$q_{iN}=P(\frac{\xi(\si_1\cdots\si_i}{N})\in\tri_i$. Let $\ep>0$.

a) If $\sum\limits_{N\geq1}2^{N\ep}
2^{N(p_1+\cdots+p_i)}q_{1N}\cdots q_{iN} < \infty$, for some $i,
1\leq i \leq n$ then for a.e. tree sequence, a.s. eventually,
$\mu_N(\tri) = 0.$

b) If $\sum\limits_{N\geq1}2^{N\ep}
2^{-N(p_1+\cdots+p_i)}q_{1N}^{-1}\cdots q_{iN}^{-1} < \infty$, for
each $i= 1, \cdots, n$, then for a.e. tree sequence the following is
true: a.s. eventually,
\[(1-\ep)q_{1N}\cdots q_{nN} \leq \mu_N(\tri) \leq (1+\ep)q_{1N}\cdots
q_{nN}.\]
\end{cor}

Now the proof of existence of the asymptotic free energy for  this
model is routine and for almost every tree sequence the expression
for free energy will be same as that of the deterministic tree model
where $\frac{k(i,N)}{N}\ra p_i$ for $1\leq i\leq n$. To state
precisely, let us we assume that the sequence
$\left\{\frac{\xi(\si_1\cdots\si_i)}{N}\right\}$ satisfies LDP with
good rate function $\mathcal{I}_i$ for each $i$, then we have the
following.
\begin{thm}
With the setup as in the above paragraph, for almost every tree
sequence, almost surely,
\[\dsty\lim_{N} \frac{1}{N} \log
Z_N(\be)=\log2-\dinf_{\widetilde{x}\in\Psi}\left\{\dsum_{i=1}^n\left(\mathcal{I}_i(x_i)-\be
a_i x_i\right)\right\},\]
where
\[\Psi = \left\{\widetilde{x}\in \mathbb{R}^n :
\dsum_{i=1}^k \mathcal{I}_i(x_i) \leq \dsum_{i=1}^k p_i
\log2, \hspace{1ex} 1\leq k\leq n\right\}.\]
\end{thm}

\begin{rem}\label{r3.3.1} Going back to Theorem \ref{t2.3.3}, let $(T_N)$ and
$(\tilde{T}_N)$ be two sequences of trees. Suppose there are numbers
$C>c>0$ such that for each $i$, $c\leq
\frac{\tilde{s}_{iN}}{s_{iN}}\leq C$ and $c\leq
\frac{\tilde{B}_{iN}}{B_{iN}}\leq C$. Then it is easy to see that,
hypothesis of Theorem \ref{t2.3.3}(b) holds for $(T_N)$ iff it holds
for $(\tilde{T}_N)$. Accordingly, the conclusion of Theorem
\ref{t2.3.3}(b) holds for $(T_N)$ iff it holds for $(\tilde{T}_N)$.
Same remark applies for Theorem \ref{t2.3.3}(a).
\end{rem}

\section{Bolthausen - Kistler GREM}
In 2006, Bolthausen and Kistler proposed a model where they tried to
go beyond the natural {\em ultrametricity} of the GREM model. To
recall, a metric $d$ is {\em ultrametric} if in the metric property
one replaces the triangle inequality by $d(x,z)\leq
\max(d(x,y),d(y,z))$. In all of the above GREM model one can define
a metric on the configuration space $\Si_N$ of the $N$ particle
system through the covariance structure of the Hamiltonian. To be
precise for two configurations $\si$ and $\tau$ in $\Si_N$,
\begin{equation}\label{e3.4.9-a}
d(\si,\tau)= \sqrt{E(H_N(\si)-H_N(\tau))^2}.
\end{equation}
Here as usual, $H_N(\si)=\dsum_{i=1}^n a_i \xi(\si_1\cdots \si_i)$
with the usual GREM notation. Let $\si=(\si_1,\cdots,\si_n)$ and
$\tau=(\tau_1,\cdots,\tau_n)$. If $\si=\tau$ then $d(\si,\tau)= 0$.
If $\si_i=\tau_i$ for $1\leq i\leq k< n$ but $\si_{(k+1)}\neq
\tau_{(k+1)}$ then $d(\si,\tau)=\sqrt{2\dsum_{i=k+1}^n
a_i^2E\xi(\si_1\cdots\si_i)^2}$, assuming that the $\xi$'s are
symmetric with finite variance. The distance between any two
configuration will be maximum, when they differ at the first level
of the tree. The longer the initial segment of $\si$ and $\tau$
coincide, the closer they are. Also it is quite easy to verify that
this metric indeed is an ultrametric. To see this, it is enough to
see the level of difference among the configurations. Suppose, we
have any three configurations $\si,\tau$ and $\eta$ in $\Si_N$. Let
$k_i, k_2$ and $k_3$ be the maximum non-negative integers so that
$\si_i=\tau_i$ for $1\leq i\leq k_1$; $\si_i=\eta_i$ for $1\leq
i\leq k_2$ and $\tau_i=\eta_i$ for $1\leq i\leq k_3$ respectively.
To show $d(\si,\tau)\leq \max(d(\si,\eta),d(\eta,\tau))$, we only
need to show that $k_1\geq\min(k_2, k_3)$. Without loss of
generality if we assume $k_2\leq k_3$ then
$\si=\eta_1\cdots\eta_{k_2}\si_{k_2+1}\cdots\si_n$ and
$\tau=\eta_1\cdots\eta_{k_3}\tau_{k_3+1}\cdots\tau_n$. If $k_1<k_2$,
then $\si_{k_1+1}=\eta_{k_1+1}=\tau_{k_1+1}$ will contradict the
maximality of $k_1$.

We will denote the model of Bolthausen and Kistler as BK-GREM. The
set up in the BK-GREM is the following: For a fixed number
$n\in\mathbb{N}$, they consider the set $I=\{1,2,\cdots,n\}$ and a
collection of non-negative real numbers $\{a_J\}_{J\subset I}$ such
that $\dsum_{J\subset I} a_J =1$ with $a_\varnothing=0$. There may
be subsets $J$ of $I$ for which $a_J=0$, so they consider
$\mathcal{P}_I$ as that collection of subsets $J$ of $I$ for which
$a_J>0$ that is, $\mathcal{P}_I=\{J: \,a_J>0\}$. Like in the usual
GREM, they fix $n$ positive real numbers $\ga_i$ for $i\in I$ so
that $\dsum_{i=1}^n \ga_i =1$ and split the configurations space
$\Si_N=\{-1,1\}^N$ in to products
$\Si_{\ga_1N}\times\Si_{\ga_2N}\times\cdots\times\Si_{\ga_nN}$ with
$\Si_{\ga_iN}=\{-1,1\}^{\ga_iN}$ for each $i$. Since $\ga_iN$ may
not be an integer one needs to use $[\ga_iN]$ instead of $\ga_iN$,
and $\prod_{i=1}^n \Si_{[\ga_iN]}$ as the configuration space etc.
Since we shall soon reformulate this model we do not elaborate on
these points. So a configuration $\si$ can be written as
$(\si_1,\cdots,\si_n)$. For $J=\{j_1,\cdots,j_k\}\subset I$, denote
$\Si_{J,N}$ for $\dprod_{l=1}^k \Si_{\ga_{j_l}N}$ and $\si_J$ for the
projected configuration $(\si_j)_{j\in J}\in \Si_{J,N}$. In this
setup, the random Hamiltonian is defined as
$$H_N(\si)=\dsum_{J\in\mathcal{P}_I}\xi_J(\si_J),$$
where for $J\in \mathcal{P}_I$ and $\si_J\in\Si_{J,N}$ the random
variables $\xi_J(\si_J)$ are independent centered Gaussian random
variables with variance $a_JN$. It is quite easy to verify that, in
this model if we define the metric on the configurations space by
the same formula as in (\ref{e3.4.9-a}) then the metric will not be
always an ultrametric. But yet they have shown that the limiting
free energy is again a GREM free energy. To be precise, define a
{\em chain} $(A_0, A_1,\cdots,A_k)$ to be an increasing sequence of
subsets of $I$ with $\varnothing=A_0\subset A_1\subset \cdots\subset
A_k=I$. What they have shown is that for any BK-GREM there exists a
chain $(A_0, A_1,\cdots,A_k)$ and positive constants
$\widetilde{a}_i$ for $1\leq i\leq k$ with $\dsum_{i=1}^k
\widetilde{a}_i = 1$ such that the following holds: the free energy
of this BK-GREM is same as that of a Gaussian $k$ level GREM energy
where random variables at the $i$-th level have variance
$\widetilde{a}_iN$. This means, once again the limiting free energy
does not go beyond the GREM one. We will present here an alternative
and elegant proof of this thanks to large deviation results. To do
that we will reformulate the model in the next section.

\subsection{Reformulation}
We formulate BK-GREM as follows. Fix a set $I=\{1,2,\cdots,n\}$ with
$n\geq1$ . Let $N\geq n$ be the number of particles, each of which
can have two states/spins $+1, -1$; so that the configuration space
is $\Si_N=2^N$. Consider a partition of $N$ into integers $k(i,N),
1\leq i \leq n$ with each $k(i,N) \geq 1$ and $\dsum_i k(i,N)=N$. We
will as usual think of $2^N$ as $\prod_{i\in I} 2^{k(i,N)}$ and
$\si\in 2^N$ is $\si_1\cdots\si_n$ where $\si_i\in 2^{k(i,N)}$. Let
$S$ be the collection of non-empty subset of $I$. For each element
$s$ in $S$ we denote $2^{K_{sN}}=\dsty\prod_{i\in s}2^{k(i,N)}$.
With this notation $2^N=2^{K_{IN}}$. The map $\si\in 2^{K_{IN}}\ra \si(s)\in
2^{K_{sN}}$ is the projection map via $s$. For $s=\{ i_1, i_2,
\cdots,i_k\}\in S$ where $i_1< i_2< \cdots<i_k$ and
$\si=\si_1\cdots\si_n \in 2^{K_{IN}}$, we denote $\si(s)=\si_{i_1}
\si_{i_2}\cdots \si_{i_k}\in 2^{K_{sN}}$, the projection of $\si$
via $s$. Now for fixed $N$, we have a bunch of independent random
variables $\xi(s,\si(s))$ as $s$ varies over $S$ and $\si(s)$ varies
over $2^{K_{sN}}$.

For each $\si \in 2^N$ one can think of a lattice isomorphic to the
lattice of power set of $I$ where $\si(s)$ corresponds to the edge
of the lattice joining the nodes $s=\{ i_1, i_2, \cdots,i_k\}$ and
$\{ i_1, i_2, \cdots,i_{k-1}\}$. Now for each $\si$ associate random
variables $\xi(s,\si(s))$ to each of the lattice edge $\si(s)$. We
associate weights $a_s\geq 0$ to each edges $\si(s)$. These are not
random. In a configuration $\si=\si_1\cdots\si_n$ the Hamiltonian is
defined as
\[ H_N(\si) = N\dsum_{s\in S} a_s\xi(s,\si(s)).\]

For $\be>0$ the partition function is
\[Z_N(\be) = 2^N\mbf{E}_\si e^{-\be H_N(\si)}.\]
Here $\mbf{E_\si}$ stands for expectation with respect to $\si$ when
$2^N$ has uniform distribution. In other words, $\mbf{E_\si}$ is
simply the usual average over $\si$.

Since $\xi$'s are random variables both $H_N$ and $Z_N$ are random
variables. We suppress the parameter $\om$ that comes with the
random variables $\xi$. As usual $\frac{1}{N}\log Z_N(\be)$ is the
free energy of the $N$-particle system. As $N$ changes, the
distribution of the $\xi$'s would in general depends on $N$. So
strictly speaking we should be using superscript $N$ for the random
variables. But for ease in reading we suppress the superscript. This
should be borne in mind. We assume that all our random variables are
defined on one probability space.

\subsection{LDP Approach}

In this subsection, we outline how large deviation principle can be
used. Since we will prove a more general result in the next chapter
(see section \ref{s4.2}), we refrain from giving complete details.
Let us consider the map $\Si_N \rightarrow
\mathbb{R}^{S}$ (recall $S$ is the collection of non-empty subsets
of $I$) defined by
\[\si \mapsto \xi_\si = (\xi(s,\si(s)))_{s\in S}.\]

Let $\mu_N$ be the induced probability on $\mathbb{R}^{S}$ when
$\Si_N$ has uniform distribution, that is, each $\si \in \Si_N$ has
probability $\frac{1}{2^N}$. In other words, for any Borel set $A
\subset \mathbb{R}^{S}$,
\[\mu_N(A)=\frac{1}{2^N} \#\{\si:\xi_\si \in A\}.\]

In particular, if $A$ is a box, say $\tri=\prod_{s\in S}\tri_s$, with each $\tri_s\subseteq\mathbb{R}$, then
\[\mu_N(\tri) = \frac{1}{2^N} \dsum_{\si}
\prod\limits_{s\in S} \mbf{1}_{\tri_s}(\xi(s,\si(s))).\]

Here now is the basic observation similar to that of Theorem
\ref{t2.3.3}.
\begin{thm}\label{t3.4.1}
Let $\tri = \prod_{s\in S}\tri_s \subset \mathbb{R}^{S}$. Denote
$q_{sN} = P(\xi(s,\si(s)) \in \tri_s)$ for $s\in S$. For $t\in S$ we
denote $\prod_{s\subseteq t} q_{sN}$ by $Q_{tN}$ and $\prod_{i\in t}
k(i,N)$ by $K_{tN}$.

a) If  $\sum\limits_{N \geq 1}2^{K_{tN}} Q_{tN} < \infty$, for some
$t\in S$ then a.s. eventually, $\mu_N(\tri)= 0$.

b) If for all $t\in S$, $\sum\limits_{N\geq 1} 2^{-K_{tN}}
Q_{tN}^{-1} < \infty$, then for any $\ep >0$ a.s. eventually,
\[(1-\ep)\mbf{E}\mu_N(\tri) \leq \mu_N(\tri) \leq (1+\ep)
\mbf{E}\mu_N(\tri).\]
\end{thm}
There is no new idea, we need to verify that arguments of the previous
chapter go through.
\begin{proof}

a) Let $t$ be such that $\sum\limits_{N\geq 1}2^{K_{tN}}Q_{tN} <
\infty.$ Then
\[\begin{array}{lll}
\mu_N(\tri) &=& \dfrac{1}{2^N}\sum\limits_{\si} \prod\limits_{s\in
S} \mbf{1}_{\tri_s}(\xi(s,\si(s)))\\
&\leq&\dfrac{2^{N-K_{tN}}}{2^N}\sum\limits_{\si(t)}
\prod\limits_{s\subseteq t} \mbf{1}_{\tri_s}(\xi(s,\si(s)))\\ &=&
\dfrac{1}{2^{K_{tN}}}\sum\limits_{\si(t)} \prod\limits_{s\subseteq
t} \mbf{1}_{\tri_s}(\xi(s,\si(s)))=G_N, (\mbox{say}).
  \end{array}\]

Let $A_N = \{G_N=0\}.$ Observe that $$A_N^c=
\left\{\sum\limits_{\si(t)} \prod\limits_{s\subseteq t}
\mbf{1}_{\tri_s}(\xi(s,\si(s)))\geq 1\right\}.$$ Now by Chebyshev's
inequality,
$$\mbf{P}(A_N^c) <\mbf{E}\sum\limits_{\si(t)} \prod\limits_{s\subseteq t}
\mbf{1}_{\tri_s}(\xi(s,\si(s)))=2^{K_tN}Q_{tN}.$$ Thus by assumption
and Borel-Cantelli,  $A_N$ will occur a.s. eventually. i.e. $G_N =
0$ and hence $\mu_N(\tri) = 0$.

b)\begin{align*}
&\mbox{Var}(\mu_N(\tri))\\
=& \mbf{E}(\mu_N(\tri))^2 - (\mbf{E}\mu_N(\tri))^2\\
=& \dfrac{1}{2^{2N}}\sum\limits_{\si,\tau}
\left[{\mbf{E}}\prod\limits_{s\in S}
{\mbf{1}}_{\tri_s}(\xi(s,\si(s))){\mbf{1}}_{\tri_s}(\xi(s,\tau(s)))-
Q_{IN}^2\right]\\
=& \dfrac{1}{2^{2N}}\sum\limits_{t\in S}
\sum\limits_{\substack{\si(t)=\tau(t)\\ \si_i\neq
\tau_i,\,\forall i\in t^c}} \left[{\mbf{E}}\prod\limits_{s\subseteq t}
{\mbf{1}}_{\tri_s}(\xi(s,\si(s)))\prod\limits_{s\nsubseteq t}
{\mbf{1}}_{\tri_s}(\xi(s,\si(s))){\mbf{1}}_{\tri_s}(\xi(s,\tau(s)))- Q_{IN}^2\right]\\
\leq& \dfrac{1}{2^{2N}}\sum\limits_{t\in S}\prod\limits_{s\subseteq t}q_{sN} \prod\limits_{s\nsubseteq
t}q^2_{sN}\sum\limits_{\substack{\si(t)=\tau(t)\\
\si_i\neq \tau_i,\,\forall i\in t^c}}1\\
\leq& \dfrac{1}{2^{2N}}\sum\limits_{t\in S}\frac{Q^2_{IN}}{Q_{tN}}2^{K_{tN}}2^{2(N-K_{tN})}\\
=& \sum\limits_{t\in S} \frac{Q^2_{IN}}{2^{K_{tN}}Q_{tN}}.
\end{align*}

Hence for any $\ep >0$, by Chebyshev's inequality
\[\mbf{P}(|\mu_N(\tri)-\mbf{E}\mu_N(\tri)|>\ep \mbf{E}\mu_N(\tri))<
\frac{1}{\ep^2}\sum\limits_{t\in S}\frac{1}{2^{K_{tN}}Q_{tN}}.\]
But, in view of the assumption, the sum over N of the right side is
finite. So by Borel-Cantelli lemma, a.s. eventually,
$$(1-\ep)\mbf{E}\mu_N(\tri) \leq \mu_N(\tri) \leq (1+\ep) \mbf{E}\mu_N(\tri).$$
\end{proof}

Suppose that for each element $s\in S$, we have a sequence of
probabilities $\{\la_N^s:\; N\geq n\}$ obeying LDP with a good convex
rate function $\mathcal{I}_s(x)$. We now consider reformulated BK-model where
each $\xi(s,\si(s))$ has distribution $\la_N^s$. Thus for fixed $N$,
we have a bunch of independent random variables $\xi(s,\si(s))$ as
$s$ and $\si(s)$ vary. For example, for the $N$ particle system, one
can consider for each $s\in S$; $\xi(s,\si(s))$ to be i.i.d. having
density
\begin{equation}\label{e3.4.12}
\phi(x) =
\frac{1}{2\Gamma(\frac{1}{\ga_s})}\left(\frac{\ga_s}{N}\right)^{\frac{\ga_s-1}{\ga_s}}
e^{-N\frac{|x|^{\ga_s}}{\ga_s}} \quad -\infty< x <\infty.
\end{equation}

Let us denote
\begin{equation}\label{e3.4.13}
\Psi = \{\widetilde{x}\in \mathbb{R}^{S} :
\dsum_{s\subseteq t} \mathcal{I}_s(x_s) \leq \dsum_{i\in t} p_i
\log2, \hspace{1ex} \forall t\in S\}
\end{equation}
and the map $\mathcal{J}: \mathbb{R}^S \ra \mathbb{R}$, defined
by,
\[\begin{array}{llll}
\mathcal{J}(\tilde{x}) &=& \dsum_{s\in S} \mathcal{I}_s(x_s) &\mbox{ if }
\tilde{x} \in
\Psi\\
&=& \infty &\mbox{ otherwise}.
\end{array}\]

Then with the help of Theorem \ref{t3.4.1}, one can mimic the steps in Theorem \ref{t2.5.1} to get
\begin{thm}
In the reformulated BK-GREM, let $\frac{k(i,N)}{N}\ra p_i>0$ as
$N\ra\infty$ for $1\leq i\leq n$. Then almost surely, the sequence
$\{\mu_N,\, N\geq 1\}$ satisfies LDP with rate function
$\mathcal{J}$ defined above.
\end{thm}

In this way, once again Varadhan's lemma will ensure the existence
of the limiting free energy in this case also. Thus though we don't
have ultrametricity on the configuration space, the simple LDP
technique works.
\begin{thm}\label{t3.4.6}
In reformulated BK-GREM, almost surely
$$\lim_{N\ra\infty}\frac{1}{N} \log Z_N(\be)= \log2
-\dinf_{\widetilde{x}\in \Psi}\sum_{s\in S}\left(\be
a_sx_s+\mathcal{I}_s(x_s)\right).$$
\end{thm}

\begin{rem}
One may feel that the reformulation of BK-GREM is not exactly
similar to the original version of Bolthausen and
Kistler. They consider only those subsets $s$ of $I$ for which $a_s\neq 0$
whereas in the above reformulation all the non-empty subsets are considered.
But from the above theorem it is easy to check that in calculating
infimum, $\mathcal{I}_s$ being non-negative and $\mathcal{I}_s(0)= 0$, the terms
corresponding to those $s\in S$ for which $a_s=0$ will not contribute.
\end{rem}

In \cite{BK}, Bolthausen and Kistler identified the free energy of
this model as minimum of several GREMs associated
with, what they call, chains. We shall show that there are $n!$ many
$n$ level usual GREMs hidden in the above model. The method used
above also identifies the free energy of the BK-GREM as the minimum
of the free energies of these $n!$ GREMs. This is what we do in the
next section.

\section{Hidden Tree GREMs}\label{s3.5}

In this section, we consider the BK-GREM, that is, we take $S$ to be
the set of all increasing sequences of elements of $I$ with the
Gaussian driving distributions. As mentioned earlier, this is
nothing but the Bolthausen-Kistler's model since such sequences
correspond to non-empty  subsets of $I$. Suppose now for $s\in S$
the associated weight is $a_s$. Here we evaluate the explicit
expression for the limiting free energy of BK-GREM.  Though it is
possible to consider different driving distributions for each $s\in
S$, a general closed form expression appears to be difficult. Of course,
we could also start with some more general driving distributions than
Gaussian, like distribution having density $\phi$ as in
(\ref{e3.4.12}) with $\ga>1$ at all the levels. Since in that case, there is no new
idea needed, we restrict ourselves to Gaussian case for
notational simplicity. It is worth mentioning here that, in
\cite{BK}, Bolthausen and Kistler evaluate the expression of the
limiting free energy in two steps. In the first step, they define a {\em
chain} as a sequence of strictly increasing sequences of subsets
$(A_0, A_1,\cdots,A_K)$ of $I$ so that $\varnothing=A_0\subset
A_1\subset\cdots\subset A_K=I$. For such a chain they associated a $K$
level GREM with appropriate weights calculated from the weights of the original model.
Then by second moment estimates, they have shown that the limiting
free energy of each such GREM associated to a chain is an almost
sure upper bound for the limiting free energy of their model. In the
second step, they constructed a chain in which the free energy
of the BK-GREM is attained.

Here we get the expression for the limiting free energy by
calculating
\begin{equation}
\inf_{\widetilde{x}\in\Psi} \sum_{s\in S} \left(\be a_s
x_s+\frac{1}{2}x_s^2\right),
\end{equation}
where $$\Psi=\{\widetilde{x}\in \mbb{R}^S:\; \sum_{t\subseteq
s}x_t^2\leq \sum_{i\in s} 2 p_i \log2,\; \forall s\in S\}.$$

Note that $\Psi$ is same as that of (\ref{e3.4.13}) with
$\mathcal{I}_s(x_s)=\frac{1}{2} x_s^2$. As earlier, the above
infimum is same as
\begin{equation}\label{e3.5.a18}
\inf_{\widetilde{x}\in\Psi^+} \sum_{s\in S}
\left(\frac{1}{2}x_s^2-\be a_s x_s\right),
\end{equation}
where
\begin{equation}\label{e3.5.b18}
\Psi^+=\{\widetilde{x}\in \mbb{R}^S:\; \sum_{t\subseteq s}x_t^2\leq
\sum_{i\in s} 2 p_i \log2\;\&\; x_s\geq 0,\; \forall s\in S\}.
\end{equation}
To evaluate (\ref{e3.5.a18}), we start with some notations. Here the ideas are very
much similar to that of Bolthausen and Kistler. A new idea is the
introduction of permutations of $\{1,2,\cdots,n\}$ justifying the
title of this section. For $A\subseteq I$, let us define
$$p_A=\sum_{i\in A} p_i$$ and $$w^2_A=\dsum_{\substack{s\subseteq A\\s\in
S}} a_s^2,$$ with $w^2_\varnothing=0$. Let $\mathcal{P}_0=\mathcal{P}_I$
denote the set of  permutations of $I$. With this notation, for $\pi \in
\mathcal{P}_I$ and $0\leq i < j \leq n$, denote
\begin{equation}\label{e3.5.c18}
B_{ij}^{\pi}
=\sqrt{\frac{2(p_{\pi(i+1)}+\cdots+p_{\pi(j)})\log2}{w^2_{\{\pi(1),\cdots,
\pi(j)\}}-w^2_{\{\pi(1),\cdots, \pi(i)\}}}},
\end{equation}
where, for $i=0$ the set $\{\pi(1),\cdots, \pi(i)\}$ that appears in the denominator is treated as the empty set.

Note that earlier to evaluate the explicit energy expression for
$\gamma$-GREM with $\ga>1$, in subsection \ref{s2.6.1}, we consider
only one triangular array of numbers defined as $B(j,k)$ in
(\ref{e2.5.a12}) for $1\leq j\leq k\leq n$. Now here we are
considering $n!$ many triangular arrays corresponding to each
permutation $\pi$.

Now let,
\begin{equation}\label{e3.5.d18}
\be_1=\dmin_{\pi\in\mathcal{P}_I}\dmin_{0< j\leq n}
B_{0j}^\pi=\dmin_{(\pi,j)} B_{0j}^\pi.
\end{equation}

Also note that, in subsection \ref{s2.6.1}, we define $\be_1$ in
(\ref{e2.5.b12}) as the minimum over all the entries in the first
row of the triangular array $B(j,k)$ in (\ref{e2.5.a12}). Here we
are defining $\be_1$ as the minimum over all the entries in the
first rows of all the $n!$ triangular arrays.

It may be further noted that in subsection \ref{s2.6.1}, to define
$\be_1$ if the minimum occurred at two places, we  had taken the
maximum index (see (\ref{e2.5.a12}) for the definition of $r_l$). We
now implement the same plan in the present setting also. Since the
minimum may be attained in the first lines of two different
triangular arrays, one needs to know what is meant by maximum index.
This will be done now.

Suppose the minimum occurs at two places, say at $(\theta, k)$ and
$(\varrho, l)$, that is, $\be_1= B_{0k}^{\theta}= B_{0l}^{\varrho}$.
Let $G=\{\theta(1),\cdots,\theta(k)\}$;
$H=\{\varrho(1),\cdots,\varrho(l)\}$. Let $|G\cup H|=m$ and let
$\mathcal{P}_1$ denote the class of all permutations of $I$ for
which $\{\pi(1),\cdots,\pi(m)\}=G\cup H$. So $\mathcal{P}_1\subset
\mathcal{P}_I$.

\noindent{\bf Claim:} {\em
$B_{0m}^\pi=\sqrt{\dfrac{2p_{G\cup H}\log2}{w_{G\cup H}^2}}=\be_1$, for
every $\pi\in\mathcal{P}_1$.}

To justify the claim, first of all note that $$w_{G\cup H}^2\geq
w_{G}^2+w_{H}^2-w_{G\cap H}^2,$$ whereas,
$$p_{G\cup H}=p_{G}+p_{H}-p_{G\cap H}.$$
Then for $\pi\in \mathcal{P}_1$,
\begin{align*}
\;& 2p_{G\cup H}\log2-\be_1^2 w_{G\cup H}^2\\
\leq \;& 2(p_{G}+p_{H}-p_{G\cap H})\log2 -
\be_1^2\left(w_{G}^2+w_{H}^2-w_{G\cap
H}^2\right)\\
= \;& \left(2p_{G}\log2-
\be_1^2w_{G}^2\right)+\left(2p_{H}\log2-
\be_1^2w_{H}^2\right)+\left(\be_1^2w_{G\cap
H}^2-p_{G\cap H}\right)\\
\intertext{as $\be_1= B_{0k}^{\theta}=
B_{0l}^{\varrho}$, first two terms are zero,}
= \;& \be_1^2w_{G\cap H}^2-p_{G\cap H}\\
\leq \;& 0.
\end{align*}
The last inequality follows from the fact that $\be_1$ is obtained
by taking the minimum over all possible choice of $(\pi,j)$. This shows
$\be_1^2 \geq \frac{2p_{G\cup H}\log2}{w_{G\cup H}^2}$.

Once again
$\be_1$ being the minimum over all possible choice of $(\pi,j)$, we
conclude that $\be_1^2$ actually equals $\frac{2p_{G\cup H}\log2}{w_{G\cup H}^2}$ proving
the claim.
\medskip

If the minimum in (\ref{e3.5.d18}) occurs at more than two places,
we can use induction to conclude that  there exists a unique maximal
set, say, $G_1\subseteq I$ such that the following holds. Let
$|G_1|=l_1$ and $\mathcal{P}_1=\text{ all permutations that map
$\{1,2,\cdots l_1\}$ on to $G_1$}$. Then for any $\pi\in
\mathcal{P}_1$,  $B_{0l_1}^\pi=\be_1$.

It may happen that $G_1=I$, then we will stop. Otherwise, let us
define \begin{equation*}\be_2=\dmin_{\pi\in\mathcal{P}_1}\dmin_{l_1<
j\leq n} B_{l_1j}^\pi=\dmin_{(\pi,j)} B_{l_1j}^\pi.\end{equation*}
Of course, the last minimum is only over $\pi\in \mathcal{P}_1$.

Once again going back to subsection \ref{s2.6.1}, to define $\be_2$
in (\ref{e2.5.b12}), we look only the entries from the $r_1+1$-th
row of the triangular array $B(j,k)$ in (\ref{e2.5.a12}). Here, we
are looking the entries of $l_1+1$-th rows (as $0$ corresponds the
first row) of the triangular arrays corresponding to each $\pi\in
\mathcal{P}_1$.

If possible, suppose the minimum occurs at two places, say, at
$(\theta, k)$ and at $(\varrho, l)$. So $\theta,
\varrho\in\mathcal{P}_1$; \, $l_1< k, l\leq n$ and $\be_2=
B_{l_1k}^{\theta}= B_{l_1l}^{\varrho}$. Let
$G=\{\theta(1),\cdots,\theta(k)\}$;
$H=\{\varrho(1),\cdots,\varrho(l)\}$ and $G_2=G\cup H$. Let
$|G_2|=l_2$ and denote by $\mathcal{P}_2$, the class of all
permutations of $\mathcal{P}_1$  for which
$\{\pi(1),\cdots,\pi(l_2)\}=G_2$. Since $\pi$ in $\mathcal{P}_1$
already maps $\{1,2,\cdots,l_1\}$ onto $G_1$, this extra condition
only means that $\pi$ moreover maps $\{l_1+1,\cdots,l_2\}$ onto
$G_2-G_1$. Clearly, $\mathcal{P}_2\subset \mathcal{P}_1$.

\noindent{\bf Claim:}{\em \; $B_{l_1l_2}^\pi=\sqrt{\dfrac{2p_{G_2-G_1}\log2}{w_{G_2}^2-w_{G_1}^2}}=\be_2$,
for every $\pi\in\mathcal{P}_2$.}

The justification of this claim is similar to that of the earlier
one. Once again note that, $$w_{G_2}^2\geq
w_{G}^2+w_{H}^2-w_{G\cap H}^2,$$ but
$$p_{G_2-G_1}=p_{G-G_1}+p_{H-G_1}-p_{G\cap H-G_1}.$$ So for $\pi\in
\mathcal{P}_2$, we have

$$\begin{array}{l}
2p_{G_2-G_1}\log2-\be_2^2
\left(w_{G_2}^2-w_{G_1}^2\right)\vspace{1ex}\\

\leq 2\left(p_{G-G_1}+p_{H-G_1}-p_{G\cap
H-G_1}\right)\log2 -
\be_2^2\left(w_{G}^2+w_{H}^2-w_{G\cap
H}^2-w_{G_1}^2\right)\vspace{1ex}\\

= 2\left(p_{G-G_1}+p_{H-G_1}-p_{G\cap
H-G_1}\right)\log2 - \vspace{1ex}\\

\hspace{35ex}\be_2^2\left(w_{G}^2-w_{G_1}^2+w_{H}^2-w_{G_1}^2-w_{G\cap
H}^2+w_{G_1}^2\right)\vspace{1ex}\\

=\left(2p_{G-G_1}\log2-
\be_2^2\left(w_{G}^2-w_{G_1}^2\right)\right)+\left(2p_{H-G_1}\log2-
\be_2^2\left(w_{H}^2-w_{G_1}^2\right)\right)+\vspace{1ex}\\

\hspace{35ex}\left(\be_2^2\left(w_{G\cap
H}^2-w_{G_1}^2\right)-2p_{G\cap
H-G_1}\log2\right)\vspace{1ex}\\

=\be_2^2\left(w_{G\cap H}^2-w_{G_1}^2\right)-
2p_{G\cap H-G_1}\log2\\

\leq 0.
\end{array}$$

Hence $\be_2^2 \geq \dfrac{2p_{G_2-G_1}\log2}{w_{G_2}^2-w_{G_1}^2}$
and once again $\be_2$ being the minimum over all possible choice of
$\pi\in\mathcal{P}_1$ and $l_1< j\leq n$, the only possibility left
is the equality. That is $\be_2^2 =
\dfrac{2p_{G_2-G_1}\log2}{w_{G_2}^2-w_{G_1}^2}$ and hence the claim
is proved.

If the minimum occurs at more than two places, we can use induction
to conclude that  there exists a unique maximal set, say, $G_2$,
such that $G_1\subset G_2\subseteq I$ and the following holds. Let
$|G_2|=l_2$ and $\mathcal{P}_2$ be all permutations of
$\mathcal{P}_I$ that map $\{l_1+1,\cdots,l_2\}$ onto $G_2-G_1$. Then
$B_{l_1l_2}^\pi(G_1)=\be_2$ for all $\pi\in \mathcal{P}_2$. Of
course, all the quantity $l_2$, $\be_2$ depend on $G_1$.

Proceeding by induction can summarize:

There is a (unique) integer $K$ with $1\leq K\leq n$ and for every
$i$ with $1\leq i\leq K$ there are $\be_i$, $l_i$, $G_i$ and
$\mathcal{P}_i$ satisfying the following:
\begin{enumerate}
\item{$\varnothing= G_0\subset G_1\subset\cdots\subset G_K=I$ with
$|G_i|=l_i$ so that $1\leq l_1<l_2<\cdots<l_K=n$.}

\item{$\mathcal{P}_i$ is the set of permutations $\pi$ of $I$ that
maps $\{1, 2, \cdots, l_j\}$ onto $G_j$ for each $j\leq i$ so that
$\mathcal{P}_1\supset\mathcal{P}_2\supset\cdots\supset\mathcal{P}_K$.}

\item{$\be_i=B_{l_1l_2}^\pi$ for every $\pi\in\mathcal{P}_i$ and
this common value is also same as $\dmin_{\pi\in\mathcal{P}_{i-1}}\dmin_{l_{i-1}< j\leq n}
B_{l_{i-1}j}^\pi$.}
\end{enumerate}

So note that for any $\pi\in \mathcal{P}_K$, we can trace out the
$\beta_i$ for $1\leq i \leq K$, as
\begin{equation}\be_i=\sqrt{\dfrac{2(p_{\pi(l_{i-1}+1)}+\cdots+p_{\pi(l_i)})\log2}{w^2_{\{\pi(1),\cdots,
\pi(l_i)\}}-w^2_{\{\pi(1),\cdots, \pi(l_{i-1})\}}}}.\end{equation}

Moreover, the infimum in (\ref{e3.5.a18}) reduces to the following:
\begin{equation}\label{e3.5.18}
\dinf_{\widetilde{x}\in \Psi}\sum_{s\in S}\left(\frac{1}{2}x_s^2-\be
a_sx_s\right), \end{equation} and
\begin{equation}\Psi=\left\{\widetilde{x}:\,
X_C^2\leq 2p_C\log2,\, \forall (\varnothing\neq)C\subseteq
I\right\}\end{equation} where we used the notation
\begin{align*}X_C^2=\dsum_{s\subseteq C,
s\in S}x_s^2 &\quad \text{and}\quad p_C=\sum_{i\in C}
p_i.\end{align*}

Now we prove that, if $\be_j\leq \be <\be_{j+1}$, the above infimum
is attained at $\widetilde{x}^*=(x_s^*;s\in S)\in \mbb{R}^S$ given
by
\[x_s^*=\begin{cases}
               \be_1 a_s & \mbox{if $s\subseteq G_1$}\\
               \be_2 a_s & \mbox{if $s\subseteq G_2, s\nsubseteq G_1$}\\
               \vdots &\\
               \be_j a_s & \mbox{if $s\subseteq G_j, s\nsubseteq G_{j-1}$}\\
               \be a_s &   \mbox{if $s\nsubseteq G_j$}.
              \end{cases}
\]

First of all note that $\widetilde{x}^*\in \Psi$. For, $C\subseteq
I$ implies
\begin{align*}
&\quad{X^*_C}^2\\
&=\dsum_{s\subseteq C, s\in S} {x_s^*}^2\\
\intertext{}
&=\dsum_{i=1}^{(j+1)\wedge K} \dsum_{\substack{s\subseteq G_i,s\nsubseteq G_{i-1}\\s\subseteq C, s\in S}}{x_s^*}^2\\
&=\dsum_{i=1}^j\dsum_{\substack{s\subseteq G_i,s\nsubseteq G_{i-1}\\s\subseteq C, s\in S}}\be_i^2a_s^2\;+\quad 1_{\{j+1\leq K\}}\dsum_{\substack{s\nsubseteq G_j\\s\subseteq C, s\in S}}\be^2a_s^2\\
&\leq\dsum_{i=1}^j\dsum_{\substack{s\subseteq G_i,s\nsubseteq G_{i-1}\\s\subseteq C, s\in S}}\be_i^2a_s^2\;+\quad 1_{\{j+1\leq K\}}\dsum_{\substack{s\nsubseteq G_j\\s\subseteq C, s\in S}}\be_{j+1}^2a_s^2\\
&=\dsum_{i=1}^j\dsum_{\substack{s\subseteq G_i,s\nsubseteq
G_{i-1}\\s\subseteq C, s\in
S}}\frac{2p_{G_i-G_{i-1}}\log2}{w^2_{G_i}-w^2_{G_{i-1}}}
a_s^2\;+\quad 1_{\{j+1\leq K\}}\dsum_{\substack{s\nsubseteq G_j\\s\subseteq C, s\in S}}\frac{2p_{G_{j+1}-G_j}\log2}{w^2_{G_{j+1}}-w^2_{G_j}}a_s^2\\
&\leq\dsum_{i=1}^j\dsum_{\substack{s\subseteq G_i,s\nsubseteq
G_{i-1}\\s\subseteq C, s\in S}}\frac{2p_{(C\cup
G_{i-1})-G_{i-1}}\log2}{w^2_{C\cup G_{i-1}}-w^2_{G_{i-1}}}
a_s^2\;+\quad 1_{\{j+1\leq K\}}\dsum_{\substack{s\nsubseteq G_j\\s\subseteq C, s\in S}}\frac{2p_{(C\cup G_j)-G_j}\log2}{w^2_{C\cup G_j}-w^2_{G_j}}a_s^2\\
&\leq\dsum_{i=1}^j2p_{(C\cup G_{i-1})-G_{i-1}}\log2 \;+\quad 1_{\{j+1\leq K\}}2p_{(C\cup G_j)-G_j}\log2\\
&= 2p_C\log2.
\end{align*}

Secondly, note that for any $\widetilde{x}\in \Psi$, we have
$$\dsum_{i=1}^j\dsum_{s\subseteq G_i,s\nsubseteq G_{i-1}} \be
a_s(x^*_s - x_s) \geq \dsum_{i=1}^j\dsum_{s\subseteq G_i,s\nsubseteq
G_{i-1}} x^*_s(x^*_s - x_s).$$ For, by Holder's inequality we have
$$\dsum_{i=1}^j\dsum_{s\subseteq G_i,s\nsubseteq G_{i-1}} x^*_s x_s
\leq \sqrt{\dsum_{i=1}^j\dsum_{s\subseteq G_i,s\nsubseteq G_{i-1}}
{x^*_s}^2} \sqrt{\dsum_{i=1}^j\dsum_{s\subseteq G_i,s\nsubseteq
G_{i-1}} x_s^2}\leq \dsum_{i=1}^j\dsum_{s\subseteq G_i,s\nsubseteq
G_{i-1}} {x^*_s}^2,$$ where the last inequality follows from the
fact that $\widetilde{x}\in \Psi$. Hence
$$\dsum_{i=1}^j\dsum_{s\subseteq G_i,s\nsubseteq G_{i-1}}
x^*_s(x^*_s-x_s) \geq 0.$$ Since $\be>\be_j$, we have
$\left(\frac{\be}{\be_i}-1\right)>0$ for $1\leq i\leq j$. Moreover
$\be_i$ being increasing in $i$, these numbers
$\left(\frac{\be}{\be_i}-1\right)$ are decreasing and hence we get,
$$\dsum_{i=1}^j\left(\frac{\be}{\be_i}-1\right)\dsum_{s\subseteq G_i,s\nsubseteq G_{i-1}}
x^*_s(x^*_s-x_s) \geq 0.$$ In other words using the definition of
$x^*_s$ we get the observation $$\dsum_{i=1}^j\dsum_{s\subseteq
G_i,s\nsubseteq G_{i-1}} \be a_s(x^*_s - x_s) \geq
\dsum_{i=1}^j\dsum_{s\subseteq G_i,s\nsubseteq G_{i-1}} x^*_s(x^*_s
- x_s).$$

Now by using the above inequality, we have
$$\begin{array}{l}\dsum_{s\subseteq G_j}\left(\frac{1}{2}x_s^2 - \be
a_s x_s\right)-\dsum_{s\subseteq G_j}\left(\frac{1}{2}{x_s^*}^2 -
\be a_s
x_s^*\right)\\
=\dsum_{s\subseteq G_j}\left(\frac{1}{2}x_s^2 - \be a_s
(x_s-x_s^*)-\frac{1}{2}x_s^2\right)\\
\geq \dsum_{s\subseteq G_j}\left(\frac{1}{2}x_s^2 - x_s^*
(x_s-x_s^*)-\frac{1}{2}x_s^2\right)\\
=\frac{1}{2}\dsum_{s\subseteq
G_j}\left(x_s - x_s^* \right)^2\geq 0.
\end{array}$$
Moreover, using the definition of $x_s^*$ for $s\nsubseteq G_j$, we
have
$$\begin{array}{l}\dsum_{s\nsubseteq G_j}\left(\frac{1}{2}x_s^2 -
\be a_s x_s\right)-\dsum_{s\nsubseteq G_j}\left(\frac{1}{2}{x_s^*}^2
- \be a_s
x_s^*\right)\\
=\dsum_{s\nsubseteq G_j}\left(\frac{1}{2}x_s^2 - \be a_s
x_s +\frac{1}{2}\be^2a_s^2\right)\\
=\frac{1}{2}\dsum_{s\nsubseteq G_j}\left(x_s - \be x_s\right)^2\geq
0.
\end{array}$$
Thus combining the above two inequality, $$\dsum_{s\in
S}\left(\frac{1}{2}x_s^2 - \be a_s x_s\right)-\dsum_{s\in
S}\left(\frac{1}{2}{x_s^*}^2 - \be a_s x_s^*\right)\geq0$$ and hence
the infimum occurs at $\widetilde{x}^*$.

Denote $\be_0=0$ and $\be_{K+1}=\infty$. Suppose $1\leq j\leq K$ and
 $\be\in \left[ \be_j,\be_{j+1}\right)$ then the infimum in (\ref{e3.5.18}) becomes
$$\begin{array}{l}
\dsum_{i=1}^j \dsum_{s\subseteq G_i , s\nsubseteq G_{i-1}}
\left(\frac{1}{2}\be_i^2a_s^2 -\be \be_i a_s^2\right)
+\dsum_{s\nsubseteq G_j} \left(\frac{1}{2}\be^2a_s^2 -\be^2
a_s^2\right)\\
=\dsum_{i=1}^j\frac{1}{2}\be_i^2 \dsum_{s\subseteq G_i , s\nsubseteq
G_{i-1}} a_s^2 -\be\dsum_{i=1}^j \be_i\dsum_{s\subseteq G_i ,
s\nsubseteq G_{i-1}} a_s^2 -\frac{1}{2}\be^2\dsum_{s\nsubseteq G_j}
a_s^2\\
= p_{G_j}\log2 -\be\dsum_{i=1}^j \be_i\dsum_{s\subseteq G_i ,
s\nsubseteq G_{i-1}} a_s^2 -\frac{1}{2}\be^2\dsum_{s\nsubseteq G_j}
a_s^2
\end{array}$$

We can summarize the above discussion in the following
\begin{thm}
In the Gaussian BK-GREM, almost surely,
$$\lim_N \frac{1}{N}\log Z_N(\be)=\dsum_{i\notin G_j} p_i\log2 +\be\dsum_{i=1}^j \be_i\dsum_{s\subseteq G_i
, s\nsubseteq G_{i-1}} a_s^2 +\frac{1}{2}\be^2\dsum_{s\nsubseteq
G_j} a_s^2,$$ if $\be\in \left[ \be_j,\be_{j+1}\right)$ for $0\leq
j\leq K$.
\end{thm}
We shall now describe for each $\pi\in \mathcal{P}_I$ an $n$ level
GREM. In what follows $\pi\in \mathcal{P}_I$ is fixed. For the $N$
particle system there are $2^{k(\pi(i),N)}$ furcations at the $i$th
level, below each node of the $(i-1)$th level. The weights at the
$i$-th level in this GREM are $w(\pi,i)$ which are defined by
$w(\pi,1)=a_{\pi(1)}$, and in general, for $1\leq i\leq n$
\begin{equation}\label{e3.5.21}w^2(\pi,i)=\dsum_{\substack{s\subseteq\{\pi(1),\cdots,\pi(i)\}\\
s\nsubseteq\{\pi(1),\cdots,\pi(i-1)\}}} a_s^2.\end{equation}

Let $\mathcal{E}(\pi,\be)$ be the almost sure limiting free energy
of this GREM. This exists by Theorem \ref{t2.4.6}. As done in
subsection \ref{s2.6.1}, we set $r_0^\pi=0$ and let
$$\be_{i}^\pi = \dmin_{k>r_{i-1}}B_{r^\pi_{i-1},k}^\pi$$ with $r_i^\pi = \max\{l>r_{i-l}:
B(r^\pi_{i-1},l)=\be_{i}\}$ for $1\leq i\leq K^\pi$ with
$r_{K^\pi}=n$. Also denote $\be_0^\pi=0$ and
$\be_{K^\pi+1}^\pi=\infty$. Then by Theorem \ref{t2.6.1}, we have
for $\be\in \left[ \be_j^\pi,\be_{j+1}^\pi\right)$ with $0\leq j\leq
K^\pi$, \begin{equation}\mathcal{E}(\pi,\be)=\dsum_{i=r_j^\pi+1}^n
p_{\pi(i)}\log2 + \frac{1}{2}\be^2\dsum_{i=r_j^\pi+1}^n w^2(\pi,i)
+\be\dsum_{i=1}^{r_j^\pi} \be_i^\pi w^2(\pi,i).\end{equation}

Now let us consider $\pi\in\mathcal{P}_K$. Then note that
$l_i=r_i^\pi$ for all $1\leq i\leq K$ and by definition $\be_i^\pi$
is same as $\be_i$ . Hence $\dsum_{i=r_j^\pi+1}^n
p_{\pi(i)}=\dsum_{i\notin G_j} p_i$; $\dsum_{s\subseteq G_i ,
s\nsubseteq G_{i-1}} a_s^2=w^2(\pi,i)$ and $\dsum_{s\nsubseteq G_j}
a_s^2=\dsum_{i=r_j^\pi+1}^n w^2(\pi,i)$, so that
$\mathcal{E}(\be)=\mathcal{E}(\pi,\be)$.

Thus for every $\pi\in\mathcal{P}_K$, the GREM associated in the
above paragraph has the same energy, namely, $\mathcal{E}(\be)$, the
energy of the BK-GREM.

We now go on to show that if $\pi$ is any permutation then the
energy of the GREM associated with $\pi$, namely
$\mathcal{E}(\pi,\be)$, is larger than $\mathcal{E}(\be)$. So fix a
permutation $\pi$.

Denote $H_i^\pi=\{s:\; s\subseteq\{\pi(1),\cdots,\pi(i)\} \&
s\nsubseteq\{\pi(1),\cdots,\pi(i-1)\}\}$, that is, $H_i^\pi$
consists of all subsets of $\{\pi(1),\cdots,\pi(i)\}$ that include
$\pi(i)$. Then
\begin{align*}
&\dsum_{s\in S} \left(\frac{1}{2}x_s^2 -\be a_s x_s\right)\\
=&\dsum_{i=1}^n \dsum_{\substack{s\in S, s\subseteq\{\pi(1),\cdots,\pi(i)\}\\
s\nsubseteq\{\pi(1),\cdots,\pi(i-1)\}}}\left(\frac{1}{2}x_s^2
-\be a_s
x_s\right)\\
\geq&  \dsum_{i=1}^n \left(\dsum_{\substack{s\subseteq\{\pi(1),\cdots,\pi(i)\}\\
s\nsubseteq\{\pi(1),\cdots,\pi(i-1)\}}}\frac{1}{2}x_s^2
-\be \left(\dsum_{\substack{s\subseteq\{\pi(1),\cdots,\pi(i)\}\\
s\nsubseteq\{\pi(1),\cdots,\pi(i-1)\}}}a_s^2\right)^\frac{1}{2}\left(\dsum_{\substack{s\subseteq\{\pi(1),\cdots,\pi(i)\}\\
s\nsubseteq\{\pi(1),\cdots,\pi(i-1)\}}}
x_s^2\right)^\frac{1}{2}\right)\\
\intertext{since for $C\subseteq I$,  $\dsum_{s\in S,s\subseteq
C}a_s x_s\leq \left(\dsum_{s\in S,s\subseteq
C}a_s^2\right)^{\frac{1}{2}}\left(\dsum_{s\in S,s\subseteq
C}x_s^2\right)^{\frac{1}{2}}=w_CX_C$,} =&\dsum_{i=1}^n
\left(\frac{1}{2}X_{H_i^\pi}^2-\be w(\pi,i)X_{H_i^\pi}\right).
\end{align*}

Moreover, for $\pi\in \mathcal{P}_I$, let us denote
$$\begin{array}{rl}
\Psi_\pi
&=\left\{X_{\{\pi(1),\cdots,\pi(i)\}}^2\leq
2p_{\{\pi(1),\cdots,\pi(i)\}}\log2, \quad 1\leq i \leq
n\right\}\\
&=\left\{\dsum_{i=1}^k X_{H_i^\pi}^2\leq \dsum_{i=1}^k
2p_{\pi(i)}\log2, \quad \forall 1\leq k \leq n\right\}\subset
\mbb{R}^S.\\\end{array}$$ Then $\Psi \subseteq \Psi_\pi$ for every
$\pi\in \mathcal{P}_I$. Hence for every $\pi\in \mathcal{P}_I$, we
have
$$\dinf_\Psi \dsum_{s\in S} \left(\frac{1}{2}x_s^2 -\be a_s
x_s\right) \geq \dinf_{\Psi_\pi}\dsum_{i=1}^n \left(
\frac{1}{2}X_{H_i^\pi}^2-\be w(\pi,i)X_{H_i^\pi}\right)$$ and hence
$$\begin{array}{rll}\mathcal{E}(\be)&=\log2 -\dinf_\Psi
\dsum_{s\in S} \left(\frac{1}{2}x_s^2 -\be a_s x_s\right)&\\
&\leq \log2 - \dinf_{\Psi_\pi}\dsum_{i=1}^n
\frac{1}{2}X_{H_i^\pi}^2-\be \be
w(\pi,i)X_{H_i^\pi}&=\mathcal{E}(\pi,\be).\end{array}$$ Thus we have
proved the following.
\begin{thm}\label{t3.5.2}
Almost surely,
$$\mathcal{E}(\be)=\dinf_\pi\mathcal{E}(\be,\pi).$$
\end{thm}
That is, the free energy of the Gaussian BK-GREM represents the free
energy of an  $n$ level tree GREM. In fact, it represents the
minimum out of all possible $n!$ many $n$-level Gaussian tree GREM
energies with appropriately defined weights.

\begin{rem}
A closer look of the definition of $B_{ij}^\pi$ reveals that if
$a_s=0$ for some $s\in S$ then such an $s$ plays no role in the
definition of $\be_i$s. Moreover, since $\sum_{s\in
S}\left(\frac{1}{2}x_s^2-\be a_sx_s\right)= \sum_{s\in S, a_s\neq
0}\left(\frac{1}{2}x_s^2-\be a_sx_s\right)+\sum_{s\in S, a_s=
0}\frac{1}{2}x_s^2$, in calculating infimum of $\sum_{s\in
S}\left(\frac{1}{2}x_s^2-\be a_sx_s\right)$ on $\psi$ we will be
quite justified to put $x_s=0$ for all those $s\in S$ for which
$a_s=0$. This will lead to the calculation of infimum of $\sum_{s\in
S, a_s\neq 0}\left(\frac{1}{2}x_s^2-\be a_sx_s\right)$ on $\Psi$. In
other words, we could consider $S$ to be the collection of all those
increasing sequences $s$ for which $a_s\neq 0$, instead of all
sequences. This is the setup of the original Bolthausen-Kistler
model.

Bolthausen and Kistler have shown that the free energy is the
minimum among the free energies of the tree GREMs associated with
the all possible increasing chains of subsets of $I$. What the above argument shows is that one need
not consider all chains. It is enough to consider $n!$ many $n$
level GREMs. Can we reduce $n!$? Perhaps not in general.
Incidentally, the argument also identifies all these $n$ level GREMs
which attains the minimum. In fact, the number of such $n$ level
GREM is precisely $|\mathcal{P}_K|$, cardinality of $\mathcal{P}_K$.
\end{rem}

\begin{rem}
Though the BK-model is not a hierarchial model, yet when the driving
distribution is Gaussian we are not able to get out of the tree
GREM. That is, the tree GREM is in some way hidden in this model.

Now one can raise the question whether going out of Gaussian driving
distributions leads to a BK-GREM that is not usual tree GREM (in the
sense of energy). In this regard, it is worth mentioning, that if we
consider that the driving distributions $\{\mu_N^s\}_N$ satisfying
LDP with rate function $\mathcal{I}_s(x_s)=\frac{1}{\ga}x_s^\ga$ for
some $\ga\geq 1$ and every $s\in S$ then Theorem \ref{t3.5.2}
remains true. To see this we follow the same line of proof as above
with the appropriate changes as done in section 2.6.
\end{rem}

\begin{rem}
Large deviation approach allows us to consider different driving
distributions for different $s\in S$. This can be done with BK-GREM
also and one can prove the existence of free energy. But it is not
easy to obtain explicit formula.
\end{rem}

\section{Block Tree GREM}

In the previous section we have shown that in the Gaussian BK-Model
the limiting free energy is the minimum over all possible $n!$ many
$n$-level tree GREMs with appropriate weights. Now we will conclude
this chapter by exhibiting one model which includes again $n!$ many
$n$-level GREMs and where the free energy is maximum over all those
GREMs. To define the model we will use the notation $n, N,
\si=\si_1\cdots\si_n$ with the same interpretation as that of the
earlier section. Let $a_1,\cdots,a_n$ be given non-negative weights.
For any sequence $s=\left<j_1,\cdots,j_i\right>$ of distinct
elements from $I=\{1,2,\cdots,n\}$ and for any
$\si(s)=\left<\si_{j_1},\cdots,\si_{j_i}\right>\in
2^{k(j_1,N)}\times\cdots\times 2^{k(j_i,N)}$. We  have random
variables $\xi_{\si(s)}^s$ and these are independent
$\mathcal{N}(0,N)$. Now depending on $\pi$, a permutation of $I$ and
$\si\in 2^N$, we define the Hamiltonian as
\begin{equation}H_N(\sigma,\pi)=\sum_{i=1}^n a_{\pi(i)}
\xi_{\sigma_{\pi(1)}\sigma_{\pi(2)}
\cdots\sigma_{\pi(i)}}^{\pi(1)\pi(2)\cdots\pi(i)}.\end{equation}

Note that here the configuration space has $n!\times 2^N$ many
points instead of usual $2^N$ many. We call this model as Block tree
GREM. We define the partition function corresponding to inverse
temperature $\be>0$ as
$$Z_N(\be)=\sum_{\pi\in\mathcal{P}_I} \sum_{\si\in \Si_N}e^{-\be
H_N(\sigma,\pi)},$$ and the definition of free energy is
$\frac{1}{N}\log Z_N(\be)$.

\begin{figure}[ht]
 \centering\includegraphics{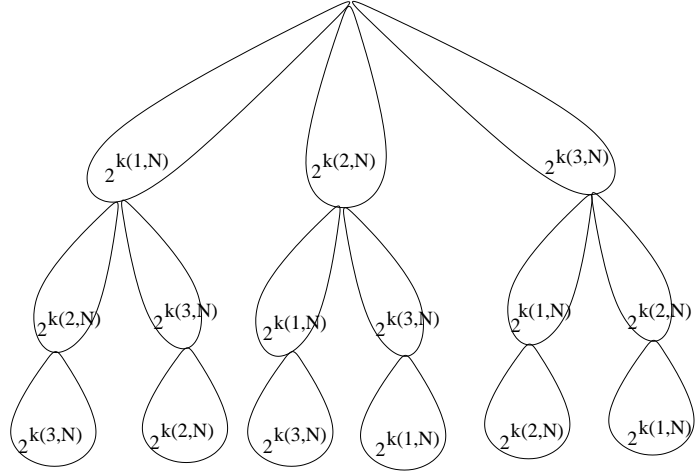}
 \caption{Block Tree GREM}\label{f3.1}
\end{figure}

So for $n=3$ the model will look like as in Figure \ref{f3.1}.

Now for each $\pi\in\mathcal{P}_I$, let us denote $$Z_N^\pi(\be)=
\sum_{\si\in \Si_N}e^{-\be H_N(\sigma,\pi)}.$$ Note that for each
$\pi$, $Z_N^\pi$ denote the partition function for the $n$-level
tree GREM with $2^{k(\pi(i),N)}$ furcations below each of node at
the $(i-1)$-th level of the tree and with the associated weight in
the $i$-th level being $a_{\pi(i)}$. So we can write
$$Z_N(\be)=\sum_{\pi\in\mathcal{P}_I} Z_N^\pi(\be).$$ Hence
$\frac{1}{N}\log Z_N(\be)=\frac{1}{N}\log\sum_{\pi\in\mathcal{P}_I}
Z_N^\pi(\be)=\frac{1}{N}\log\dmax_\pi
Z_N^\pi(\be)+\frac{1}{N}\log\sum_{\pi\in\mathcal{P}_I}\frac{Z_N^\pi(\be)}{\dmax_\pi
Z_N^\pi(\be)}$.

Since limiting free energy exists almost surely corresponding to
every $\pi$, let us denote $\mathcal{E}^\pi(\be)=\lim_N
\frac{1}{N}\log Z_N^\pi(\be)$. Hence
$$\lim_N \frac{1}{N}\log
Z_N(\be) = \lim_N \frac{1}{N}\log\dmax_\pi Z_N^\pi(\be).$$
Now
$\log$ being increasing function we can bring the $\max$ out side
$\log$ and the range of $\pi$ being finite we can push the limit
after $\max$ so that
$$\lim_N \frac{1}{N}\log Z_N(\be) =
\dmax_\pi\lim_N \frac{1}{N}\log Z_N^\pi(\be)=\dmax_\pi
\mathcal{E}^\pi(\be).$$
Thus we have the following
\begin{thm}
In the block tree GREM, the limiting free energy $\mathcal{E}(\be)$
exists almost surely and $$\mathcal{E}(\be)=\dmax_\pi
\mathcal{E}^\pi(\be).$$
\end{thm}
\begin{rem}
In the definition of weights, we fixed numbers $a_1,\cdots, a_n$ and
weighted $\xi_{\sigma_{\pi(1)}\sigma_{\pi(2)}
\cdots\sigma_{\pi(i)}}^{\pi(1)\pi(2)\cdots\pi(i)}$ with
$a_{\pi(i)}$. Instead one could fix for each $s$, a sequence of
distinct elements of $I$, a number $a_s$ and then
$\xi_{\sigma_{\pi(1)}\sigma_{\pi(2)}
\cdots\sigma_{\pi(i)}}^{\pi(1)\pi(2)\cdots\pi(i)}$ could be weighted
with $a_{\{\pi(1),\cdots,\pi(i)\}}$. Different driving distributions
for different $s$ can also be considered. Then also the above
theorem remains true.
\end{rem}

\begin{rem}
We now consider the weights ($a_s,\; s\subseteq I$) as mentioned in
the above remark. Using the notation of previous section, consider
BK-GREM with these weights. Consider the GREM associated with
$\pi\in \mathcal{P}_I$ in the BK-GREM and denotes its energy by
$\mathcal{E}(\pi,\be)$.

On the other hand, consider block tree GREM as mentioned in the
above theorem with weights
$\widetilde{a}_{\{\pi(1),\cdots,\pi(i)\}}$, where
$$\widetilde{a}_{\{\pi(1),\cdots,\pi(i)\}}^2=
\dsum_{\substack{s\subseteq{\{\pi(1),\cdots,\pi(i)\}}\\
s\nsubseteq{\{\pi(1),\cdots,\pi(i-1)\}}}}a_s^2.$$

By (\ref{e3.5.21}), we observe that for each $\pi\in \mathcal{P}_I$,
$\widetilde{a}_{\{\pi(1),\cdots,\pi(i)\}}=w(\pi,i)$ for $1\leq i\leq
n$. Moreover, for a fixed $\pi\in \mathcal{P}_I$, in both the
associated GREM model have $2^{k(\pi(i),N)}$ many edges at the
$i$-th level below each node of the $(i-1)$-th level. Hence
$\mathcal{E}^\pi(\be)=\mathcal{E}(\pi,\be)$ for each $\pi\in
\mathcal{P}_I$. Thus the limiting free energy of this block tree
GREM is larger than that of the BK-GREM.
\end{rem}


\cleardoublepage

\chapter{Word GREM with External Field}



In this concluding chapter we discuss a more general version of
random energy models, called Word GREM. This model includes
Derrida's REM and GREM, also the model of Bolthausen and Kistler.
Moreover the model is considered with external field. We apply this
analysis to analyze the free energy of REM with external field.
\smallskip

\goodbreak

\section{Word GREM}

In the previous chapters we have shown that the almost sure
existence of the limiting free energy is assured through the simple
LDP of certain empirical measures. This techniques is quite simple
and neat. In this section, we present a general setup which includes
all the models mentioned above. However, it is not just the
generalization that should be noted. More importantly, we use the
same large deviation technique which allows us to introduce external
field in the model. To our knowledge these models are so far not
discussed with external field except the REM by Derrida in
\cite{D1}. Not only that, as already mentioned in the previous
chapter this method allows consideration of different driving
distributions. This in turn leads to diverse covariance structures
for the Hamiltonian.

\section{The Model}

Let $I=\{\varsigma_1,\varsigma_2,\cdots,\varsigma_n\}$ be a set of
$n$ symbols where $n\geq 1$ is a positive integer. Let $S(I)$ be the
set of all words formed by these $n$ symbols. Let $S$ be a finite
subset of $S(I)$. So a typical word $s \in S$ of length $l$ will
look like $s=\varsigma_{i_1}\varsigma_{i_2}\cdots\varsigma_{i_l}$
where each $\varsigma_{i_j}\in I$. Occasionally we will use the
symbol $s\in S$ as a word as well as a subset of $I$ consisting of
all the symbols in $s$. Since symbols may be repeated in a word, it
is possible that two different words may correspond to the same
subset of $I$. Moreover, without loss of generality we assume that
each symbol appears in at least one word of $S$, that is
$\bigcup\limits_{s\in S} s = I$.

For $N\geq n$, the $N$ particle system has configuration space, as
usual, $\Si_N=\{+1,-1\}^N$ consisting of sequence of length $N$ with
entries +1 and $-1$. For $1\leq i\leq n$, let $k(i,N)\geq 1$ be
integers with $\dsum_{i=1}^n k(i,N)=N$ and $\frac{k(i,N)}{N}\ra
p_i>0$ as $N\ra\infty$. Clearly, $\sum_1^n p_i=1$.

For $\si=\seq{\si_1,\cdots,\si_N}\in\Si_N$, we denote
$\si^1=\seq{\si_i:i\leq k(1,N)},\;\si^2 = \seq{\si_i : k(1,N)+1 \leq
i \leq k(1,N)+k(2,N)}$, etc. Thus $\si$ can also be written as
$\si=\seq{\si^1,\cdots,\si^n}$. For each
$s=\varsigma_{i_1}\varsigma_{i_2}\cdots\varsigma_{i_l}\in S$ and
$\si=\seq{\si^1,\cdots,\si^n}$, we put,
$\si(s)=\seq{\si^{i_1},\si^{i_2},\cdots,\si^{i_l}}$,
$k(s,N)=\dsum_{i=1}^n k(i,N)1_{\{\varsigma_i\in s\}}$.

For each $s\in S$ and $\si\in \Si_N$ we have a random variables
$\xi(s,\si(s))$. These are assumed to be independent random
variables(distributions in general depend on $N$.) To make it more
precise, denote  $\Si_{iN}= \{+1,-1\}^{k(i,N)}$, for $1\leq i\leq
n$. For each
$s=\varsigma_{i_1}\varsigma_{i_2}\cdots\varsigma_{i_l}\in S$ and
$\si(s)=\seq{\si^{i_1},\cdots,\si^{i_l}}\in
\Si_{i_1N}\times\cdots\times\Si_{i_lN}$, we have one random variable
$\xi(s,\si(s))$. All these $\dsum_{s\in S} 2^{k(s,N)}$ random
variables are independent. Let us assume, for $s\in S$, all the
$\xi(s,\si(s))$ have distribution $\la_N^s$ on $\mbb{R}$, that is,
the distribution of $\xi(s,\si(s))$ depends on $s$ but not on
$\si(s)$. Let $f:\mbb{R}^{S}\ra\mbb{R}$ be a continuous function.
For the configuration $\si=\left<\si_1,\si_2,\cdots,\si_N\right>$,
the Hamiltonian of the system is defined as
\begin{equation}\label{e4.2.1}
H_N(\si,h)=Nf(\xi(\si)) +h\dsum_{i=1}^N \si_i,
\end{equation}
where $\xi(\si)=(\xi(s,\si(s)))_{s\in S}$ and $h\geq 0$ is a number
representing the intensity of the external field. The partition
function of the system is
$$Z_N=\dsum_\si e^{-\be H_N(\si,h)},$$ with $\be>0$ being the inverse
temperature. Once again the limiting
free energy is $\dlim_N \frac{1}{N} \log Z_N(\be)$.

\begin{rem}
Observe that if $S=S_1$ consists of only one word
$\varsigma_1\varsigma_2\cdots\varsigma_n$, and if $f(x)=x$ then this
is just the REM. If $S=S_n$ consists of the $n$ words
$\{\varsigma_1, \varsigma_1\varsigma_2,\cdots,\varsigma_1\varsigma_2\cdots\varsigma_n\}$, and
if $f((x_s)_s)=\dsum_{s\in S} a_s x_s$ then this is just GREM. On
the other hand, if $S=S_{2^n}$ consists of all the words
$\varsigma_{i_1}\varsigma_{i_2}\cdots\varsigma_{i_l}$ with out
repetition of symbols then $S$ can clearly be identified as the
collection of non-empty subsets of $S$. If moreover,
$f((x_s)_s)=\dsum_{s\in S} a_s x_s$ then this will lead to the
BK-GREM. Of course, one could also take $S=S_{n!}$ consisting of all
the $n!.n$ many words
$\{\varsigma_{\pi(1)}\varsigma_{\pi(2)}\cdots\varsigma_{\pi(l)}:
1\leq l \leq n\, \&\, \pi \text{ a permutation of
$\{1,2,\cdots,n\}$}\}$.
\end{rem}

Let $\overline{\si}^i$ denote the sum of the $k(i,N)$ many $+1$ and
$-1$ appearing in $\si^i$. In other words, if
$\si=\left<\si_1,\si_2,\cdots,\si_N\right>$ then $\ol{\si}^i$ is the
sum of all $\si_j$ where $j$ satisfies $k(1,N)+\cdots+k(i-1,N)+1\leq j\leq
k(1,N)+\cdots k(i,N)$. Then, note from (\ref{e4.2.1}) that
\begin{equation}\label{e4.2.1a} \frac{H_N(\si,h)}{N}=f(\xi(\si))
+\frac{h}{N}\dsum_{i=1}^n \ol{\si}^i,
\end{equation}
and \begin{equation}\label{e4.2.1b} Z_N(\be,h)=2^N E_\si e^{-N\be
\frac{H_N(\si,h)}{N}},
\end{equation}
where $E_\si$ is the expectation with respect to the uniform
probability on the configuration space.

Under certain assumptions we shall show that the limit $\dlim_N
\frac{1}{N} \log Z_N(\be)$ exists almost surely and is a non-random
quantity. The essential assumptions are the following: Firstly the
distributions of $\xi$ should have exponential decay and secondly
$\frac{k(i,N)}{N}$ converges.

\subsection*{Notations}
We start with some notations which we will use in the rest of the
chapter. A typical points in $\mbb{R}^{S}\times\mbb{R}^n$ will be
denoted by $\left((x_s, s\in S),(y_i, i\leq n)\right)$ or simply as
$(x_{_S},y_{_I})$. In what follows, $\square=\prod_{s\in
S}\tri_s\times\prod_{i=1}^n\itri_i$ is a box in
$\mbb{R}^S\times\mbb{R}^n$ where $\tri_s$ for each $s\in S$ and
$\itri_i$ for $i\leq n$ are open subintervals of $\mbb{R}$.

For $A\subseteq I$, let $\mathcal{S}_A=\{s\in S: s\subseteq A\}$. So
note that $\mathcal{S}_I=S$. We will denote
$Q_{AN}=\prod_{s\in\mathcal{S}_A} q_{sN}$ where
$q_{sN}=\la_N^s(\tri_s)$. If $\mathcal{S}_A$ is empty for some $A$,
we put $Q_{AN}=1$. Also we will denote $Q_{sN}$ with the same
understanding as above considering $s$ as a subset of $I$. Strictly
speaking we should denote $Q_{AN}$ and $q_{sN}$ as $Q_{AN}(\square)$
and $q_{sN}(\tri_s)$ respectively, but for ease of writing we are
not doing so. If
$A=\{\varsigma_{i_1},\varsigma_{i_2},\cdots,\varsigma_{i_m}\}$, then
sometimes we need only the indices $\{i_1, i_2,\cdots,i_m\}$ and we
will denote them by $[A]$. For $A\subseteq I$, we denote
$k(A,N)=\dsum_{i=1}^n k(i,N)1_{\{\varsigma_i\in A\}}=\dsum_{i\in
[A]} k(i,N)$ and $\al_{AN}=\frac{1}{2^{k(A,N)}}\dsum_{\left<\si^i:\;
i\in [A]\right>}\dprod_{i\in
[A]}1_{\itri_i}(\frac{\overline{\si}^i}{N})$. We want to point out
once again that time to time we will consider $s\in S$ as a subset
of $I$. For example, if
$s=\varsigma_{i_1}\varsigma_{i_2}\cdots\varsigma_{i_l}\in S$, we
will use the notation
$\al_{sN}=\frac{1}{2^{k(s,N)}}\sum_{\si(s)}\prod_{j=1}^l1_{\itri_{i_j}}
\left(\frac{\ol{\si}^{^{i_j}}}{N}\right)$.

\section{A large deviation principle}\label{s4.2}

For each $s\in S$, let us consider a probability $\la^s$ on
$\mbb{R}$. If $X$ is distributed like $\la^s$, let us denote
$\Lambda_s(\rho)=\log E e^{\rho X}$ and
$\mathcal{D}_{\Lambda_s}=\{\rho:\,\Lambda_s(\rho)< \infty\}$. Note
that $0\in\mathcal{D}_{\Lambda_s}$, but we want that $0\in
\mathcal{D}_{\Lambda_s}^0$. So from now on we will focus our
attention on those $\la^s$ for which 0 is an interior point in
$\mathcal{D}_{\Lambda_s}$. As $0\in \mathcal{D}^0_{\la^s}$, the mean
$\ol{x}_s=\int xd\la^s(dx)$ exists and is finite quantity for each
$s\in S$. Now if $X_1^s,X_2^s,\cdots$ are i.i.d. random variables
having distribution $\la^s$, we will consider $\la_N^s$ to be the
law of $\frac{1}{N}(X_1^s+X_2^s+\cdots+X_N^s)$. Now by Cramer's
theorem (Theorem \ref{p0.3.6}) the sequence $\{\la_N^s\}$ satisfies
large deviation principle with a good, convex rate function
$\mathcal{I}_s$ given by
$\mathcal{I}_s(x)=\dsup_{\rho\in\mbb{R}}\{\rho x-
\Lambda_s(\rho)\}$. Note that this is a convex, good, non-negative
lower semicontinuous function. Moreover, by property of good rate
function, $\mathcal{I}_s(\ol{x}_s)=0$ for every $s\in S$ so that the
set $\mathcal{I}_s(x)<\al$ is non-empty for every $\al>0$. We also
want to point out that the functions $\mathcal{I}_s$ are increasing
on $[\ol{x}_s,\infty)$ and decreasing on $(-\infty,\ol{x}_s]$.

Once again by Cramer's theorem, the arithmetic averages of i.i.d.
mean zero, $\pm 1$ valued random variables satisfy LDP with rate
function $\mathcal{I}_0$ where $\mathcal{I}_0(y)=\infty$ for
$|y|>1$; $\mathcal{I}_0(\pm 1)=\log2$ and for $-1< y < 1$,
\begin{equation}\label{e4.2.0}
\begin{array}{rl}\mathcal{I}_0(y)&= y\tanh^{-1}y
-\log\cosh(\tanh^{-1}y)\vspace{1ex}\\
&=\frac{1+y}{2}\log(1+y)+\frac{1-y}{2}\log(1-y).
\end{array}\end{equation}

Let us define the map from $\Si_N\ra \mbb{R}^S\times \mbb{R}^n$ as
follows:
$$\si=\left(\si^1,\si^2,\cdots,\si^n\right)\mapsto \left(\left(\xi(s,\si(s)),\,
s\in S\right),\left(\frac{\overline{\si}^i}{N},\, 1\leq i\leq
n\right)\right),$$ where $\ol{\si}^i$ is the sum of the entries of
$\si^i$.

Thus for each $\om$ (sample point of the random variables $\xi$,
which is suppressed so far), this map transports the uniform
probability on $\Si_N$ to $\mbb{R}^S\times\mbb{R}^n$. Denote this
induced random probability by $\mu_N$. Hence from (\ref{e4.2.1b}),
we have,
\begin{equation}\label{e4.2.1c}
\frac{1}{N}\log Z_N(\be,h)=\log2- \frac{1}{N}\log\int_{\mbb{R}^S\times
\mbb{R}^n}e^{-N\be\left(f(x_S)+h\sum\limits_{i=1}^n
y_i\right)}d\mu_N(x_S,y_I).
\end{equation}

\begin{prop}\label{p4.2.1}
If for some $A\subseteq I$, $\sum_{N\geq n} 2^{k(A,N)} Q_{AN}
\al_{AN}<\infty$ then almost surely eventually $\mu_N(\square)=0$.
\end{prop}

\begin{proof} Let $A$ be such that $\sum_{N\geq
n}2^{k(A,N)}Q_{AN}\al_{AN}<\infty$. Then
\begin{align*}
\mu_N(\square)&=\frac{1}{2^N}\dsum_\si \dprod_{s\in S}
1_{\tri_s}\left(\xi(s,\si(s))\right)\dprod_{i\leq
n}1_{\itri_i}\left(\frac{\ol{\si}^i}{N}\right)\\
&\leq\frac{1}{2^N}\dsum_\si \dprod_{s\in \mathcal{S}_A}
1_{\tri_s}\left(\xi(s,\si(s))\right)\dprod_{i\in [A]}
1_{\itri_i}\left(\frac{\ol{\si}^i}{N}\right)\\
&=\frac{1}{2^{k(A,N)}}\dsum_{\si^i : i\in [A]} \dprod_{s\in
\mathcal{S}_A} 1_{\tri_s}\left(\xi(s,\si(s))\right)\dprod_{i\in [A]}
1_{\itri_i}\left(\frac{\ol{\si}^i}{N}\right)
\end{align*}

As a consequence,

\begin{align*} P(\mu_N(\square)>0)&=P\left(\dsum_{\si_i
: i\in [A]} \dprod_{s\in \mathcal{S}_A}
1_{\tri_s}\left(\xi(s,\si(s))\right)\dprod_{i\in
[A]}1_{\itri_i}\left(\frac{\ol{\si}^i}{N}\right)\geq 1\right)\\
&\leq Q_{AN}\dsum_{\si^i : i\in [A]} \dprod_{i\in
[A]}1_{\itri_i}\left(\frac{\ol{\si}^i}{N}\right)\\
&=2^{k(A,N)}Q_{AN}\al_{AN}.
\end{align*}

The hypothesis and Borel-Cantelli lemma completes the proof.
\end{proof}

Let us note that, $E\mu_N(\square)=Q_{IN} \al_{IN}$.
\begin{prop}\label{p4.2.2}
If for all non-empty $A\subseteq I$, $\sum_{N\geq
n}2^{-k(A,N)}Q_{AN}^{-1}\al_{AN}^{-1} < \infty$ then for all
$\ep>o$, almost surely eventually, $$(1-\ep)Q_{IN}\al_{IN}\leq
\mu_N(\square)\leq(1+\ep)Q_{IN}\al_{IN}.$$ That is $$(1-\ep)E\mu_N(\square)\leq
\mu_N(\square)\leq(1+\ep)E\mu_N(\square).$$
\end{prop}
\begin{proof}
Note that
\begin{align*}
&\mbox{Var}(\mu_N(\square))\\
\begin{split}=&\frac{1}{2^{2N}}\dsum_\si\dsum_\tau E\left( \dsty\prod_{s\in S}
1_{\tri_s}\left(\xi(s,\si(s))\right)1_{\tri_s}\left(\xi(s,\tau(s))\right)\right)\dsty\prod_{i\leq
n}1_{\itri_i}\left(\frac{\ol{\si}^i}{N}\right)1_{\itri_i}   \left(\frac{\ol{\tau}^i}{N}\right)\vspace{1ex}\\
&\hspace{65ex}-Q_{IN}^2\al_{IN}^2\end{split}\vspace{2ex}\\
\begin{split}\leq&\frac{1}{2^{2N}}\dsum_{\substack{A\subseteq I\\A\neq
\phi}}\dsum_\si\dsum_{\substack{\tau_i=\si_i,
\forall i\in [A]\\\tau_i\neq \si_i, \forall i\in [A^c]
}}E\left(\dsty\prod_{s\in S}
1_{\tri_s}\left(\xi(s,\si(s))\right)1_{\tri_s}\left(\xi(s,\tau(s))\right)\right)\times\vspace{1ex}\\
&\hspace{55ex}\dsty\prod_{i\leq
n}1_{\itri_i}\left(\frac{\ol{\si}^i}{N}\right)1_{\itri_i}\left(\frac{\ol{\tau}^i}{N}\right)\end{split}\vspace{2ex}\\
\intertext{(since $Q_{IN}^2\al_{IN}^2$ cancels the terms
corresponding to $\si_i \neq \tau_i,\,\forall i\in [I]$)}
=&\frac{1}{2^{2N}}\dsum_{\substack{A\subseteq I\\A\neq
\phi}}\frac{Q_{IN}^2}{Q_{AN}}\dsum_\si\dsum_{\substack{\tau_i=\si_i,
\forall i\in [A]\\\tau_i\neq \si_i, \forall i\in [A^c]
}}\dsty\prod_{i\in
[A]}1_{\itri_i}\left(\frac{\ol{\si}^i}{N}\right)\dsty\prod_{i\in
[A^c]}1_{\itri_i}\left(\frac{\ol{\si}^i}{N}\right)1_{\itri_i}\left(\frac{\ol{\tau}^i}{N}\right)\vspace{2ex}\\
\intertext{(by definition of $Q_{AN}$, $A\subseteq I$)}
=&\dsum_{\substack{A\subseteq I\\A\neq \phi}}
\frac{Q_{IN}^2}{Q_{AN}}\,\frac{1}{2^{k(A,N)}}\,\frac{\al_{IN}^2}{\al_{AN}}.
\end{align*}

Now by Chebyshev's inequality for any $\ep>0$
\[P(|\mu_N(\square)-E\mu_N(\square)|>\ep E\mu_N(\square))<
\frac{1}{\ep^2}\dsum_{A\subseteq I}
\frac{1}{2^{k(A,N)}Q_{AN}\al_{AN}}.\] Once again Borel-Cantelli
lemma and the hypothesis yield that a.s. eventually,
$$(1-\ep)E\mu_N(\square) \leq \mu_N(\square) \leq (1+\ep)
E\mu_N(\square).$$
Hence the proof.
\end{proof}

\begin{thm}\label{t4.3.3}
For a.e. $\omega$, the sequence $\{\mu_N(\omega), N\geq n\}$
satisfies LDP with rate function $\mathcal{J}$ given as follows:

$\mathcal{D}_\mathcal{J}=\left\{(x_{_S},y_{_I}):\; \forall A\subseteq I,
\dsum_{t\in S_A} \mathcal{I}_t(x_t)+\dsum_{i \in [A]}
p_i\mathcal{I}_0\left(\frac{y_i}{p_i}\right) \leq \dsum_{k\in [A]}
p_k\log2 \right\}$

and $$\begin{array}{rll} \mathcal{J}(x_{_S},y_{_I})& = \dsum_{s\in
S} \mathcal{I}_s(x_s) + \dsum_{i \in [I]}
p_i\mathcal{I}_0\left(\frac{y_i}{p_i}\right)& \mbox{if }
(x_{_S},y_{_I})\in \mathcal{D}_\mathcal{J}\\
&=\infty&\mbox{otherwise}
\end{array}$$
\end{thm}

\begin{proof} In what follows, $A$ denotes a non empty subset of $I$.

First of all note that, as $\mathcal{I}_0$ and $\mathcal{I}_s$ for
$s\in S$ are convex, good rate functions, $\mathcal{D}_\mathcal{J}$ is a convex compact set.

Now let $\square=\prod_{s\in S}\tri_s\times\prod_{i=1}^n\itri_i$ be
an open box in $\mbb{R}^S\times\mbb{R}^n$ where $\tri_s$ for each
$s\in S$ and $\itri_i$ for $i\leq n$ are subintervals of $\mbb{R}$
with rational end points.

\noindent{\bf\ul{Step 1}}\quad Suppose that closure of $\square$ is
disjoint with $\mathcal{D}_\mathcal{J}$, that is
$\mathcal{D}_\mathcal{J}\cap \ol{\square}=\phi$. In other words, for
every $(x_S,y_I)\in\ol{\square}$, there exists an $A\subseteq
I$(depending on $(x_S,y_I)$) so that $\dsum_{t\in S_A}
\mathcal{I}_t(x_t)+\dsum_{i \in [A]}
p_i\mathcal{I}_0\left(\frac{y_i}{p_i}\right) > \dsum_{k\in [A]}
p_k\log2$. We shall show, almost surely eventually $\mu_N(\square)=0$.

Note that as $\mathcal{I}_0$ and $\mathcal{I}_s$ are lower
semicontinuous functions for every $s\in S$ and $\ol{\itri}_i$ and
$\ol{\tri}_s$ are compact sets, we can get $(x_S^0,y_I^0)\in
\ol{\square}$ so that $\mathcal{I}(y_i^0)=\mathcal{I}_0(\itri_i)$
for $1\leq i\leq n$ and $\mathcal{I}(x_s^0)=\mathcal{I}_s(\tri_s)$
for every $s\in S$.

For this point $(x_S^0,y_I^0)\in \ol{\square}$ there exists an
$A\subseteq I$ so that $\dsum_{t\in S_A}
\mathcal{I}_t(x_t^0)+\dsum_{i \in [A]}
p_i\mathcal{I}_0\left(\frac{y_i^0}{p_i}\right) > \dsum_{k\in [A]}
p_k\log2$. We will prove that for this $A$ the hypothesis of
Proposition \ref{p4.2.1} is satisfied and hence for this $\square$
almost surely eventually $\mu_N(\square)=0$ leading to
$\lim\limits_{N\ra \infty}\frac{1}{N} \log \mu_N(\square) =
-\infty.$

Since  $\{\la_N^s\}_N$ satisfies LDP with rate function
$\mathcal{I}_s$, we have $$\limsup \frac{1}{N} \log
\la_N^s(\tri_s)\leq -\mathcal{I}_s(\ol{\tri}_s).$$ Let $\ep>0$, to
be chosen later. For all large $N$, $$\frac{1}{N} \log
\la_N^s(\tri_s)<-\mathcal{I}_s(\ol{\tri}_s)+\ep=-\mathcal{I}_s(x_s^0)+\ep,$$
that is $q_{sN}=\la_N^s(\tri_s)<e^{-N(\mathcal{I}_s(x_s^0)-\ep)}$
eventually. And this is true for every $s\in S_A$. So eventually
$$Q_{AN}<e^{-N\dsum_{s\in S_A}(\mathcal{I}_s(x_s^0)-\ep)}.$$

Similarly, the law of $\frac{\overline{\si}^i}{N}$ satisfies LDP
with rate $p_i\mathcal{I}_0(\frac{y_i}{p_i})$ and hence we will have
eventually
$$\al_{AN}<e^{-N\dsum_{i\in
[A]}(p_i\mathcal{I}_0(\frac{y_i^0}{p_i})-\ep)}.$$

Thus $$ 2^{k(A,N)} Q_{AN} \al_{AN}<e^{-N\left[\dsum_{s\in
S_A}(\mathcal{I}_s(x_s^0)-\ep)+\dsum_{i\in
[A]}\left(p_i\mathcal{I}_0(\frac{y_i^0}{p_i})-\ep-\frac{k(i,N)}{N}\log2\right)\right]}.$$
Now as $\frac{k(i,N)}{N}\ra p_i$ and we have strict inequality in
$\dsum_{t\in S_A} \mathcal{I}_t(x_t^0)+\dsum_{i \in [A]}
p_i\mathcal{I}_0\left(\frac{y_i^0}{p_i}\right) > \dsum_{i\in [A]}
p_i\log2$, we can choose an $\ep$ so that $$\sum_{N\geq n}
2^{k(A,N)} Q_{AN} \al_{AN}<\infty.$$ Hence by Proposition
\ref{p4.2.1} we have, almost surely $\mu_N(\square)=0$. It is not
difficult to see now that almost surely $\mu_N$ is eventually
supported on a compact set.

\noindent{\bf\ul{Step 2}}\quad Let us now consider a $\square$ which
has non-empty intersection with $\mathcal{D}_\mathcal{J}$. We show
that for this $\square$, almost surely,

\begin{align}
-\left[\dsum_{s\in S} \mathcal{I}_s(\tri_s)+\dsum_{1\leq i \leq
n}p_i\mathcal{I}_0(\frac{1}{p_i}\itri_i)\right]&\leq \liminf\limits_{N\ra \infty}\frac{1}{N} \log \mu_N(\square)\notag\\
&\leq \limsup\limits_{N\ra \infty}\frac{1}{N} \log \mu_N(\square)\label{e4.2.a2a}\\
&\leq -\left[\dsum_{s\in S} \mathcal{I}_s(\ol{\tri}_s)+\dsum_{1\leq
i \leq
n}p_i\mathcal{I}_0(\frac{1}{p_i}\ol{\itri}_i)\right].\notag\end{align}

Using LDP, we have
$$\liminf_N\frac{1}{N}\log \la_N^s(\tri_s)\geq -
\mathcal{I}_s(\tri_s).$$ Hence for $\ep>0$ eventually,
\begin{equation}\label{e4.2.2a}
\frac{1}{N}\log \la_N^s(\tri_s)>-\mathcal{I}_s(\tri_s)-\ep
\end{equation}
for every $s\in S$. Moreover, eventually,
\begin{equation}\label{e4.2.2b}
\frac{1}{N}\log P\left(\frac{\ol{\si}^i}{N}\in\itri_i\right)> -
p_i\mathcal{I}_0\left(\frac{\itri_i}{p_i}\right)-\ep
\end{equation} for every
$i\leq n$. Hence for every $A\subseteq I$,
\begin{align*}
2^{-K(A,N)}Q_{AN}^{-1}\al_{AN}^{-1}&=e^{-N\left(\dsum_{i\in
[A]}\frac{k(i,N)}{N}\log 2
 +\dsum_{s\in S_A}\frac{1}{N}\log\la_N^s(\tri_s)
 +\dsum_{i\in [A]}\frac{1}{N}\log P(\frac{\ol{\si}^i}{N}\in\itri_i)\right)}\\
&< e^{-N\left[\dsum_{i\in [A]}\frac{k(i,N)}{N}\log 2
 -\dsum_{s\in S_A}(\mathcal{I}_s(\tri_s)+\ep)
 -\dsum_{i\in
 [A]}\left(p_i\mathcal{I}_0\left(\frac{{\itri}_i}{p_i}\right)+\ep\right)\right]},
 \end{align*}
by (\ref{e4.2.2a}) and (\ref{e4.2.2b}).

As $\mathcal{D}_\mathcal{J}$ is a convex set and $\square$ is an
non-empty open set, there exists at least one $(x^0_S,y^0_I)$ in
$\mathcal{D}_\mathcal{J}^-\cap\square$, where
$$ \mathcal{D}_\mathcal{J}^-=\left\{(x_{_S},y_{_I}):\; \forall A\subseteq I,
\dsum_{t\in S_A} \mathcal{I}_t(x_t)+\dsum_{i \in [A]}
p_i\mathcal{I}_0\left(\frac{y_i}{p_i}\right) < \dsum_{k\in [A]}
p_k\log2 \right\}.$$ Being a point in $\mathcal{D}_\mathcal{J}^-$,
for every $A\subseteq I$, we have
$$\sum_{s\in S_A}\mathcal{I}_s(x^0_s)+\sum_{i\leq
n}p_i\mathcal{I}_0\left(\frac{y^0_i}{p_i}\right)<\sum_{i\in
A}p_i\log2.$$ That is $$\sum_{s\in
S_A}\mathcal{I}_s(\tri_s)+\sum_{i\leq
n}p_i\mathcal{I}_0\left(\frac{\itri_i}{p_i}\right)<\sum_{i\in
A}p_i\log2.$$ The above being a strict inequality, we can choose
$\ep$ depending on $\square$ so that for every $A\subseteq I$ the
quantity
$$2^{-K(A,N)}Q_{AN}^{-1}\al_{AN}^{-1}$$
is summable over $N$. Now Proposition \ref{p4.2.2} yields
(\ref{e4.2.a2a}). This completes step 2.

Towards step 3, let $\mathcal{A}$ be the collection of all open
boxes $\square$ with rational corner points satisfying either
$\ol{\square}\cap \mathcal{D}_\mathcal{J}=\emptyset$ or $\square\cap
\mathcal{D}^-_\mathcal{J}\neq \emptyset$. This collection is so rich
that they form a base for the topology of
$\mbb{R}^S\times\mbb{R}^I$. Note that $\mathcal{A}$ being a
countable family,  out side a grand null set, for every $\square$ in
$\mathcal{A}$ conclusions of Step 1 and Step 2 hold. In the next two
steps, we show
\begin{align}
 \mathcal{J}(x_S,y_I)&=\dsup_{\square: (x_S,y_I)\in\square}
\{-\liminf_N\frac{1}{N}\log \mu_N(\square)\}\label{e4.2.a5}\\
&=\dsup_{\square: (x_S,y_I)\in\square}
\{-\limsup_N\frac{1}{N}\log \mu_N(\square)\}.\label{e4.2.b5}
\end{align}

\noindent{\bf\ul{Step 3}}\quad  Let $(x_S,y_I)\notin
\mathcal{D}_\mathcal{J}$. Then $\mathcal{D}_\mathcal{J}$ being a
closed set we can find a $\square\in \mathcal{A}$ containing
$(x_S,y_I)$  so that $\ol{\square}$ does not intersect with
$\mathcal{D}_\mathcal{J}$. By Step 1, $\mu_N(\square)=0$ eventually
so that $\dlim_N \frac{1}{N}\log\mu_N(\square)=\infty$. Also by
definition of $\mathcal{J}$, we have $\mathcal{J}(x_S,y_I)=\infty$.
Hence the above equalities hold when $(x_S,y_I)\notin
\mathcal{D}_\mathcal{J}$.

\noindent{\bf\ul{Step 4}}\quad Now let $(x_S^0,y_I^0)\in
\mathcal{D}_\mathcal{J}$ and
$\mathcal{A}_{(x_S^0,y_I^0)}=\{\square\in\mathcal{A}:\,(x_S^0,y_I^0)\in\square\}$.
Then, as observed in Step 2, for every $\square\in\mathcal{A}_{(x_S^0,y_I^0)}$ we have
eventually
\begin{align}
-\left[\dsum_{s\in S} \mathcal{I}_s(\tri_s)+\dsum_{1\leq i \leq
n}p_i\mathcal{I}_0(\frac{1}{p_i}\itri_i)\right]&\leq \liminf\limits_{N\ra \infty}\frac{1}{N} \log \mu_N(\square)\notag\\
&\leq \limsup\limits_{N\ra \infty}\frac{1}{N} \log \mu_N(\square)\label{e4.2.2}\\
 &\leq
-\left[\dsum_{s\in S} \mathcal{I}_s(\ol{\tri}_s)+\dsum_{1\leq i \leq
n}p_i\mathcal{I}_0(\frac{1}{p_i}\ol{\itri}_i)\right].\notag\end{align}

From the first part of the
above inequality we have,
$$\liminf_N\frac{1}{N}\log \mu_N(\square) \geq -\mathcal{J}(\square).$$
And hence
\begin{equation}\label{e4.2.5}
\sup_{\square\in
\mathcal{A}_{(x_S^0,y_I^0)}}\{-\liminf_N\frac{1}{N}\log
\mu_N(\square)\}\leq \sup_{\square\in \mathcal{A}_{(x_S^0,y_I^0)}}
\mathcal{J}(\square)\;\leq\; \mathcal{J}(x_S^0,y_I^0).
\end{equation}

On the other hand, for every $\square \in
\mathcal{A}_{(x_S^0,y_I^0)}$ using the right side inequality of
(\ref{e4.2.2}), we have
$$\limsup_N\frac{1}{N}\log \mu_N(\square)\leq -\mathcal{J}(\ol{\square}).$$

Let $\mathcal{A}_{(x_S^0,y_I^0)}'=\{\square_k \in \mathcal{A}: k\geq
1\}$ be a subclass of $\mathcal{A}$ so that
$\ol{\square}_{k+1}\subset\square_k$ for every $k$ and
$\cap_k\square_k=\{(x_S^0,y_I^0)\}.$ Then
\begin{align}\label{e4.2.6}
\sup_{\square\in
\mathcal{A}_{(x_S^0,y_I^0)}}\{-\limsup_N\frac{1}{N}\log
\mu_N(\square)\}&\geq \sup_{\square\in \mathcal{A}_{(x_S^0,y_I^0)}}
\mathcal{J}(\ol{\square})\notag\\
&\geq \sup_{\square\in \mathcal{A}_{(x_S^0,y_I^0)}'}
\mathcal{J}(\ol{\square})\notag\\
&=\lim_k \mathcal{J}(\ol{\square}_k)\notag\\ &=
\mathcal{J}(x_S^0,y_I^0).
\end{align}
The last equality follows as $\mathcal{J}$ is a good lower
semicontinuous function (see Proposition \ref{p0.3.2}).

Thus Step 3 and Step 4 complete the proof of (\ref{e4.2.a5}) and (\ref{e4.2.b5}).

Now proof of Theorem \ref{t4.3.3} is completed by appealing to
Proposition \ref{p0.3.3} and observing that $\{\mu_N\}$ is
eventually supported on a compact set.
\end{proof}

\begin{rem}
A closer look at the above proof shows that the convexity of the
rate functions $I_s$ is not an essential condition. Good rate
functions on real line whose graph look like `U` will suffice. That
is we could take those non-negative functions $\mathcal{I}$ for
which we will get at most two points $\ul{x}\leq \ol{x}$ so that
$\mathcal{I}$ is zero on $[\ul{x}, \ol{x}]$; strictly decreasing on
$(-\infty, \ul{x}]\cap \{\mathcal{I}\in(0,\infty)\}$ and strictly
increasing on $[\ol{x}, \infty)\cap \{\mathcal{I}\in(0,\infty)\}$.
With these conditions we can find $(x^0_S,y^0_I)$ as stated in Step
1 and Step 2 for the proof to go through keeping all other Steps as
it is.

For each $s$, we started with a distribution $\la^s$ and appealed to
Cramer's theorem (see the first paragraph of this section). Instead,
we could start with $(\la_N^s,\, N\geq n)$ and assume that
$\{\la_N^s\}$ satisfies LDP with rate function $\mathcal{I}_s$
having the properties mentioned in the above paragraph.
\end{rem}

Since almost every sequence of probabilities $\mu_N$ is eventually
supported on a compact set, we could use Varadhan's lemma with any
continuous function. This will lead to the following (see equation
(\ref{e4.2.1c})):
\begin{thm}
In the word GREM, almost surely $$\lim_{N} \frac{1}{N}\log Z_N(\be,
h)= \log2 -\inf_{\mathcal{D}_\mathcal{J}} \left\{\be f(x_S)+ \be
h\sum_{i=1}^ny_i + \mathcal{J}(x_S,y_I) \right\}.$$
\end{thm}

Though we have taken any real valued continuous function $f$ on
$\mathbb{R}^S$, it is customary to consider $f(x_S)=\dsum_{s\in S}
a_sx_S$ so that the Hamiltonian becomes $H_N(\si)=N\dsum_{s\in S}
a_s \xi(s, \si(s)) + h \ol{\si}$ where $a_s, \; s\in S$ are
non-negative weights and $h>0$ is the strength of the external
field. It is also customary to consider Gaussian driving
distribution. In this setup if $\la^s$ is standard normal, the above
theorem will be applicable and will reduce to the following
\begin{cor}\label{c4.3.5}
In the Gaussian word GREM with external field, almost surely, the limiting free energy is\vspace{1ex}
\begin{center}
$\log2 - \dinf_{\mathcal{D}_\mathcal{J}} \left\{\dsum_{s\in
S}\left(\frac{1}{2}x_s^2+\be a_sx_s\right) + \dsum_{i=1}^n
\left(\frac{p_i+y_i}{2}\log\frac{p_i+y_i}{p_i}+\frac{p_i-y_i}{2}\log\frac{p_i-y_i}{p_i}+
\be h y_i\right)\right\},$\vspace{1ex}
\end{center}
\noindent where $\mathcal{D}_\mathcal{J}$ is the set consisting of
$(x_{_S},y_{_I})\in\mbb{R}^S\times\mbb{R}^n$ such that $\forall
A\subseteq I$,\vspace{1ex}
\begin{center}
$\dsum_{s\in S_A} \frac{1}{2}x_s^2+\dsum_{i \in [A]}
\left(\frac{p_i+y_i}{2}\log\frac{p_i+y_i}{p_i}+
\frac{p_i-y_i}{2}\log\frac{p_i-y_i}{p_i}\right) \leq \dsum_{k\in [A]} p_k\log2.$
\end{center}
\end{cor}

Hence the use of large deviation techniques not only ensures the
almost sure existence of the limiting free energy, the calculation
of free energy of the system is then reduced to that of an
optimization problem. Of course, it is not always possible to solve
this optimization problem to arrive at a closed form expression when
the external field is present or different driving distributions are
considered for different $s\in S$. Even for $n=2$ with Gaussian
driving distribution it is difficult to obtain a closed form
expression. The only case where we will get some 'closed' form
expression is the case for $n=1$ - that is the case of REM with
external fields. This we will consider in the next section. But with
no external field the situation is not that worse. In some of the
cases, the method of calculation of the infimum will just reduce to
what we did in section \ref{s3.5}. In that case, though the model,
to start with, was not an $n$ level tree GREM, it reduced (as far as
the free energy is concerned) to an $n$ level tree GREM with
appropriate weights [see \S \ref{s3.5}]. It is quite conceivable
that the present complicated model may always be equivalent to a
tree GREM. We do not think so.

\goodbreak

\section{REM with external field}
As mentioned in the last section, there is no general technique of
obtaining a formula for the word GREM as well as tree GREM energy
with external field. We will discuss here the simple REM with
external field. Let us consider the word GREM where $S$ consists of
only one word, that is, $S$ consists of the word
$\varsigma_1\varsigma_2\cdots\varsigma_n$ where
$I=\{\varsigma_1,\varsigma_2,\cdots,\varsigma_n\}$ (see beginning of
the previous section). In such a case the word GREM reduces exactly
to the usual REM with external field. Thus the Hamiltonian is
$$H_N(\si)=aN\xi_\si +h\sum_{i=1}^N \si_i$$
where $\xi_\si$ are i.i.d. random variables (for each fixed $N$) and
$h, a$ are positive constants. Moreover for the Gaussian REM,
$\xi_\si$ are $\mathcal{N}(0,\frac{1}{N})$. Then by Corollary
\ref{c4.3.5}, the limiting free energy for the REM with external
field exists almost surely and is given by
$$\begin{array}{ll}\mathcal{E}(\be,h)&=
\log2-\dinf_{\mathcal{D}_\mathcal{J}}\left\{\frac{x^2}{2}+\frac{1+y}{2}
\log(1+y)+\frac{1-y}{2}\log(1-y)+\be(ax+hy)\right\}\\
&=\log2-\dinf_{\mathcal{D}_\mathcal{J}}\left\{\frac{x^2}{2}+\mathcal{I}_0(y)+\be(ax+hy)\right\},\end{array}$$
where $\mathcal{I}_0$ is given by (\ref{e4.2.0}) and
$$\begin{array}{ll}\mathcal{D}_\mathcal{J}&=\{(x,y):
\frac{x^2}{2}+\frac{1+y}{2}\log(1+y)+\frac{1-y}{2}\log(1-y)\leq
\log2\}\\
&=\{(x,y): \frac{x^2}{2}+\mathcal{I}_0(y)\leq \log2\}.\end{array}$$

In other words, $$\mathcal{E}(\be,h)=
\log2-\dinf_{\mathcal{D}_\mathcal{J}^+}f(x,y),$$ where
$f(x,y)=\left\{\frac{x^2}{2}+\mathcal{I}_0(y)-\be(ax+hy)\right\}$
and $\mathcal{D}_\mathcal{J}^+$ equals all points of
$\mathcal{D}_\mathcal{J}$ with both coordinates non-negative.

To calculate the above infimum, first fix $\be, h$ and $y$ with
$0\leq y\leq 1$. Then the range of $x$ is $0\leq x\leq
\sqrt{2[\log2-\mathcal{I}_0(y)]}$. It is easy to see that if
$\mathcal{I}_0(y)\leq \log2-\frac{1}{2}\be^2a^2$ then the $\dinf_x
f(x,y)$ is attained for $x=\be a$ and if $\mathcal{I}_0(y)>
\log2-\frac{1}{2}\be^2a^2$ then the infimum is attained for
$x=\sqrt{2[\log2-\mathcal{I}_0(y)]}$. Since $\mathcal{I}_0$ is a
non-negative function, the set $\{\mathcal{I}_0(y)\leq
\log2-\frac{1}{2}\be^2a^2\}$ will be non-empty only when
$\be\leq\frac{1}{a}\sqrt{2\log2}$. For
$\be>\frac{1}{a}\sqrt{2\log2}$, we always have $\mathcal{I}_0(y)>
\log2-\frac{1}{2}\be^2a^2$ so that the infimum is attained for
$x=\sqrt{2[\log2-\mathcal{I}_0(y)]}$. Substituting these values of
$x$ in $f(x,y)$ we obtain the following expression for the infimum
of $f(x,y)$ over $x$. First we need a notation. For
$\be\leq\frac{1}{a}\sqrt{2\log2}$, let $c_\be$ be the solution of
\begin{equation}\label{e4.4.1}
\mathcal{I}_0(c_\be)=\log2-\frac{1}{2}\be^2a^2.
\end{equation}
Then
\begin{equation}\label{e4.4.2}
\varphi(y)=\dinf_{0\leq x\leq\sqrt{2[\log2-\mathcal{I}_0(y)]}}
f(x,y)=
\begin{cases} g_1(y) &\text{if
$\be\leq\frac{1}{a}\sqrt{2\log2}$} \quad\text{and $y\leq
c_\be$},\\
g_2(y) &\text{if $\be\leq\frac{1}{a}\sqrt{2\log2}$} \quad\text{and
$y> c_\be$},\\
g_2(y) &\text{if $\be>\frac{1}{a}\sqrt{2\log2}$},
\end{cases}
\end{equation}
where $$g_1(y)=-\frac{1}{2}\be^2a^2+\mathcal{I}_0(y)-\be h y$$ and
$$g_2(y)= \log2-\be a\sqrt{2[\log2-\mathcal{I}_0(y)]}-\be h y.$$

Since $$g_1'(y)=\tanh^{-1}(y)-\be h,$$ we have $g_1'(0)=-\be h$ and
$$g_1'(y)\lesseqqgtr 0 \Leftrightarrow y\lesseqqgtr \tanh(\be h).$$

On the other hand, as
$$g_2'(y)=\frac{\be a\tanh^{-1}(y)}{\sqrt{2[\log2-\mathcal{I}_0(y)]}}-\be
h,$$ we have $g_2'(0)=-\be h$ and by (\ref{e4.4.1}),
$g_2'(c_\be)=-\be h+\tanh^{-1}(c_\be)$. Thus $g_2'(c_\be)\leq 0$ iff
$c_\be\leq \tanh(\be h)$. Moreover, $$g_2'(y)\lesseqqgtr
0\Leftrightarrow
\frac{a\tanh^{-1}(y)}{\sqrt{2[\log2-\mathcal{I}_0(y)]}}\lesseqqgtr
h.$$

Let $y_0$ be the non-negative solution of
\begin{equation}\label{e4.4.3}
\frac{a\tanh^{-1}(y)}{\sqrt{2[\log2-\mathcal{I}_0(y)]}}= h.
\end{equation}
Such a solution always exists since $\log2-\mathcal{I}_0(y)\ra 0$ as
$y\ra 1$.

Since $\mathcal{I}_0$ is a strictly increasing function of $[0,1]$,
from equations (\ref{e4.4.1}) and (\ref{e4.4.3}), we note that
$$\tanh(\be h)\lesseqqgtr y_0 \Leftrightarrow  y_0\lesseqqgtr c_\be.$$

Now if $\be\leq \frac{1}{a}\sqrt{2\log2}$ and $y_0\leq c_\be$ then
$\tanh(\be h)\leq y_0\leq c_\be$ and the function $\varphi$ in
(\ref{e4.4.2}) is decreasing up to $y=\tanh(\be h)$ and then
increasing. In such case, the $\dinf_{0\leq y\leq1}\varphi(y)$ will
occur at $y=\tanh(\be h)$ so that $$\dinf_{0\leq
y\leq1}\varphi(y)=-\frac{1}{2}\be^2a^2-\log\cosh(\be h).$$

On the other hand, if $\be\leq \frac{1}{a}\sqrt{2\log2}$ and $y_0>
c_\be$ then $c_\be<y_0<\tanh(\be h)$ and the function $\varphi$ is
decreasing up to $y=y_0$ and then increasing. In such case, the
$\dinf_{0\leq y\leq1}\varphi(y)$ will occur at $y=y_0$ so that
$$\dinf_{0\leq y\leq1}\varphi(y)=\log2-\be a x_0-\be h y_0,$$ where
$x_0=\sqrt{2[\log2-\mathcal{I}_0(y)]}=\frac{a\tanh^{-1}y_0}{h}$.

Finally, if $\be> \frac{1}{a}\sqrt{2\log2}$ then the function
$\varphi$ is decreasing up to $y=y_0$ and then increasing. Hence in
this case, the $\dinf_{0\leq y\leq1}\varphi(y)$ will occur at
$y=y_0$ so that
$$\dinf_{0\leq y\leq1}\varphi(y)=\log2-\be a x_0-\be h y_0,$$ where
$x_0=\sqrt{2[\log2-\mathcal{I}_0(y)]}=\frac{a\tanh^{-1}y_0}{h}$.

We can summarize the above discussion in the following:
\begin{thm}\label{t4.4.1}
In the Gaussian REM with external field, the limiting free energy
exists almost surely and given by
$$\mathcal{E}(\be,h)=
\begin{cases}
\log2+\frac{\be^2 a^2}{2}+\log\cosh(\be h)&\text{if $\be\leq
\frac{1}{a}\sqrt{2\log2}$ and $y_0\leq c_\be$}\\
\be\left( ax_0+h y_0\right)&\text{otherwise},
\end{cases}$$
where $y_0$ be the non-negative solution of
$\frac{a\tanh^{-1}(y)}{\sqrt{2[\log2-\mathcal{I}_0(y)]}}= h$,
$\mathcal{I}_0(y)= y\tanh^{-1}y -\log\cosh(\tanh^{-1}y)$, $c_\be$ is
the solution of $\mathcal{I}_0(c_\be)=\log2-\frac{1}{2}\be^2a^2$ and
$x_0=\frac{a\tanh^{-1}y_0}{h}$.
\end{thm}
Note that the case `otherwise' in the theorem above consists of
$\be\leq \frac{1}{a}\sqrt{2\log2}$ and $y_0> c_\be$ or if $\be>
\frac{1}{a}\sqrt{2\log2}$.

Theorem \ref{t4.4.1} provides yet another justification for the
phase diagram (FIG. 3) in \cite{D1} of Derrida.

\cleardoublepage

\bibliographystyle{amsplain}
\providecommand{\bysame}{\leavevmode\hbox to3em{\hrulefill}\thinspace}
\providecommand{\MR}{\relax\ifhmode\unskip\space\fi MR }
\providecommand{\MRhref}[2]{%
  \href{http://www.ams.org/mathscinet-getitem?mr=#1}{#2}
}
\providecommand{\href}[2]{#2}

\end{document}